\documentclass[10pt]{amsart}
\usepackage[margin=0.8in]{geometry}
\usepackage{amsfonts}
\usepackage{amsmath}
\usepackage{amssymb}
\usepackage{amsthm}
\usepackage{mathtools}
\usepackage{slashed}

\usepackage[
backend=biber,
style=alphabetic,
]{biblatex}
\usepackage{color}
\usepackage{float}
\usepackage{graphicx} % Required for inserting images
\usepackage{hyperref}
\usepackage{tikz}

\title{Local and Global Results on Three Dimensional Rarefaction Waves in Spherical Symmetry}
\author{Ruotong Zhang}
\date{}

\newtheorem{theorem}{Theorem}[section]
\newtheorem{corollary}[theorem]{Corollary}

\newtheorem{proposition}[theorem]{Proposition}

\theoremstyle{definition}
\newtheorem{definition}[theorem]{Definition}

\theoremstyle{remark}
\newtheorem{remark}[theorem]{Remark}

\numberwithin{equation}{section}

\begin{document}

\maketitle

\tableofcontents

\section{Introduction}

\subsection{Review of Shocks and Rarefactions}

In the study of compressible Euler equations, which mainly concerns gas dynamics, there are two important classes of phenomena: shock and rarefaction. Shock is a type of strong discontinuity that can emerge from regular solutions and can persist in further evolution. It corresponds to the concentration and intersection of characteristics. Rarefaction, on the other hand, is a type of resolving behavior that can smooth out initial singularities. It corresponds to the dispersion of characteristics. These two phenomena were first rigorously studied by Riemann in the pioneering work \cite{Riemann1860} in 1860. He studied the solution to the 1-D compressible Euler equations with the initial data being two constant states separated by a discontinuity. He provided explicit expressions for the solutions. For generic initial data, four patterns may emerge, which are combinations of shocks and (centered) rarefactions. This problem is now called the Riemann problem, or the shock tube problem. There are shock tube experiments that verify these theoretical predictions. Since they are observable, they should be stable. One can check the textbooks \cite{Courant-Friedrichs-Book-Supersonic-Shock} \cite{Smoller-Book} or the original paper for the calculations. Later, Lax generalized this work to hyperbolic conservation laws in \cite{Lax-Conservation-Law}. The theory in 1-D then gradually developed and is quite fruitful, with the help of the bounded variable (BV) function space that suits the problem very well. There is a famous encyclopedic book \cite{Dafermos-Book-Conservation-Law} by Dafermos that collects many results in the field of hyperbolic conservation laws.

In the multidimensional case, the theory is less complete than in 1-D. The BV space is no longer suitable and people usually use $L^2$ energy estimates.

For the story of shock, there are several main topics: how does the singularity form from regular solutions; how does the singularity (blowup of derivatives) really become the shock, which is a discontinuity across a surface; how does the shock evolve and stay stable.

The answer to the first question starts from Sideris \cite{Sideris-Singularity}, where the stable blowup for classical solutions to the compressible Euler equations is shown via a contradiction argument, which does not provide detailed information on the structure of the singularity. Later in the series of work \cite{Alinhac-Shock-Radial1} \cite{Alinhac-Shock-Radial2} \cite{Alinhac-Shock-Radial3}, Alinhac proved the formation of singularities of 2-D compressible Euler with radial symmetry, and he provided precise estimates up to the blowup time. Then in \cite{Alinhac-Shock1} \cite{Alinhac-Shock2}, Alinhac extended the result to a class of quasilinear wave equations in 2-D without symmetry, and showed that the blowup is caused by the collapse of characteristic hypersurface folations. A major breakthrough in the formation of shock type singularities is due to Christodoulou. In the monograph \cite{Christodoulou-Monograph-Shock}, Christodoulou showed the formation of singularities for irrotational and isentropic fluids satisfying the relativistic compressible Euler equations in 3-D without symmetry. The blowup is caused by the collapse of characteristic hypersurface foliations, and he also describes the geometry of the boundary of the maximal development of the data. There are many works inspired by the framework used in this monograph, where many parts of it come from the work on the nonlinear stability of Minkowski space \cite{Christodoulou-Klainerman-Stability-Minkowski} by Christodoulou and Klainerman.

The second question, which is called the shock development problem, is more tricky to answer. It requires one to go beyond singularity and obtain again a piecewise regular solution containing a discontinuity surface. One breakthrough is made by Christodoulou in \cite{Christodoulou-Shock-Development}. He is able to construct the shock surface from regular solutions in the restricted scheme, where one forces the entropy and vorticity to have zero jump across the shock.

The third question, which is called the shock front problem, was first addressed by Majda in \cite{Majda-Shock-Evolve1} \cite{Majda-Shock-Evolve2}. The stability of shock waves are proved in this work, which means if the initial data contains a shock, this structure will be preserved in the evolution at least for short time. There are many works extending this result in many directions later.

For the story of (centered) rarefaction, there are much less results. The first major results on the construction of multidimensional rarefaction waves are the work \cite{Alinhac-Existence} \cite{Alinhac-Uniqueness} by Alinhac. The local existence and uniqueness of the rarefaction wave solutions are proved by a well-designed Nash-Moser type scheme, but the estimates obtained are in weighted space, which degenerate near the rarefaction fronts, so we don't obtain full control of the solutions. The scheme in this work was then used in many different settings. Recently Luo and Yu proved the local existence and uniqueness of multidimensional rarefaction waves in \cite{Luo-Yu-1} \cite{Luo-Yu-2}, and obtained full control of the norms in the rarefaction region. The geometric framework is inspired by \cite{Christodoulou-Monograph-Shock}, although they are quite different in technical parts. In a sense, rarefaction waves resolving the initial singularity can be viewed as an inversion of the formation of shocks.

Our work is a follow-up of \cite{Luo-Yu-1} \cite{Luo-Yu-2} in the spherical symmetric case. We deal with nontrivial curvature of the initial singular surface, which is not present in their work. However, we avoided lots of technical difficulties due to spherical symmetry. The new thing in this work is that we are able to push the same framework to prove global in time results for data that are good enough. The results obtained are similar to the results in \cite{Wang-Global-Exterior} by Wang. The estimates obtained in this work can be used as a rough reference for the transversal parts of the full problem without symmetry.

\subsection{Problem Setting and Main Results}

The flow considered in this article is three-dimensional homentropic continuous flow of perfect gas with spherical symmetry, governed by Euler equations and the gamma law. The equations are as follows:
$$
\begin{cases}
\partial_t \rho + \mathbf{v} \cdot \nabla \rho = - \rho \nabla \cdot \mathbf{v} \\
\partial_t \mathbf{v} + \mathbf{v} \cdot \nabla \mathbf{v} = - \frac{\nabla p}{\rho} \\
\frac{p}{\rho^\gamma} = K = Const.
\end{cases}
$$ \\
where $\mathbf{v} = v \mathbf{e_r}$ is radial, and all the quantities involved rely only on $r$, the radial direction. \\
Under this assumption, pressure $p$ is now a function of $\rho$, $p = K \rho^\gamma$. \\
The speed of sound $c$ also has an explicit formula in $\rho$, which is $c^2 = (\frac{d p}{d \rho})_{\mathrm{Const\ entropy}} = K \gamma \rho^{\gamma - 1}$.

So, there are essentially two variables involved: the radial velocity $v$ (positive means outgoing) and the density $\rho$.

The specific solution we would like to find is the rarefaction wave solution with a given exterior background flow, which can be seen as the analogue of the Riemann problem in spherical symmetry.

We write down the rough version of the two main results here.

\begin{theorem}[Local Existence of Rarefaction Wave Solutions]
Given a smooth spherical symmetric exterior solution $(v, c)$ to the compressible Euler equation, one can construct a piece of (centered) rarefaction wave solution on a fan-shaped region originating from an initial sphere. The solution approaches the 1-D rarefaction solution when we approach the initial time.
\end{theorem}

\begin{theorem}[Global Existence of Rarefaction Wave Solutions]
If the exterior solution $(v, c)$ is sufficiently close to the constant state and decays to the constant state at proper rate, the rarefaction wave solutions we construct can be extended globally in time. The characteristics separate from each other of size $\ln t$.
\end{theorem}

\begin{figure}[H]
\centering
\scalebox{0.6}{\begin{tikzpicture}

% coordinate axis
\coordinate (origin) at (0,0);
\coordinate (r-axis) at (8,0);
\coordinate (t-axis) at (0,6);
\node at (r-axis) [anchor=north]{r};
\node at (t-axis) [anchor=east]{t};
\draw [thick,->] (origin) -- (r-axis);
\draw [thick,->] (origin) -- (t-axis);
\draw [dotted,semithick] (0,5) -- (8,5);
\node at (0,5) [anchor=east]{$t^*$};
% C_0
\draw [semithick] (3,0) node [anchor=north]{$r_0$} .. controls +(1,2) and +(-1,-2) .. +(3,5) node [anchor=south]{$C_0$};
% Rarefaction
\foreach \X in {1,...,5}
{\draw [dotted,semithick] (3,0) .. controls +(1-0.3*\X,2) and +(-1+0.15*\X,-2) .. +(3-0.7*\X,5);};
% Rarefaction End
\draw [semithick] (3,0) .. controls +(1-0.3*6,2) and +(-1+0.15*6,-2) .. +(3-0.7*6,5) node [anchor=south]{$C_{u^*}$};;
% Exterior
\foreach \X in {1,...,2}
{\draw [dotted,semithick] (3+0.7*\X,0) .. controls +(1-0.1*\X,2-0.2*\X) and +(-1+0.1*\X,-2+0.2*\X) .. +(3,5);};
\foreach \X in {3,...,5}
{\draw [dotted,semithick] (3+0.7*\X,0) .. controls +(1-0.1*\X,2-0.2*\X) and +(-1+0.1*\X,-2+0.2*\X) .. +(3+0.7*2-0.7*\X,5+1.1*2-1.1*\X);};
% Text
\node at (3.5,4) {Rarefaction};
\node at (6,2) {Exterior};

% Coordinates
\coordinate (origin) at (10,0);
\coordinate (r-axis) at (16,0);
\coordinate (t-axis) at (10,6);
\coordinate (singular) at (11,1);
\node at (r-axis) [anchor=north]{r};
\node at (t-axis) [anchor=east]{t};
\draw [thick,->] (origin) -- (r-axis);
\draw [thick,->] (origin) -- (t-axis);
\coordinate (rS) at (11,0);
\coordinate (tS) at (10,1);
\node at (rS) [anchor=north]{1};
\node at (tS) [anchor=east]{1};
\draw [dotted,semithick] (rS) -- (singular);
\draw [dotted,semithick] (tS) -- +(6,0);
% Figure
\draw [semithick] (singular) .. controls +(1,1.2) and +(-1,-1) .. +(4.5,4.5) node [anchor=south]{$C_0$};
\draw [dashed,semithick] (singular) .. controls +(-0.2,1.2) and +(-1+0.12,-1) .. +(2.7,4.5) node [anchor=south]{$C_{u^*}$};
\foreach \x in {1,...,5}
{\draw [dotted,semithick] (singular) .. controls +(1-0.2*\x,1.2) and +(-1+0.02^\x,-1) .. +(4.5-0.3*\x,4.5);};
\fill [gray,opacity=0.5] (singular) .. controls +(1,1.2) and +(-1,-1) .. +(4.5,4.5) -- (16,5.5) -- (16,1) -- cycle;
% Text
\node at (13.5,4.5) {$\kappa \sim \ln t$};
\node at (14,2) {Close to Constant};

\end{tikzpicture}}
\caption{Local and Global}
\label{fig:local-global}
\end{figure}
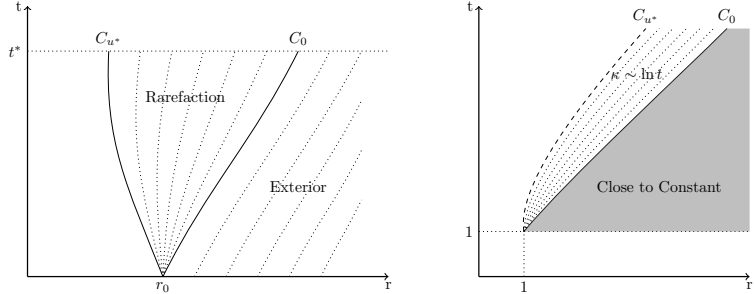

See figure \ref{fig:local-global} for a depiction of the two results. The rigorous versions can be found in \ref{thm:local} and \ref{thm:global}.

We list the steps to prove these results. They correspond to the subsections of the proof.

The proof of the local result can be divided into the following steps.

\begin{itemize}
    \item Determine the data (transversal derivatives) of the rarefaction wave solution on the boundary $C_0$. This is done by solving systems of ODEs along $C_0$. The choice of foliation is also encoded in this process.

    \item Use the transversal derivatives to construct non-singular approximate Cauchy data on time slice $\{t = \delta\}$ by the Taylor polynomials, $\delta$ small. Solve from $\{t = \delta\}$ to obtain approximate solutions.

    \item Prove uniform priori energy estimates in $\delta$ to obtain uniform life span and norm bounds of the approximate solutions.

    \item Send $\delta \to 0$ and perform a convergence argument to obtain the exact rarefaction wave solution.
\end{itemize}

The proof of the global result can be divided into the following steps.

\begin{itemize}
    \item Solve data on $C_0$ explicitly for the constant state exterior solution to obtain the expected decay rate for all the quantities.

    \item Prove uniform priori energy estimates in $T$ (large time) for a class of data on $C_0$ with loosened decay rate requirements. This leads to global existence.

    \item Show that the requirements of data on $C_0$ are indeed verified for exterior solutions that are close to the constant state and have the proper decay rate.
\end{itemize}

\subsection{Review of the Riemann Problem, Riemann invariants} \label{subsec:1D}

Now we review the classical 1D Riemann problem, which will give us some insights of the 3D spherical symmetric problem. We will focus on the 1-family rarefaction wave case.

Consider the case of piston generated rarefaction wave, a tight-fitting piston is suddenly withdrawn from a stationary gas, contained in a uniform tube, at the steady speed $V_p$. Suppose that the piston is located at $x=0$ at $t=0$ , and moves in the $-x$ or left direction.

The equations in 1-dimension are:
\begin{equation*}
\begin{cases}
    \partial_t \rho + v \partial_x \rho = - \rho \partial_x v \\
    \partial_t v + v \partial_x v = - \frac{\partial_x p}{\rho} \\
    \frac{p}{\rho^\gamma} = s = Const.
\end{cases}
\end{equation*}

The initial data are $(v_0 = 0, \rho_1, p_1)$, and the boundary condition on the piston is $v = - V_p$ at $x = - V_p t$.

\begin{figure}[H]
\centering
\scalebox{0.6}{\begin{tikzpicture}

% Coordinates
\coordinate (origin) at (6,0);
\coordinate (x-axis-L) at (0,0);
\coordinate (x-axis-R) at (10,0);
\coordinate (t-axis-T) at (0,5);
\coordinate (t-axis-B) at (0,4);
\node at (x-axis-R) [anchor=north]{x};
\node at (t-axis-T) [anchor=east]{t};
\node at (origin) [anchor=north]{0};
\draw [thick,->] (x-axis-L) -- (x-axis-R);
\draw [thick,->] (t-axis-B) -- (t-axis-T);
\coordinate (RFR) at (9,4);
\coordinate (RFL) at (1.5,4);
% Rarefaction Boundary
\draw [semithick] (origin) -- (RFL) node [anchor=south]{Piston};
\draw [dashed,semithick] (origin) -- (RFR);
% Null Foliation
\foreach \x in {1,...,14}
{\draw [dotted,semithick] (origin) -- +(3-0.5*\x,4);};
\foreach \x in {1,...,2}
{\draw [dotted,semithick] (origin)++(0.5*\x,0) -- +(3,4);};
\foreach \x in {3,...,7}
{\draw [dotted,semithick] (origin)++(0.5*\x,0) -- +(4-0.5*\x,5.33-0.67*\x);};
\fill [gray,opacity=0.5] (origin) -- (RFR) -- ++(1,0) -- ++(0,-4) -- cycle;
% Text
\node at (6,3) {Rarefaction};
\node at (9,2) {Constant};

\end{tikzpicture}}
\caption{1-D Rarefaction}
\label{fig:1-D}
\end{figure}
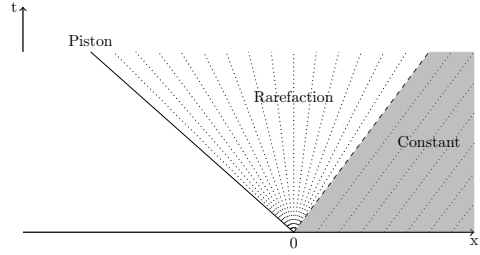

The shape of the solution should be a region moving with the piston with data $(- V_p, \rho_2, p_2)$ and an undisturbed stationary region with data $(0, \rho_1, p_1)$ connected by a fan-shaped rarefaction wave region with data $(v(\frac{x}{t}), \rho(\frac{x}{t}), p(\frac{x}{t}))$. See figure \ref{fig:1-D}.

In the fan-shaped region, the velocity $v$ and the speed of sound $c$ are linear functions of $\frac{x}{t}$.
\begin{equation*}
    v = \frac{2}{\gamma + 1} (\frac{x}{t} - c_1), \qquad c = \frac{\gamma - 1}{\gamma + 1} (\frac{x}{t} - c_1) + c_1
\end{equation*}

The quantities on the two edges are continuous, so this gives the rarefaction wave solution. If $V_p$ is too large, then $c$ will drop to zero and there will be a vacuum region on the left side.

The Riemann invariants are a transformation of the system in order to make it easier to solve. For the Euler equations above, it can be written as: \\
$\begin{pmatrix}
\rho \\
v
\end{pmatrix}_t + 
\begin{pmatrix}
v & \rho \\
\frac{c^2}{\rho} & v
\end{pmatrix}
\begin{pmatrix}
\rho \\
v
\end{pmatrix}_x = 0$

Consider $J := f(\rho, v)$ and try to make $J$ satisfy $\partial_t J + \lambda(\rho, v) \partial_x J = 0$. This is equivalent to: \\
$\begin{pmatrix}
\partial_\rho f & \partial_v f
\end{pmatrix}
\begin{pmatrix}
\rho \\
v
\end{pmatrix}_t + \lambda(\rho, v)
\begin{pmatrix}
\partial_\rho f & \partial_v f
\end{pmatrix}
\begin{pmatrix}
\rho \\
v
\end{pmatrix}_x = 0$ \\
then using the Euler equations to get: \\
$\begin{pmatrix}
\partial_\rho f & \partial_v f
\end{pmatrix}
\begin{pmatrix}
v & \rho \\
\frac{c^2}{\rho} & v
\end{pmatrix}
\begin{pmatrix}
\rho \\
v
\end{pmatrix}_x = \lambda(\rho, v) 
\begin{pmatrix}
\partial_\rho f & \partial_v f
\end{pmatrix}
\begin{pmatrix}
\rho \\
v
\end{pmatrix}_x$ \\
So, by taking a transpose we can see $\begin{pmatrix}
\partial_\rho f \\
\partial_v f
\end{pmatrix}$ is the eigenvector of the matrix 
$\begin{pmatrix}
v & \frac{c^2}{\rho} \\
\rho & v
\end{pmatrix}$, and $\lambda(\rho, v)$ is its eigenvalue. In this case, $\lambda = v \pm c$ and $\begin{pmatrix}
\partial_\rho f \\
\partial_v f
\end{pmatrix} = k
\begin{pmatrix}
\frac{c}{\rho} \\
\pm 1
\end{pmatrix}$, so $f$ can be chosen to be $k(\int \frac{c}{\rho} d\rho \pm v) = k(\frac{2}{\gamma - 1} c \pm v)$. Following Riemann's convention, we choose $k = \frac{1}{2}$. These are the Riemann invariants.

There are two Riemann invariants: 
\begin{equation*}
\begin{cases}
    w = \frac{1}{\gamma - 1} c - \frac{1}{2} v \\
    \underline{w} = \frac{1}{\gamma - 1} c + \frac{1}{2} v
\end{cases} \qquad \begin{cases}
    c = \frac{\gamma - 1}{2} (\underline{w} + w) \\
    v = \underline{w} - w
\end{cases}
\end{equation*}

And they satisfy the following equations.
\begin{equation*}
    \begin{cases}\partial_t w + (v - c) \partial_x w = 0 \\
\partial_t \underline{w} + (v + c) \partial_x \underline{w} = 0\end{cases}
\end{equation*}

We can also write down their explicit formula from the expression of $v, c$.
\begin{equation*}
\begin{cases}
    w = \frac{1}{\gamma - 1} c_1 \\
    \underline{w} = \frac{2}{\gamma + 1} (\frac{x}{t} - c_1) + \frac{1}{\gamma - 1} c_1
\end{cases}
\end{equation*}

In the case of 3-dimensional spherical symmetric case, similar Riemann invariants can be chosen, but the equations will be slightly changed.

Let $\mathbf{v} = \frac{v}{r} \mathbf{r}$, then the Euler equations can be written as: 
\begin{equation*}
\begin{cases}
    \partial_t \rho + v \partial_r \rho + \rho \partial_r v = - 2 \frac{\rho v}{r} \\
    \partial_t v + v \partial_r v + \frac{c^2}{\rho} \partial_r \rho = 0 \\
    \frac{p}{\rho^\gamma} = s = Const.
\end{cases}
\end{equation*}

There is only one extra term coming from $\nabla \cdot \mathbf{v}$.

Following the same process above, Riemann invariants can still be chosen as $w = \frac{1}{\gamma - 1} c - \frac{1}{2} v, \underline{w} = \frac{1}{\gamma - 1} c + \frac{1}{2} v$, the same expression as above. The Euler equations can be rewritten as follows.
\begin{equation*}
\begin{cases}
    \partial_t w + (v - c) \partial_r w = - \frac{c v}{r}, \\
    \partial_t \underline{w} + (v + c) \partial_r \underline{w} = - \frac{c v}{r}
\end{cases}
\end{equation*}

In terms of the vector fields $L_{(t, r)} = \partial_t + (v + c) \partial_r, \underline{L}_{(t, r)} = \frac{\kappa}{c} (\partial_t + (v - c) \partial_r)$ that we will use afterwards, the Euler equations can be rewritten as follows.
\begin{equation*}
\begin{cases}
    \underline{L} w = - \frac{\kappa v}{r}, \\
    L \underline{w} = - \frac{c v}{r}
\end{cases}
\end{equation*}

\section{Preparations}

In this section, we present the acoustical coordinates on which we will be working, define the necessary objects, and set up the framework of energy estimates in the acoustical coordinates.

The main part of this section follows the monograph \cite{Christodoulou-Miao-Shock} on the formation of shocks by Christodoulou and Miao. We adapt their general acoustical geometry settings to the spherical symmetric case, where we have major simplifications and can write many expressions in a more explicit form.

\subsection{Acoustical Metric, Coordinates and Frames}

We use $(t, x^i)$ or $(t, r), r = |x|$ to represent the Galilean coordinates in spherical symmetry. To avoid confusion, we put the coordinates in the subscript to emphasize that the expression or object comes from that coordinate system.

Assume that there has been a sufficiently regular solution for the Euler equations (which will be the case in the rarefaction wave region of the solution which we would like to construct), then we can define a Lorentz metric $g = g_{\mu\nu}$, which is the acoustical metric.
\begin{equation*}
    g = - c^2 dt^2 + \sum_i (dx^i - v^i dt)^2
\end{equation*}

It has the following properties:

When restricted to the spatial plane of constant time $\{(t, x) \mid t = t_0\}$, the restricted metric is the standard Euclidean metric on $\mathbb{R}^3$.

The null vectors are the directions in which the sound (linear perturbation of constant states) travels.

The future pointing vector field $B_{(t, x)} = \frac{\partial}{\partial t} + \sum_i v^i \frac{\partial}{\partial x^i} = \frac{\partial}{\partial t} + v \frac{\partial}{\partial r}$ is orthogonal to $\Sigma_t$, which is the material derivative, the direction (in spacetime) in which the fluid travels. \\

Then we define the vector fields and the outgoing null foliation that we need.

We first define the outgoing null vector field $L_{(t, r)} = \frac{\partial}{\partial t} + (v + c) \frac{\partial}{\partial r}$. It is chosen such that the time component in the Galilean coordinates is $1$. We don't need to solve for geodesic equations, because the spherical symmetry provides the canonical outgoing and ingoing directions.

Then we construct the outgoing null foliation using $L$. The foliation is parameterized by a radial function $u = u(t, r)$, called the acoustical function. The function $u$ is first manually set to be a smooth decreasing radial function $u(t_0, r)$ with non-zero gradient on an initial time slice $\{t = t_0\}$, then propagated to the entire spacetime region along the vector field $L$ by solving equation $Lu = 0$ using the characteristic line method.

Now we have the acoustical coordinates $(t, u, \theta^1, \theta^2)$, where $(\theta^1, \theta^2)$ are (local) angular coordinates which satisfy $L \theta^A = 0$. Sometimes we also abbreviate it to $(t, u)$, since we have spherical symmetry. We use $X_A = \partial_{\theta^A, (t, u, \theta^A)}$ to represent the angular coordinate vector fields.

We introduce several notations for geometric objects in the acoustical coordinates. $\Sigma_{t_0} = \{(t, u, \theta^1, \theta^2) \mid t = t_0\}, C_{u_0} = \{(t, u, \theta^1, \theta^2) \mid u = u_0\}, S_{t_0, u_0} = \Sigma_{t_0} \cap C_{u_0}$. They correspond to the time slice, the outgoing null cone, a 2-sphere respectively.

Then we define an ingoing spacelike vector field $T_{(t, u)} = \partial_u$ in the acoustical coordinates, whose expression in Galilean coordinates is $T_{(t, r)} = \partial_{u, (t, u)} r \partial_r$. It is tangent to the time slices $\Sigma_t$, and its length (evaluated in Euclidean metric or acoustical metric) is denoted as $\kappa$. $\kappa = \sqrt{g(T, T)} = - \partial_{u, (t, u)} r$. Intuitively, this function $\kappa$ represents the inverse of the spatial density of the foliation. We also define $\mu = - \frac{1}{(g^{-1}) (du, dt)} = c \kappa$, which is the inverse temporal density of the foliation. $\mu$ will be used during the derivation of the energy estimate framework, but in the end when we are doing the estimates, we always use the relation $\mu = c \kappa$ and only deal with $\kappa$.

Then we define the ingoing null vector field $\underline{L} = 2 T + \frac{\kappa}{c} L$. It is normalized such that we have $g(L, \underline{L}) = - 2 \mu$.

The vector fields $L, T$ are just the coordinate vector fields in the acoustical coordinates, $L_{(t, u)} = \partial_t, T_{(t, u)} = \partial_u$. So we have simple commutator formulas.
\begin{equation*}
    [L, T] = 0, \qquad [L, \underline{L}] = L (\frac{\kappa}{c}) L
\end{equation*}

We also define the restricted metric on $S_{t, u}$ to be $\slashed{g}_{A B} = r^2 g_{\mathbb{S}^2}$.

We use $D$ for the covariant derivative of $g$, and use $\slashed{D}$ for the restricted covariant derivative on $S_{t, u}$ of $\slashed{g}$.

Now we write the expression of $g$ and its inverse $g^{- 1}$ in the acoustical coordinates $(t, u, \theta^A)$, as well as the volume form.
\begin{equation*}
\begin{aligned}
    & g = - 2 \mu dt du + \kappa^2 du^2 + \slashed{g}_{A B} d\theta^A d\theta^B \\
    & g^{- 1} = - \frac{1}{2 \mu} (L^\mu \underline{L}^\nu + L^\nu \underline{L}^\mu) + \slashed{g}^{A B} X_A^\mu X_B^\nu \\
    & d\mu_{g} = \mu d\mu_{\slashed{g}} du dt
\end{aligned}
\end{equation*}

Then we define the second fundamental forms $\chi, \underline{\chi}$.
\begin{equation*}
\begin{aligned}
    & \chi(X_A, X_B) = \frac{1}{2} (\mathcal{L}_L g)(X_A, X_B) = g(D_{X_A} L, X_B) \\
    & \underline{\chi} = \frac{1}{2} (\mathcal{L}_{\underline{L}} g)(X_A, X_B) = g(D_{X_A} \underline{L}, X_B) \\
\end{aligned}
\end{equation*}

They have the following expression since we are in spherical symmetry.
\begin{equation*}
\begin{aligned}
    & \chi = \frac{L r}{r} \slashed{g} = \frac{v + c}{r} \slashed{g} && \mathrm{tr} \chi = 2 \frac{v + c}{r} \\
    & \underline{\chi} = \frac{\underline{L} r}{r} \slashed{g} = \frac{\kappa}{c} \frac{v - c}{r} \slashed{g} && \mathrm{tr} \underline{\chi} = 2 \frac{\kappa}{c} \frac{v - c}{r} \\
\end{aligned}
\end{equation*}

We also list the connection table for the frame $(L, \underline{L}, X_A)$.
\begin{equation*}
\begin{cases}
\begin{aligned}
    & D_L L = \mu^{-1} (L \mu) L && D_{X_A} L = {\chi_A}^B X_B \\
    & D_{\underline{L}} L = - L(c^{-1} \kappa) L && D_{X_A} \underline{L} = {\underline{\chi}_A}^B X_B \\
    & D_L \underline{L} = 0 && D_L X_A = {\chi_A}^B X_B (= D_{X_A} L) \\
    & D_{\underline{L}} \underline{L} = (\mu^{-1} \underline{L} \mu + L(c^{-1} \kappa)) \underline{L}  && D_{\underline{L}} X_A = {\underline{\chi}_A}^B X_B (= D_{X_A} \underline{L}) \\
\end{aligned} \\
    D_{X_A} X_B = \slashed{D}_{X_A} X_B + \frac{1}{2} \mu^{-1} \underline{\chi}_{AB} L + \frac{1}{2} \mu^{-1} \chi_{AB} \underline{L}
\end{cases}
\end{equation*}

We list the Euler equations and the transport equation for $r$ in the acoustical coordinates.
\begin{equation} \label{eq:full_equation_in_acoustic}
\begin{cases}
    \underline{L} w = - \frac{\kappa v}{r} \\
    L \underline{w} = - \frac{c v}{r} \\
    L r = v + c
\end{cases}
\end{equation}

We will mainly work in the acoustical coordinates and deal with this system of equations.

\subsection{Wave Operator, Energy Estimates}

Now we can express the wave operator $\Box_g$ in the acoustical coordinates and write it in the form of $L, \underline{L}, T$.
\begin{equation*}
    \Box_g f = g^{\mu \nu} (D^2 f)_{\mu \nu} = - \mu^{-1} (D^2 f) (L, \underline{L}) + \slashed{g}^{AB} (D^2 f) (X_A, X_B)
\end{equation*}

We compute the two parts.
\begin{equation*}
\begin{aligned}
    (D^2 f) (X_A, X_B) & = (\slashed{D}^2 f) (X_A, X_B) - \frac{1}{2} \mu^{-1} \underline{\chi}_{AB} (Lf) - \frac{1}{2} \mu^{-1} \chi_{AB} (\underline{L}f) \\
    & = - \frac{1}{2} \mu^{-1} \underline{\chi}_{AB} (Lf) - \frac{1}{2} \mu^{-1} \chi_{AB} (\underline{L}f)
\end{aligned}
\end{equation*}

\begin{equation*}
    (D^2 f) (L, \underline{L}) = L(\underline{L}(f)) - (D_L \underline{L})f = L (\underline{L} (f))
\end{equation*}

The vanishing terms are due to spherical symmetry. Then we obtain the following expressions.
\begin{equation*}
\begin{aligned}
    \mu \Box_g f & = - \frac{1}{2} \mathrm{tr} \underline{\chi} (Lf) - \frac{1}{2} \mathrm{tr} \chi (\underline{L}f) - L(\underline{L}(f)) \\
    & = - \frac{1}{2} \mathrm{tr} \underline{\chi} (Lf) - \frac{1}{2} \mathrm{tr} \chi (\underline{L}f) - \underline{L}(L(f)) - L(\frac{\kappa}{c}) (Lf) \\
    & = - \frac{2 \kappa v}{c r} Lf - 2 \frac{v + c}{r} Tf - \frac{\kappa}{c} LLf - 2 LTf - L(\frac{\kappa}{c}) Lf
\end{aligned}
\end{equation*}

Then we can use the standard stress-energy-momentum tensor to get the energy estimates for a general scalar wave equation $\Box_g \psi = \tilde{\rho}$.

For a scalar function $\psi$, we define its energy-momentum tensor to be the following.
\begin{equation*}
    T_{\mu \nu} = \partial_\mu \psi \partial_\nu \psi - \frac{1}{2} g_{\mu \nu} g^{\alpha \beta} \partial_\alpha \psi \partial_\beta \psi
\end{equation*}

And we can compute its divergence.
\begin{equation*}
\begin{aligned}
    D^\mu T_{\mu \nu} & = D^\mu (\partial_\mu \psi \partial_\nu \psi) - \frac{1}{2} D^\mu (g_{\mu \nu} g^{\alpha \beta} \partial_\alpha \psi \partial_\beta \psi) \\
    & = \Box_g \psi \partial_\nu \psi + \partial_\mu \psi g^{\mu \xi} D_\xi \partial_\nu \psi - \frac{1}{2} g^{\alpha \beta} D_\nu \partial_\alpha \psi \partial_\beta \psi - \frac{1}{2} g^{\alpha \beta} \partial_\alpha \psi D_\nu \partial_\beta \psi \\
    & = \Box_g \psi \partial_\nu \psi \\
    & = \tilde{\rho} \partial_\nu \psi
\end{aligned}
\end{equation*}

Now we take the timelike multiplier $X$ to be a combination of $L, \underline{L}$.
\begin{equation*}
    X = a L + b \underline{L}
\end{equation*}

And we define the energy current.
\begin{equation*}
    P_\mu = - T_{\mu \nu} X^\nu
\end{equation*}

\begin{equation*}
\begin{aligned}
    D^\mu P_\mu & = - (D^\mu T_{\mu \nu}) X^\nu - T_{\mu \nu} D^\mu X^\nu \\
    & = - \tilde{\rho} X \psi - T_{\mu \nu} D^\mu X^\nu \\
    & = Q
\end{aligned}
\end{equation*}

Now we list their components in the frame $(L, \underline{L}, X_A)$.
\begin{equation*}
\begin{cases}
\begin{aligned}
    & T_{L L} = |L\psi|^2 && T_{L \underline{L}} = \mu |\slashed{d}\psi|_{\slashed{g}}^2 = 0 && T_{\underline{L} \underline{L}} = |\underline{L}\psi|^2 \\
    & T_{L X_A} = L\psi \slashed{d}_A \psi = 0 && T_{\underline{L} X_A} = \underline{L}\psi \slashed{d}_A \psi = 0 &&
\end{aligned} \\
    T_{X_A X_B} = \slashed{d}_A \psi \slashed{d}_B \psi + \frac{1}{2 \mu} \slashed{g}_{A B} L\psi \underline{L}\psi - \frac{1}{2} \slashed{g}_{A B} |\slashed{d}\psi|_{\slashed{g}}^2 = \frac{1}{2 \mu} \slashed{g}_{A B} L\psi \underline{L}\psi
\end{cases}
\end{equation*}

\begin{equation*}
\begin{cases}
    P^L = \frac{b}{2 \mu} |\underline{L}\psi|^2 + \frac{a}{2} |\slashed{d}\psi|_{\slashed{g}}^2 = \frac{b}{2 \mu} |\underline{L}\psi|^2 \\
    P^{\underline{L}} = \frac{a}{2 \mu} |L\psi|^2 + \frac{b}{2} |\slashed{d}\psi|_{\slashed{g}}^2 = \frac{a}{2 \mu} |L\psi|^2 \\
    P^{X_A} = - (a L\psi + b \underline{L}\psi) \slashed{d}^A \psi = 0 \\
    P^t = P^L + c^{-1} \kappa P^{\underline{L}} = \frac{b}{2 \mu} |\underline{L}\psi|^2 + c^{-1} \kappa \frac{a}{2 \mu} |L\psi|^2 \\
    P^u = 2 P^{\underline{L}} = \frac{a}{\mu} |L\psi|^2
\end{cases}
\end{equation*}

Now we apply the divergence theorem on the region $(t, u, \theta^A) \in [t_1, t_2] \times [u_1, u_2] \times \mathbb{S}^2$ to obtain the following identity.
\begin{equation*}
\begin{cases}
    \mathcal{E}^{[u_1, u_2]}(t_2) - \mathcal{E}^{[u_1, u_2]}(t_1) + \mathcal{F}^{[t_1, t_2]}(u_2) - \mathcal{F}^{[t_1, t_2]}(u_1) = \int_{\{t_1 \leq t \leq t_2, u_1 \leq u \leq u_2\}} Q d\mu_g \\
    \mathcal{E}^{[u_1, u_2]}(t) = \int_{[u_1, u_2]} \int_{S_{t, u}} \mu P^t d\mu_\slashed{g} du \\
    \mathcal{F}^{[t_1, t_2]}(u) = \int_{[t_1, t_2]} \int_{S_{t, u}} \mu P^u d\mu_\slashed{g} dt
\end{cases}
\end{equation*}

Now we substitute all the frame expressions into $P^\mu, Q$ to obtain the energy estimate framework we need. We summarize them into the following proposition.

\begin{proposition}[Energy Estimate Framework]
Let $\psi$ to be a scalar function defined on $(t, u) \in [t_1, t_2] \times [u_1, u_2]$ satisfying the wave equation $\Box_g \psi = \tilde{\rho}$. We have the following energy estimate corresponding to the multiplier $X = a L + b \underline{L}$.
\begin{equation} \label{eq:energy_estimate}
\begin{cases}
    \Box_g \psi = \tilde{\rho} \\
    \mathcal{E}^{[u_1, u_2]}(t_2) - \mathcal{E}^{[u_1, u_2]}(t_1) + \mathcal{F}^{[t_1, t_2]}(u_2) - \mathcal{F}^{[t_1, t_2]}(u_1) = \int_{\{t_1 \leq t \leq t_2, u_1 \leq u \leq u_2\}} Q d\mu_g \\
    \mathcal{E}^{[u_1, u_2]}(t) = \int_{[u_1, u_2]} \int_{S_{t, u}} \frac{1}{2} (b |\underline{L}\psi|^2 + a \frac{\kappa}{c} |L\psi|^2) d\mu_\slashed{g} du \\
    \mathcal{F}^{[t_1, t_2]}(u) = \int_{[t_1, t_2]} \int_{S_{t, u}} a |L\psi|^2 d\mu_\slashed{g} dt \\
    Q = Q_0 + Q_1 + Q_2 + Q_3 \\
    Q_0 = - \tilde{\rho} (a L\psi + b \underline{L}\psi) \\
    Q_1 = \frac{1}{2 \mu} Lb |\underline{L}\psi|^2 \\
    Q_2 = \frac{1}{2 \mu} (\underline{L}a - a L(\frac{\kappa}{c})) |L\psi|^2 \\
    Q_3 = - \frac{1}{2 \mu} (a \mathrm{tr}\chi + b \mathrm{tr}\underline{\chi}) L\psi \underline{L}\psi
\end{cases}
\end{equation}
\end{proposition}

We will apply this framework to $(w, \underline{w})$ during the proof.

\section{Local Solution for General Background}

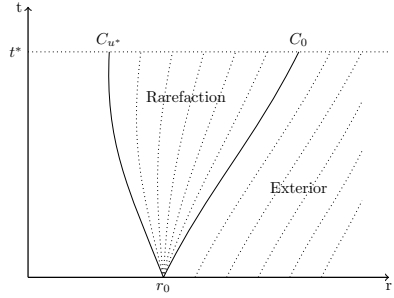
\begin{figure}[H]
\centering
\scalebox{0.6}{\begin{tikzpicture}

% coordinate axis
\coordinate (origin) at (0,0);
\coordinate (r-axis) at (8,0);
\coordinate (t-axis) at (0,6);
\node at (r-axis) [anchor=north]{r};
\node at (t-axis) [anchor=east]{t};
\draw [thick,->] (origin) -- (r-axis);
\draw [thick,->] (origin) -- (t-axis);
\draw [dotted,semithick] (0,5) -- (8,5);
\node at (0,5) [anchor=east]{$t^*$};
% C_0
\draw [semithick] (3,0) node [anchor=north]{$r_0$} .. controls +(1,2) and +(-1,-2) .. +(3,5) node [anchor=south]{$C_0$};
% Rarefaction
\foreach \X in {1,...,5}
{\draw [dotted,semithick] (3,0) .. controls +(1-0.3*\X,2) and +(-1+0.15*\X,-2) .. +(3-0.7*\X,5);};
% Rarefaction End
\draw [semithick] (3,0) .. controls +(1-0.3*6,2) and +(-1+0.15*6,-2) .. +(3-0.7*6,5) node [anchor=south]{$C_{u^*}$};;
% Exterior
\foreach \X in {1,...,2}
{\draw [dotted,semithick] (3+0.7*\X,0) .. controls +(1-0.1*\X,2-0.2*\X) and +(-1+0.1*\X,-2+0.2*\X) .. +(3,5);};
\foreach \X in {3,...,5}
{\draw [dotted,semithick] (3+0.7*\X,0) .. controls +(1-0.1*\X,2-0.2*\X) and +(-1+0.1*\X,-2+0.2*\X) .. +(3+0.7*2-0.7*\X,5+1.1*2-1.1*\X);};
% Text
\node at (3.5,4) {Rarefaction};
\node at (6,2) {Exterior};

\end{tikzpicture}}
\caption{Local Rarefaction}
\label{fig:local2}
\end{figure}

In this section, we will show the local existence of the rarefaction wave solution for general background flow. It can be formulated into the following theorem:

\begin{theorem}[Local Existence of Rarefaction Solution] \label{thm:local}
Given any smooth spherical symmetric functions $(v_{\mathrm{ext}}, c_{\mathrm{ext}})$ solving the following Euler equations in the exterior region $D_{\mathrm{ext}}$.
\begin{equation*}
\begin{cases}
    \partial_t w + (v - c) \partial_r w = - \frac{c v}{r} \\
    \partial_t \underline{w} + (v + c) \partial_r \underline{w} = - \frac{c v}{r}
\end{cases} \qquad \begin{cases}
    w = \frac{1}{\gamma - 1} c - \frac{1}{2} v \\
    \underline{w} = \frac{1}{\gamma - 1} c + \frac{1}{2} v
\end{cases}
\end{equation*}

Where $w, \underline{w}$ are the Riemann invariants defined in subsection \ref{subsec:1D}, and the region $D_\mathrm{ext}$ is the exterior of an outgoing characteristic surface $C_0$.
\begin{equation*}
\begin{cases}
    D_{\mathrm{ext}} = \{(t, r) \mid t \in [0, t_0], r \geq r_{C_0}(t)\} \\
    C_0 = \{(t, r) \mid r = r_{C_0}(t)\} \\
    \frac{d}{dt} r_{C_0}(t) = (v_{\mathrm{ext}} + c_{\mathrm{ext}})(t, r_{C_0}(t)) \\
    r_{C_0}(0) = r_0, \quad c_{\mathrm{ext}}(0, r_0) = c_0, \quad v_{\mathrm{ext}}(0,r_0) = v_0
\end{cases}
\end{equation*}

Then there exists a rarefaction wave extension $(w_{\mathrm{rf}}, \underline{w}_{\mathrm{rf}})$ of the solution to a fan-shaped region $D_\mathrm{rf} = \{(t, r) \mid t \in [0, t^*], r \in [r_{C_{u^*}}(t), r_{C_0}(t)]\}, t^* \leq t_0$, foliated by an outgoing null foliation parameterized by the acoustical function $u \in [0, u^*]$. Level sets $C_u = \{(t, r) \mid r = r_{C_u}(t)\}$ are outgoing characteristic surfaces. See figure \ref{fig:local2}.

The rarefaction wave solution satisfies the following properties.
\begin{itemize}
    \item Continuity across the boundary. $(w_{\mathrm{rf}}, \underline{w}_{\mathrm{rf}}) = (w_{\mathrm{ext}}, \underline{w}_{\mathrm{ext}})$ on $C_0 \setminus \{(0, r_0)\}$.
    
    \item Discontinuity of the transversal derivative across the boundary. More precisely, we have $\partial_r w_{\mathrm{rf}} = \partial_r w_{\mathrm{ext}}, \partial_r \underline{w}_{\mathrm{rf}} \neq \partial_r \underline{w}_{\mathrm{ext}}$ on $C_0 \setminus \{(0, r_0)\}$.

    \item Smoothness in acoustical coordinates and in radial coordinates. $(w_{\mathrm{rf}}, \underline{w}_{\mathrm{rf}}, r)$ are smooth functions of $(t, u)$ in $[0, t^*] \times [0, u^*]$. They have estimates \ref{eq:local-base-final_estimates} and \ref{eq:local-high-final_estimates}. Inversely, $(w_{\mathrm{rf}}, \underline{w}_{\mathrm{rf}}, u)$ are smooth functions of $(t, r)$ in $D_{\mathrm{rf}} \setminus \{(0, r_0)\}$.

    \item The two pieces of solutions together form a continuous, piecewise smooth weak solution of the Euler equations with initial singularity.
\end{itemize}

It also has the following limiting behaviors at the initial singularity.
\begin{itemize}
    \item All outgoing null surfaces originate from the initial singularity. $r_{C_u}(0) = r_0$.

    \item The foliation behaves like $\frac{x}{t}$ in the 1-D case when approaching the initial singularity. $\lim\limits_{t \to 0} \frac{r_{C_u}(t) - r_0}{t} = (v_{\mathrm{rf}} + c_{\mathrm{rf}})_{(t, u)}(0, u) = (v_{\mathrm{ext}} + c_{\mathrm{ext}})_{(t, r)}(0, r_0) - u$.

    \item The rarefaction solution behaves like the 1-D pattern when approaching the initial singularity.
    \begin{equation*}
    \begin{cases}
        w_{\mathrm{rf}}(0, u) = w_{\mathrm{ext}}(0, r_0) \\
        \underline{w}_{\mathrm{rf}}(0, u) = \underline{w}_{\mathrm{ext}}(0, r_0) - \frac{2}{\gamma + 1} u
    \end{cases} \qquad \begin{cases}
        c_{\mathrm{rf}}(0, u) = c_0 - \frac{\gamma - 1}{\gamma + 1} u \\
        v_{\mathrm{rf}}(0, u) = v_0 - \frac{2}{\gamma + 1} u
    \end{cases}
    \end{equation*}
\end{itemize}

$u^* \in (0, \frac{\gamma + 1}{\gamma - 1} c_0)$ can be arbitrarily chosen. The life span $t^*$ of the rarefaction solution depends on $u^*$ and up to $4$ derivatives of the exterior solution.
\end{theorem}

\begin{remark}
We are not dealing with the center $r = 0$ and the vacuum $c = 0$. The construction relies on $r \sim r_0, c \sim c_0, v \lesssim c_0$.
\end{remark}

\subsection{Determine Data on the Boundary $C_0$}

We first do some priori calculations. Assuming that we already have a desired center rarefaction wave solution with enough regularity, we want to see what quantities in the rarefaction wave region are determined by the exterior solution on their common boundary $C_0$.

The desired solution consists of two parts: an exterior solution in $D_{\mathrm{ext}}$, smooth up to the boundary, and a rarefaction wave solution in a fan-shaped region $D_{\mathrm{rf}}$, smooth up to the boundary, except at the initial singular point. Their boundaries consist of two outgoing null cones $C_0, C_{u^*}$. In addition, the rarefaction wave region has a smooth outgoing null foliation parameterized by $u$. The acoustical coordinates are established in the rarefaction wave region as above.

Recall \ref{eq:full_equation_in_acoustic}. The Euler equations written in the acoustical frames and the Riemann invariants take the following form: 
\begin{equation*}
\begin{cases}
    \underline{L} w = - \frac{\kappa v}{r} \\
    L \underline{w} = - \frac{c v}{r}
\end{cases}
\end{equation*}

Decompose $\underline{L} = 2 T + \frac{\kappa}{c} L$, and apply $T$ to the equation of $\underline{w}$ to obtain the equations for $w, \underline{w}$. 
\begin{equation*}
\begin{cases}
    T w = - (\frac{v}{2 r} + \frac{L w}{2 c}) \kappa \\
    L T \underline{w} = \frac{\gamma - 1}{2} (- 2 \frac{\underline{w}}{r} T \underline{w} + 2 \frac{w}{r} T w - \frac{\underline{w}^2 - w^2}{r^2} \kappa) 
\end{cases}
\end{equation*}

There is also the transport equation for $\kappa$, obtained by applying $T$ to the transport equation for $r$: 
\begin{equation*}
    L \kappa = - \frac{\gamma + 1}{2} T \underline{w} + \frac{3 - \gamma}{2} T w
\end{equation*}

We fix $u = 0$ and take $t$ as the only variable, then these equations can be rewritten as an ODE system for $(T w, T \underline{w}, \kappa)$ on $C_0$, because $(w, \underline{w}, r)$ are determined by the exterior solution due to the continuity across the boundary $C_0$. 
\begin{equation} \label{eq:ODE-first}
\begin{cases}
    T w = - (\frac{v}{2 r} + \frac{L w}{2 c}) \kappa \\
    L T \underline{w} = \frac{\gamma - 1}{2} (- 2 \frac{\underline{w}}{r} T \underline{w} + 2 \frac{w}{r} T w - \frac{\underline{w}^2 - w^2}{r^2} \kappa) \\
    L \kappa = - \frac{\gamma + 1}{2} T \underline{w} + \frac{3 - \gamma}{2} T w
\end{cases}
\end{equation}

The first equation gives a representation of $T w$ by $\kappa$, so the ODE part only involves $(T \underline{w}, \kappa)$. We still put them together to emphasize that $(T w, T \underline{w}, \kappa)$ are of the same order in a sense.

By standard ODE results, $(T w, T \underline{w}, \kappa)$ can be solved on $C_0$ after fixing the initial data at $t = 0$. We assume that the initial data are given for now and continue developing the scheme for solving high order derivatives, namely $(T^n w, T^n \underline{w}, T^{n - 1} \kappa)$ on $C_0$.

After solving $(Tw, T\underline{w}, \kappa)$, they can also be treated as known smooth functions of $t$, just like $(w, \underline{w}, r)$. So by differentiating the equations (applying $T$ on both sides) repeatedly, we can get similar ODE systems for $(T^n w, T^n \underline{w}, T^{n - 1} \kappa)$, after solving $(T^{\leq n - 1} w, T^{\leq n - 1} \underline{w}, T^{\leq n - 2} \kappa)$. 
\begin{equation} \label{eq:ODE-general}
\begin{cases}
    T^n w = - (\frac{v}{2 r} + \frac{Lw}{2 c}) T^{n - 1} \kappa + \sum\limits_{\substack{n_1 + n_2 = n - 1 \\ n_i \geq 0,\ n_2 \leq n - 2}} \mathrm{C}_{n_1, n_2} T^{n_1} (\frac{v}{r} + \frac{Lw}{c}) T^{n_2} \kappa \\
    L T^n \underline{w} = \frac{\gamma - 1}{2} (- 2 \frac{\underline{w}}{r} T^n \underline{w} + 2 \frac{w}{r} T^n w - \frac{\underline{w}^2 - w^2}{r^2} T^{n - 1} \kappa) \\
    \qquad \qquad + \sum\limits_{\substack{n_1 + n_2 = n - 1 \\ n_i \geq 0,\ n_2 \leq n - 2}} \mathrm{C}_{n_1, n_2} (- 2 T^{n_1} \frac{\underline{w}}{r} T^{n_2 + 1} \underline{w} + 2 T^{n_1} \frac{w}{r} T^{n_2 + 1} w - T^{n_1} \frac{\underline{w}^2 - w^2}{r^2} T^{n_2} \kappa) \\
    L T^{n - 1} \kappa = - \frac{\gamma + 1}{2} T^n \underline{w} + \frac{3 - \gamma}{2} T^n w
\end{cases}
\end{equation}

The terms in the summation are all lower order terms. They are combinations of $(L^{\leq 1} T^{\leq n - 1} w, T^{\leq n - 1} \underline{w}, T^{\leq n - 2} \kappa)$, which have already been solved, therefore can be treated as known.

By solving these ODE systems recursively using the manually set initial data at $(t, u) = (0, 0)$, the quantities $(T^n w, T^n \underline{w}, T^{n - 1} \kappa)$, which represent the information of the Taylor expansion in $u$ of the solution on $C_0$, are uniquely determined. They are all smooth functions of $t$ since the coefficients of the ODE system are smooth in each step. \\

We now turn to the choice of initial data at the singular point $(t = 0, u = 0)$.

The way we choose the initial data comes from the 1-D case, following the intuition that the 1-D picture will dominate as we approach the initial singularity. Recall that in the 1-D case, we have 
\begin{equation*}
\begin{cases}
    w = \frac{1}{\gamma - 1} c_1 \\
    \underline{w} = \frac{2}{\gamma + 1} (\frac{x}{t} - c_1) + \frac{1}{\gamma - 1} c_1
\end{cases}
\end{equation*}

Notice that the Riemann invariants are functions of $\frac{x}{t}$. If we set $u = - \frac{x}{t} + \mathrm{Const.}$, then we will get a natural choice of the initial data at $(t, u) = (0, 0)$.
\begin{equation} \label{eq:initial_data}
\begin{cases}
    \begin{aligned}
    & T w = 0, \quad && T^n w = 0 \\
    & T \underline{w} = - \frac{2}{\gamma + 1}, \quad && T^n \underline{w} = 0 \\
    & \kappa = 0, \quad && T^{n - 1} \kappa = 0
    \end{aligned}, \quad (n \geq 2)
\end{cases}
\end{equation}

Technically, we should not set the initial data for $T^n w$ according to the structure of the system. We still list them above for the sake of completeness, and one can verify that the initial data satisfy the equations for $T^n w$. Using these initial data, we can solve the ODE systems and obtain smooth functions $(T^n w, T^n \underline{w}, T^{n - 1} \kappa)$ of $t$ on $C_0$. \\

We also need to check which quantities vanish at $t = 0$. The initial data that we manually choose give $Tw, \kappa, T^{n + 1}w, T^{n + 1}\underline{w}, T^n \kappa = 0, (n \geq 1, t = 0)$. Note that the initial data for $Tw, T^{n + 1}w$ should be computed by applying $T$ on $Tw = - (\frac{v}{2 r} + \frac{Lw}{2 c}) \kappa$ repeatedly. We can also directly read from the equation $LT^{n - 1} \kappa = - \frac{\gamma + 1}{2} T^n \underline{w} + \frac{3 - \gamma}{2} T^n w, (n \geq 1)$ that $L\kappa = 1, LT^n\kappa = 0, (n \geq 1, t = 0)$.

Putting these together, we obtain the vanishing orders of the quantities at $t = 0$.
\begin{equation} \label{eq:local-vanishing_C0}
\begin{cases}
    \begin{aligned}
    & T w = O(t), && T^n w = O(t) \\
    & T \underline{w} = - \frac{2}{\gamma + 1} + O(t), && T^n \underline{w} = O(t) \\
    & \kappa = t + O(t^2), && T^{n - 1} \kappa = O(t^2)
    \end{aligned}, \quad (n \geq 2)
\end{cases}
\end{equation}

We summarize the discussions above into the following proposition.

\begin{proposition}[Construct Data on $C_0$] \label{prop:initial_C_0}
Given a smooth, spherical symmetric background solution $(c, v)$ to the Euler equations in the region $D_{\mathrm{ext}} = \{t \in [0, t_0], r \geq r_{C_0}(t)\}$ and an outgoing null hypersurface $C_0 = \{r = r_{C_0}(t) \}$, we can uniquely determine smooth functions $(T^n w, T^n \underline{w}, T^{n - 1} \kappa), (n \geq 1)$ of $t$ on $C_0$. They satisfy the propagation equations \ref{eq:ODE-first}, \ref{eq:ODE-general}, and the initial data condition \ref{eq:initial_data}. They also have vanishing orders \ref{eq:local-vanishing_C0}.
\end{proposition}

\subsection{Construct Local Approximate Solutions} \label{subsec:local-localexist}

After acquiring the data on $C_0$, we now begin to construct a family of local solutions of the Euler equations, which will converge to the rarefaction solution.

\begin{figure}[H]
\centering
\scalebox{0.6}{\begin{tikzpicture}

% Coordinate Axis
\coordinate (origin1) at (0,0);
\coordinate (r-axis1) at (6,0);
\coordinate (t-axis1) at (0,6);
\node at (r-axis1) [anchor=north]{r};
\node at (t-axis1) [anchor=east]{t};
\draw [thick,->] (origin1) -- (r-axis1);
\draw [thick,->] (origin1) -- (t-axis1);
\coordinate (origin2) at (10,0);
\coordinate (r-axis2) at (16,0);
\coordinate (t-axis2) at (10,6);
\node at (r-axis2) [anchor=north]{r};
\node at (t-axis2) [anchor=east]{t};
\draw [thick,->] (origin2) -- (r-axis2);
\draw [thick,->] (origin2) -- (t-axis2);
% C_0
\coordinate (singular1) at (2.5,0);
\node at (singular1) [anchor=north]{$r_0$};
\draw [semithick] (singular1) .. controls +(0.5,1.2) and +(-0.5,-0.8) .. ++(1,2) .. controls +(0.5,0.8) and +(-0.5,-1) .. +(2,3) node [anchor=south]{$C_0$};
\coordinate (singular2) at (12.5,0);
\node at (singular2) [anchor=north]{$r_0$};
\draw [semithick] (singular2) .. controls +(0.5,1.2) and +(-0.5,-0.8) .. ++(1,2) .. controls +(0.5,0.8) and +(-0.5,-1) .. +(2,3) node [anchor=south]{$C_0$};
% Sigma_t
\draw [dotted,semithick] (0,2) node [anchor=east]{$\delta$} -- (6,2);
\draw [dotted,semithick] (0,5) node [anchor=east]{$t^*$} -- (6,5);
\draw [dotted,semithick] (10,5) node [anchor=east]{$t^*$} -- (16,5);
% Solution
\fill [gray,opacity=0.5] (3.5,2) .. controls +(0.5,0.8) and +(-0.5,-1) .. +(2,3) -- (1.5,5) .. controls +(0,-1) and +(-0.3,1) .. (2,2) -- cycle;
\draw [semithick] (1.5,5) node [anchor=south]{$C_{u^*}$} .. controls +(0,-1) and +(-0.3,1) .. (2,2);
\fill [gray,opacity=0.5] (singular2) .. controls +(0.5,1.2) and +(-0.5,-0.8) .. ++(1,2) .. controls +(0.5,0.8) and +(-0.5,-1) .. +(2,3) -- (11.5,5) .. controls +(0,-1) and +(-1,2) .. (singular2) -- cycle;
\draw [semithick] (11.5,5) node [anchor=south]{$C_{u^*}$} .. controls +(0,-1) and +(-1,2) .. (singular2);
% Text
\node at (3,4) {Approximate};
\node at (13,4) {Exact};
\draw [semithick,->] (7,1.5) -- (9,1.5);
\node at (8,2) {$\delta \to 0$};
\node at (8,1) {Convergence};

\end{tikzpicture}}
\caption{Convergence}
\label{fig:convergence}
\end{figure}
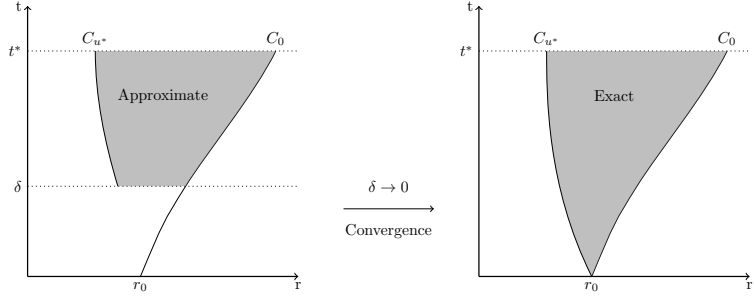

There are no local existence results for problems with singular initial data, so we manually construct the initial data on $\Sigma_\delta$ and solve from $t = \delta$. In the end, we will send $\delta \to 0$ and make a convergence argument to obtain the exact rarefaction wave solution. See figure \ref{fig:convergence}.

The construction is straightforward. We use the data on $C_0$ provided by proposition \ref{prop:initial_C_0} to construct Taylor polynomials on $\{t = \delta, u \leq u^*\}$.
\begin{equation} \label{eq:initial_Taylor}
\begin{cases}
\begin{aligned}
    & w_N(\delta, u) = w(\delta, 0) + \sum\limits_{n = 1}^N \frac{T^n w(\delta, 0)}{n!} u^n && (L w)_N(\delta, u) = - \frac{c_N v_N}{r_N} - \frac{2 c_N}{\kappa_N} T w_N \\
    & \underline{w}_N(\delta, u) = \underline{w}(\delta, 0) + \sum\limits_{n = 1}^N \frac{T^n \underline{w}(\delta, 0)}{n!} u^n && (L \underline{w})_N(\delta, u) = - \frac{c_N v_N}{r_N} \\
    & r_N(\delta, u) = r(\delta, 0) - \sum\limits_{n = 0}^{N - 1} \frac{T^n \kappa(\delta, 0)}{n!} u^n && (L r)_N(\delta, u) = L r_N(\delta, u) = c_N + v_N \\
\end{aligned}
\end{cases}
\end{equation}

We call these data the $N$-th order approximate data at time $t = \delta$.

The $L$-derivatives of the approximation data come from the equations. We note that now $(L \psi)_N \neq \partial_t \psi_N, \psi = w, \underline{w}$, they contain higher order terms in $u$. That is why we put parentheses.

We list the $L$-derivatives because we will use the wave equation results for the local well-posedness and energy estimates.

We now state the local well-posedness proposition. It is not sharp, but it is enough for our case.

\begin{proposition}[Local well-posedness of Cauchy-Characteristic Mixed Initial Value Problem]
Fix smooth data $(T^{\geq 0} w_{C_0}, T^{\geq 0} \underline{w}_{C_0}, T^{\geq 0} r_{C_0}) \in C_t^\infty$ on $C_0 = \{(t, u) \mid u = 0\} = \{(t, r) \mid r = r_{C_0}(t)\}$ given by proposition \ref{prop:initial_C_0}. Fix $N \geq 3$.

Given any data $(w_{\Sigma_\delta}, \underline{w}_{\Sigma_\delta}, r_{\Sigma_\delta}) \in H_u^{N + 1}([0, u^*])$ satisfying the following compatibility condition, there exists a unique local solution $(w, \underline{w}, r)$ of the Euler equations in the region $(t, u) \in [\delta, \delta + \tau] \times [0, u^*]$.
\begin{equation*}
    (T^{\leq N} w_{\Sigma_\delta}, T^{\leq N} \underline{w}_{\Sigma_\delta}, T^{\leq N} r_{\Sigma_\delta})(\delta, 0) = (T^{\leq N} w_{C_0}, T^{\leq N} \underline{w}_{C_0}, T^{\leq N} r_{C_0})(\delta)
\end{equation*}

The solution has the following regularity.
\begin{equation*}
\begin{cases}
    w \in C_t^0 H_u^{N + 1} \cap C_t^1 H_u^N \\
    \underline{w} \in C_t^0 H_u^{N + 1} \cap C_t^1 H_u^N \\
    r \in C_t^1 H_u^{N + 1}
\end{cases}
\end{equation*}

The life span $\tau$ depends on the data on $C_0$ (already fixed) and also on the following quantities.
\begin{equation*}
\begin{cases}
    \max\limits_{\Sigma_\delta} c, \qquad \min\limits_{\Sigma_\delta} c > 0, \qquad \max\limits_{\Sigma_\delta} |v|, \qquad \max\limits_{\Sigma_\delta} r, \qquad \min\limits_{\Sigma_\delta} r > 0, \qquad \max\limits_{\Sigma_\delta} \kappa, \qquad \min\limits_{\Sigma_\delta} \kappa > 0 \\
    u^*, \delta, \qquad \|(w_{\Sigma_\delta}, \underline{w}_{\Sigma_\delta}, r_{\Sigma_\delta})\|_{H_u^{N + 1}([0, u^*])}, \qquad \|((L w)_{\Sigma_\delta}, (L \underline{w})_{\Sigma_\delta}, (L r)_{\Sigma_\delta})\|_{H_u^N([0, u^*])}
\end{cases}
\end{equation*}

$((L w)_{\Sigma_\delta}, (L \underline{w})_{\Sigma_\delta}, (L r)_{\Sigma_\delta})$ are given by the equations as in \ref{eq:initial_Taylor}.
\end{proposition}

\begin{figure}[H]
\centering
\scalebox{0.6}{\begin{tikzpicture}

% Coordinates
\coordinate (origin1) at (4,0);
\coordinate (u-axis-R1) at (8,0);
\coordinate (u-axis-L1) at (0,0);
\coordinate (t-axis1) at (4,6);
\node at (u-axis-L1) [anchor=north]{u};
\node at (t-axis1) [anchor=east]{t};
\draw [thick,->] (u-axis-R1) -- (u-axis-L1);
\draw [thick,->] (origin1) -- (t-axis1);
\coordinate (origin2) at (10,0);
\coordinate (r-axis2) at (18,0);
\coordinate (t-axis2) at (10,6);
\node at (r-axis2) [anchor=north]{r};
\node at (t-axis2) [anchor=east]{t};
\draw [thick,->] (origin2) -- (r-axis2);
\draw [thick,->] (origin2) -- (t-axis2);
\coordinate (origin3) at (24,0);
\coordinate (u-axis3) at (20,0);
\coordinate (t-axis3) at (24,6);
\node at (u-axis3) [anchor=north]{u};
\node at (t-axis3) [anchor=east]{t};
\draw [thick,->] (origin3) -- (u-axis3);
\draw [thick,->] (origin3) -- (t-axis3);
% More Coordinates
\coordinate (uL1) at (1,0);
\coordinate (deltaL1) at (1,2);
\coordinate (deltaM1) at (4,2);
\coordinate (deltaR1) at (7,2);
\coordinate (TM1) at (4,5);
\coordinate (TR1) at (7,5);
\coordinate (singular2) at (12.5,0);
\coordinate (delta2) at (10,2);
\coordinate (deltaL2) at (11.5,2);
\coordinate (deltaM2) at (13.5,2);
\coordinate (deltaR2) at (15.5,2);
\coordinate (tL2) at (11.3,3);
\coordinate (tR2) at (14.8,3);
\coordinate (TM2) at (15.5,5);
\coordinate (TR2) at (17.5,5);
\coordinate (deltaL3) at (21,2);
\coordinate (deltaM3) at (24,2);
\coordinate (uL3) at (21,0);
\coordinate (tL3) at (21,3);
\coordinate (tM3) at (24,3);
% Figure1
\draw [semithick] (deltaL1) -- (deltaR1) -- (TR1) -- (TM1);
\draw [dotted,semithick] (deltaL1) -- (uL1);
\fill [gray,opacity=0.5] (deltaM1) -- (deltaR1) -- (TR1) -- (TM1) -- cycle;
% Figure2
\draw [semithick] (singular2) .. controls +(0.5,1.2) and +(-0.5,-0.8) .. (deltaM2) .. controls +(0.5,0.8) and +(-0.5,-1) .. (TM2) node [anchor=south]{$C_0$};
\draw [dotted,semithick] (delta2) -- +(7,0);
\draw [semithick] (deltaL2) -- (deltaR2) .. controls +(0.5,0.7) and +(-0.5,-1.1) .. (TR2) -- (TM2);
\fill [gray,opacity=0.5] (deltaR2) .. controls +(0.5,0.7) and +(-0.5,-1.1) .. (TR2) -- (TM2) .. controls +(-0.5,-1) and +(0.5,0.8) ..(deltaM2) -- cycle;
\draw [dotted,semithick] (deltaR2) .. controls +(-0.2,0.2) and +(0.15,-0.2) .. (tR2) -- (tL2) .. controls +(0.05,-0.4) and +(-0.1,0.15) .. (deltaL2);
\fill [gray,opacity=0.5] (deltaR2) .. controls +(-0.2,0.2) and +(0.15,-0.2) .. (tR2) -- (tL2) .. controls +(0.05,-0.4) and +(-0.1,0.15) .. (deltaL2) -- cycle;
% Figure3
\draw [semithick] (deltaL3) -- (deltaM3);
\draw [dotted,semithick] (tM3) -- (tL3) -- (uL3);
\fill [gray,opacity=0.5] (deltaL3) -- (deltaM3) -- (tM3) -- (tL3) -- cycle;
% Text
\node at (uL1) [anchor=north]{$u^*$};
\node at (origin1) [anchor=north]{0};
\node at (deltaM1) [anchor=north east]{$\delta$};
\node at (2.5,1.7) {Data};
\node at (5.5,3.5) {Taylor};
\node at (delta2) [anchor=east]{$\delta$};
\node at (15.5,3.5) {Taylor};
\node at (13,2.5) {LWP};
\node at (uL3) [anchor=north]{$u^*$};
\node at (deltaM3) [anchor=north east]{$\delta$};
\node at (22.5,2.5) {LWP};
\draw [semithick,->] (8,2.5) -- (9,2.5);
\node at (8.5,3) {to $(t,r)$};
\draw [semithick,->] (18,2.5) -- (20,2.5);
\node at (19,3) {back to $(t,u)$};

\end{tikzpicture}}
\caption{Local Well-Posedness}
\label{fig:local_well_posedness}
\end{figure}
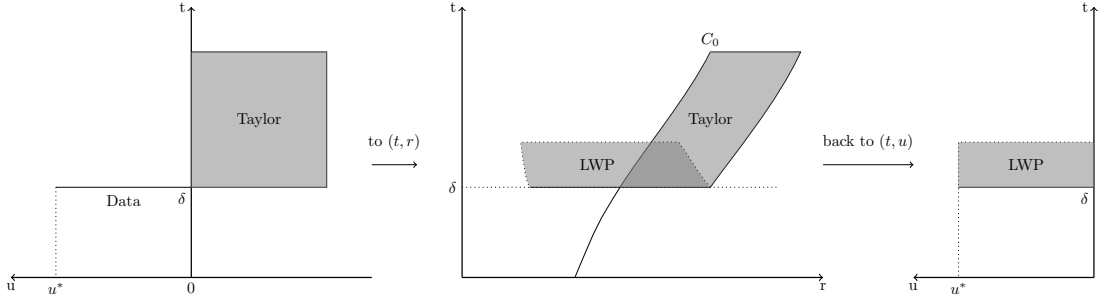

This is a local well-posedness result for a Cauchy-characteristic mixed initial data problem. We prove it by adopting Rendall's method to reduce this to a standard local well-posedness result for a Cauchy initial data problem for a quasilinear wave equation. See figure \ref{fig:local_well_posedness}.

We first construct an approximate solution to the Euler equations on the other side of $C_0$ by the Taylor expansion similar to \ref{eq:initial_Taylor}.
\begin{equation*}
\begin{cases}
    w_{C_0, N}(t, u) = w_{C_0}(t, 0) + \sum\limits_{n = 1}^N \frac{T^n w_{C_0}(t, 0)}{n!} u^n \\
    \underline{w}_{C_0, N}(t, u) = \underline{w}_{C_0}(t, 0) + \sum\limits_{n = 1}^N \frac{T^n \underline{w}_{C_0}(t, 0)}{n!} u^n \\
    r_{C_0, N}(t, u) = r_{C_0}(t, 0) - \sum\limits_{n = 0}^{N - 1} \frac{T^n \kappa_{C_0}(t, 0)}{n!} u^n, \qquad u \in [-u_*, 0] \\
\end{cases}
\end{equation*}

The range $u_*$ is chosen to ensure that the following bounds hold.
\begin{equation*}
    \min\limits_{[\delta, t_0] \times [- u_*, 0]} c_{C_0, N}(t,u) \geq \frac{1}{2} \min\limits_{t \in [\delta, t_0]} c_{C_0}(t,0) > 0, \qquad \min\limits_{[\delta, t_0] \times [- u_*, 0]} \kappa_{C_0, N}(t,u) \geq \frac{1}{2} \min\limits_{t \in [\delta, t_0]} \kappa_{C_0}(t,0) > 0
\end{equation*}

Now we transform to the $(t, x)$ coordinates.
\begin{equation*}
\begin{aligned}
& \begin{cases}
    (w_{\Sigma_\delta}, \underline{w}_{\Sigma_\delta}, r_{\Sigma_\delta}) \in H_u^{N + 1}([0, u^*]) \\
    ((L w)_{\Sigma_\delta}, (L \underline{w})_{\Sigma_\delta}) \in H_u^N([0, u^*]) \\
    (w_{C_0, N}, \underline{w}_{C_0, N}, r_{C_0, N}) \in C_{t, u}^\infty ([\delta, t_0] \times [- u_*, 0]) \\
\end{cases} \\
\Longrightarrow & \begin{cases}
    (w_{\Sigma_\delta}, \underline{w}_{\Sigma_\delta}) \in H_r^{N + 1}([r_{\Sigma_\delta}(u^*), r_{\Sigma_\delta}(0)]) \\
    ((\partial_t + (v_{\Sigma_\delta} + c_{\Sigma_\delta}) \partial_r) w_{\Sigma_\delta}, (\partial_t + (v_{\Sigma_\delta} + c_{\Sigma_\delta}) \partial_r) \underline{w}_{\Sigma_\delta}) \in H_r^N([r_{\Sigma_\delta}(u^*), r_{\Sigma_\delta}(0)]) \\
    (\partial_t w_{\Sigma_\delta}, \partial_t \underline{w}_{\Sigma_\delta}) \in H_r^N([r_{\Sigma_\delta}(u^*), r_{\Sigma_\delta}(0)]) \\
    (w_{C_0, N}, \underline{w}_{C_0, N}) \in C_{t, r}^\infty (\{t \in [\delta, t_0], r \in [r_{C_0}(t), r_{C_0, N}(t, - u_*)]\})
\end{cases} \\
\end{aligned}
\end{equation*}

Now we analyze the Euler equations and these functions. We copy the Euler equations in the radial coordinates here.
\begin{equation*}
\begin{cases}
    \partial_t w + (v - c) \partial_r w = - \frac{c v}{r}, \\
    \partial_t \underline{w} + (v + c) \partial_r \underline{w} = - \frac{c v}{r}
\end{cases}
\end{equation*}

By the definition \ref{eq:initial_Taylor} of $((L w)_{\Sigma_\delta}, (L \underline{w})_{\Sigma_\delta})$ we know that $(w_{\Sigma_\delta}, \underline{w}_{\Sigma_\delta}, \partial_t w_{\Sigma_\delta}, \partial_t \underline{w}_{\Sigma_\delta})$ satisfy the Euler equations.

$(w_{C_0, N}, \underline{w}_{C_0, N})$ doesn't satisfy the Euler equations, but the error vanishes for up to $N - 1$ derivatives on $C_0$.

Putting the two parts together, we obtain a $H_r^{N + 1} \times H_r^N$ initial data on $\Sigma_\delta$ and we consider the quasilinear wave equations obtained by differentiating the Euler equations with extra error term on the right hand side coming from $(w_{C_0, N}, \underline{w}_{C_0, N})$.

The standard $H^s$ well-posedness result for the Cauchy initial value problem of quasilinear wave equations now applies. Transforming back to $(t, u)$ coordinates, we obtain the local well-posedness result we need.

By the local well-posedness theorem, we have a solution $\Psi_{\delta, N} = (w, \underline{w}, r)_{\delta, N}$ defined on $(t, u) \in [\delta, t^*(\delta, u^*, N)] \times [0, u^*]$. The life span $t^*(\delta, u^*, N)$ of the solution guaranteed by the local well-posedness theorem relies on $\delta, u^*, N$.

\subsection{Energy Estimates Uniform in $\delta$}

Now we need to set up energy estimates in the region $(t, u) \in [\delta, t^*(\delta, u^*, N)] \times [0, u^*]$ for $\Psi_{\delta, N}$, which will be used to show that if $u^* \leq u_0^*$ is small enough (depending on the initial and boundary data) and $\delta \leq \delta_0(N)$ is small enough (depending on the data and N), then the solutions can be extended to a uniform time $t^*$ for all $(\delta, u^*, N)$ satisfying these smallness conditions. Then we can get a sequence of uniformly bounded approximate solutions on $[\delta, t^*] \times [0, u^*]$, which will lead to the exact solution by convergence.

The multiplier is chosen to be $X = \frac{\kappa}{c} L + \underline{L}$, that is, $a = \frac{\kappa}{c}, b = 1$ in the energy estimate framework \ref{eq:energy_estimate}.

The commutator is $T$.

In this case, the energy estimate takes the following form.
\begin{equation} \label{eq:energy_estimate}
\begin{cases}
    \Box_g \psi = \tilde{\rho}_\psi \\
    \mathcal{E}_\psi^{[0, u^*]}(t^*) - \mathcal{E}_\psi^{[0, u^*]}(\delta) + \mathcal{F}_\psi^{[\delta, t^*]}(u^*) - \mathcal{F}_\psi^{[\delta, t^*]}(0) = \int_{\{\delta \leq t \leq t^*, 0 \leq u \leq u^*\}} Q_\psi d\mu_g \\
    \mathcal{E}_\psi^{[0, u^*]}(t) = \frac{1}{2} \int_{[0, u^*]} \int_{S_{t, u}} |\underline{L}\psi|^2 + \frac{\kappa^2}{c^2} |L\psi|^2 d\mu_\slashed{g} du \\
    \mathcal{F}_\psi^{[\delta, t^*]}(u) = \int_{[\delta, t^*]} \int_{S_{t, u}} \frac{\kappa}{c} |L\psi|^2 d\mu_\slashed{g} dt \\
    Q_\psi = Q_{0, \psi} + Q_{1, \psi} + Q_{2, \psi} \\
    Q_{0, \psi} = - \tilde{\rho}_\psi (\frac{\kappa}{c} L\psi + \underline{L}\psi) \\
    Q_{1, \psi} = \frac{1}{\mu} T(\frac{\kappa}{c}) |L\psi|^2 \\
    Q_{2, \psi} = - \frac{1}{\mu} \frac{2 \kappa v}{c r} L\psi \underline{L}\psi
\end{cases}
\end{equation}

We now prove a general energy estimate. During the proof, $N$ stands for the number of times we apply the commutator instead of the $N$-th order approximate data.

The estimate contains two parts. One is the basic energy estimate of the $N = 3$ case, which determines the smallness we mentioned above and gives the uniform $t^*$. Another is a recursive process that increases the regularity of the solution, provided that the initial data are regular enough (that is, taking enough order of approximation).

We first list the assumptions we need in both parts.
\begin{equation} \label{eq:localness}
\begin{cases}
\begin{aligned} 
    & t^* \leq 1 && |c - c_0| \leq \frac{1}{2} c_0 \\
    & u^* \leq 1 && |v - v_0| \leq \frac{1}{2} c_0 \\
    &&& |r - r_0| \leq \frac{1}{2} r_0
\end{aligned}
\end{cases}
\end{equation}

Here, $(c_0, v_0, r_0)$ are the values of the background solution at the initial singular point $(t = 0, u = 0)$.

These assumptions say that we are dealing with local solutions, staying away from the vacuum and the center.

We also rephrase the vanishing orders \ref{eq:local-vanishing_C0} into quantitative bounds for the smooth data on $C_0$.
\begin{equation} \label{eq:C_0_bounds}
\begin{cases}
\begin{aligned}
    & |L w| \leq M_0 && |L \underline{w}| \leq M_0 && |L \kappa - 1| \leq M_0 t \\
    & |T w| \leq M_0 t && |T \underline{w} + \frac{2}{\gamma + 1}| \leq M_0 t && |\kappa - t| \leq M_0 t^2 \\
    & |L T^n w| \leq M_n && |L T^n \underline{w}| \leq M_n && |L T^n \kappa| \leq M_n t \\
    & |T^{n + 1} w| \leq M_n t && |T^{n + 1} \underline{w}| \leq M_n t && |T^n \kappa| \leq M_n t^2, \quad n \geq 1 \\
\end{aligned}
\end{cases}
\end{equation}

\begin{remark}
The $(t^*, u^*)$ bounds listed in \ref{eq:localness} are loose. We will put more and more restrictions on them during the proofs.
\end{remark}

Then we state the two propositions.

\begin{proposition}[Base Energy Estimates] \label{prop:local-base_energy}
Suppose that we have a regular solution $(w, \underline{w}, r)$ of the Euler equations defined on $(t, u) \in [\delta, t^*] \times [0, u^*]$. They coincide with the data on the boundary $C_0$. They belong to the following regularity classes.
\begin{equation} \label{eq:local-base-regularity_class}
\begin{cases}
    w \in C_t^0 H_u^{4} \cap C_t^1 H_u^3 \\
    \underline{w} \in C_t^0 H_u^4 \cap C_t^1 H_u^3 \\
    r \in C_t^1 H_u^4
\end{cases}
\end{equation}

And they have the following bounds for the initial data.
\begin{equation} \label{eq:local-base-initial_bounds}
\begin{cases}
    \mathcal{E}_{T^3 w}^{[0, u^*]}(\delta) + \mathcal{E}_{T^3 \underline{w}}^{[0, u^*]}(\delta) \leq E_{\mathrm{Data}, 3} \delta^2 \\
    \int_{[0, u^*]} \int_{S_{\delta, u}} |T^3 \kappa|^2 d\mu_{\slashed{g}} du \leq E_{\mathrm{Data}, 3} \delta^4 \\
    \mathcal{F}_{T^3 w}^{[\delta, t^*]}(0) + \mathcal{F}_{T^3 \underline{w}}^{[\delta, t^*]}(0) \leq F_{\mathrm{Data}, 3} {t^*}^2 \\
    \begin{aligned}
        & |c(\delta, u) - c_0| \leq \frac{1}{4} c_0 && |v(\delta, u) - v_0| \leq \frac{1}{4} c_0 && |r(\delta, u) - r_0| \leq \frac{1}{4} r_0 \\
        & |L \kappa(\delta, u) - 1| \leq \frac{1}{4} && |\kappa(\delta, u) - \delta| \leq \frac{1}{4} \delta && \\
    \end{aligned}
\end{cases}
\end{equation}

Suppose that we have the following smallness condition on $(t^*, u^*)$. There exist constants $\mathrm{C}_1(c_0, v_0, r_0), \mathrm{C}_2(c_0, v_0, r_0)$ such that we have:
\begin{equation} \label{eq:local-base-t_u_smallness}
    t^* \leq \frac{1}{\mathrm{C}_1 \max\{M_0, M_1, M_2, \sqrt{E_{\mathrm{Data}, 3}}, \sqrt{F_{\mathrm{Data}, 3}}\}}, \quad u^* \leq \frac{1}{\mathrm{C}_2}
\end{equation}

Then we have a uniform bound of the energy and flux in the region, as well as the uniform bound for the $L^\infty$-norm of the lower order quantities. ($n = 1, 2$)
\begin{equation} \label{eq:local-base-final_estimates}
\begin{cases}
    \mathcal{E}_{T^3 w}^{[0, u^*]}(t) + \mathcal{E}_{T^3 \underline{w}}^{[0, u^*]}(t) \leq E_3' t^2 \\
    \max\limits_{u \in [0, u^*]} \mathcal{F}_{T^3 w}^{[\delta, t]}(u) + \mathcal{F}_{T^3 \underline{w}}^{[\delta, t]}(u) \leq E_3' t^2 \\
    \int_{[0, u^*]} \int_{S_{t, u}} |T^3 \kappa|^2 d\mu_{\slashed{g}} du \leq E_3' t^4 \\
    \begin{aligned}
        & |c - c_0| \leq \frac{1}{2} c_0 && |v - v_0| \leq \frac{1}{2} v_0 && |r - r_0| \leq \frac{1}{2} r_0 \\
        & |L w| \leq M_0' && |L \underline{w}| \leq M_0' && |L \kappa - 1| \leq M_0' t \leq \frac{1}{2} \\
        & |T w| \leq M_0' t && |T \underline{w} + \frac{2}{\gamma + 1}| \leq M_0' t && |\kappa - t| \leq M_0' t^2 \leq \frac{1}{2} t \\
        & |L T^n w| \leq M_n' && |L T^n \underline{w}| \leq M_n' && |L T^n \kappa| \leq M_n' t \\
        & |T^{n + 1} w| \leq M_n' t && |T^{n + 1} \underline{w}| \leq M_n' t && |T^n \kappa| \leq M_n' t^2 \\
    \end{aligned}
\end{cases}
\end{equation}

$E_3', M_0', M_1', M_2'$ depend on $c_0, v_0, r_0, M_0, M_1, M_2, \sqrt{E_{\mathrm{Data}, 3}}, \sqrt{F_{\mathrm{Data}, 3}}$.
\end{proposition}

Then there is the propagation of regularity.

\begin{proposition}[Propagation of Regularity] \label{prop:local-regularity}
Let $N \geq 4$. For the solution in proposition \ref{prop:local-base_energy}, assume that the initial data on $\Sigma_\delta$ have higher regularity and is compatible with the data on $C_0$, satisfying the following bounds.
\begin{equation} \label{eq:local-high_initial_bounds}
\begin{cases}
    \mathcal{E}_{T^N w}^{[0, u^*]}(\delta) + \mathcal{E}_{T^N \underline{w}}^{[0, u^*]}(\delta) \leq E_{\mathrm{Data}, N} \delta^2 \\
    \int_{[0, u^*]} \int_{S_{\delta, u}} |T^N \kappa|^2 d\mu_{\slashed{g}} du \leq E_{\mathrm{Data}, N} \delta^4 \\
    \mathcal{F}_{T^N w}^{[\delta, t^*]}(0) + \mathcal{F}_{T^N \underline{w}}^{[\delta, t^*]}(0) \leq F_{\mathrm{Data}, N} {t^*}^2
\end{cases}
\end{equation}

Also, assume that we have already obtained the $L^2$ and $L^\infty$-bound of lower order terms in the whole region $(t, u) \in [\delta, t^*] \times [0, u^*]$. ($1 \leq n \leq N - 2$)
\begin{equation} \label{eq:local-high-lower_order_bounds}
\begin{cases}
    \mathcal{E}_{T^{N - 1} w}^{[0, u^*]}(t) + \mathcal{E}_{T^{N - 1} \underline{w}}^{[0, u^*]}(t) \leq E_{N - 1}' t^2 \\
    \max\limits_{u \in [0, u^*]} \mathcal{F}_{T^{N - 1} w}^{[\delta, t]}(u) + \mathcal{F}_{T^{N - 1} \underline{w}}^{[\delta, t]}(u) \leq E_{N - 1}' t^2 \\
    \int_{[0, u^*]} \int_{S_{t, u}} |T^{N - 1} \kappa|^2 d\mu_{\slashed{g}} du \leq E_{N - 1}' t^4 \\
    \begin{aligned}
        & |L w| \leq M_0' && |L \underline{w}| \leq M_0' && |L \kappa - 1| \leq M_0' t \\
        & |T w| \leq M_0' t && |T \underline{w} + \frac{2}{\gamma + 1}| \leq M_0' t && |\kappa - t| \leq M_0' t^2 \\
        & |L T^n w| \leq M_n' && |L T^n \underline{w}| \leq M_n' && |L T^n \kappa| \leq M_n' t \\
        & |T^{n + 1} w| \leq M_n' t && |T^{n + 1} \underline{w}| \leq M_n' t && |T^n \kappa| \leq M_n' t^2 \\
    \end{aligned}
\end{cases}
\end{equation}

Then we can obtain the higher order $L^2$ and $L^\infty$ bounds.
\begin{equation} \label{eq:local-high-final_estimates}
\begin{cases}
    \mathcal{E}_{T^N w}^{[0, u^*]}(t) + \mathcal{E}_{T^N \underline{w}}^{[0, u^*]}(t) \leq E_N' t^2 \\
    \max\limits_{u \in [0, u^*]} \mathcal{F}_{T^N w}^{[\delta, t]}(u) + \mathcal{F}_{T^N \underline{w}}^{[\delta, t]}(u) \leq E_N' t^2 \\
    \int_{[0, u^*]} \int_{S_{t, u}} |T^N \kappa|^2 d\mu_{\slashed{g}} du \leq E_N' t^4 \\
    \begin{aligned}
        & |L T^{N - 1} w| \leq M_{N - 1}' && |L T^{N - 1} \underline{w}| \leq M_{N - 1}' && |L T^{N - 1} \kappa| \leq M_{N - 1}' t \\
        & |T^N w| \leq M_{N - 1}' t && |T^N \underline{w}| \leq M_{N - 1}' t && |T^{N - 1} \kappa| \leq M_{N - 1}' t^2 \\
    \end{aligned}
\end{cases}
\end{equation}

The new bounds $E_N', M_{N - 1}'$ depend on $E_{\mathrm{Data}, N}, F_{\mathrm{Data}, N}, E_{N - 1}', M_n', (n \leq N - 2)$.
\end{proposition}

In the following proof, we will frequently use the less and similar notation $\lesssim$, which means $\leq$ up to a constant depending on $(c_0, v_0, r_0, N)$.

\subsection{Proof of Proposition \ref{prop:local-base_energy}} \label{subsec:local-proof-local_base_energy}

The proof is done by a continuity argument.

Let $M_{\leq 2} = \max\{M_0, M_1, M_2\}$.

First, from the bounds \ref{eq:local-base-initial_bounds} for the initial data and the bounds \ref{eq:C_0_bounds} for the data on $C_0$, we have the following $L^\infty$-norm bound by integration from the boundary $C_0$.
\begin{equation*}
\begin{aligned}
    |T^2 \kappa(\delta, u) - T^2 \kappa(\delta, 0)| & \leq \int_{[0, u]} |T^3 \kappa| du \leq \sqrt{\int_{[0, u]} |T^3 \kappa|^2 du} \sqrt{\int_{[0, u]} du} \\
    & \lesssim \sqrt{E_{\mathrm{Data}, 3} \delta^4} \sqrt{u^*} \leq \sqrt{E_{\mathrm{Data}, 3}} \delta^2 \\
\end{aligned}
\end{equation*}

\begin{equation*}
\begin{aligned}
    & |T^2 \kappa(\delta, u)| \lesssim (\sqrt{E_{\mathrm{Data}, 3}} + M_{\leq 2}) \delta^2 \\
    & |T \kappa(\delta, u)| \leq |T \kappa(\delta, 0)| + \int_{[0, u]} |T^2 \kappa| du \lesssim (\sqrt{E_{\mathrm{Data}, 3}} + M_{\leq 2}) \delta^2 \\
    & |\kappa(\delta, u) - \kappa(\delta, 0)| \leq \int_{[0, u]} |T \kappa| du \lesssim (\sqrt{E_{\mathrm{Data}, 3}} + M_{\leq 2}) \delta^2 \\
    & |\kappa(\delta, u) - \delta| \lesssim (\sqrt{E_{\mathrm{Data}, 3}} + M_{\leq 2}) \delta^2 \\
\end{aligned}
\end{equation*}

Similarly, we obtain the bounds for $w, \underline{w}$.
\begin{equation*}
\begin{aligned}
    & |L T^2 w(\delta, u)| \leq |L T^2 w(\delta, 0)| + \int_{[0, u]} |L T^3 w| du \lesssim (\sqrt{E_{\mathrm{Data}, 3}} + M_{\leq 2}) \\
    & |T^3 w(\delta, u)| \leq |T^3 w(\delta, 0)| + \int_{[0, u]} |T^4 w| du \lesssim (\sqrt{E_{\mathrm{Data}, 3}} + M_{\leq 2}) \delta \\
\end{aligned}
\end{equation*}

Here we have used the fact that $\kappa \sim t$.

Performing the same arguments repeatedly, we obtain the $L^\infty$-bounds for the quantities on $\Sigma_\delta$. ($n = 1, 2$)
\begin{equation*}
\begin{cases}
\begin{aligned}
    & |L T^n w(\delta, u)| \lesssim (\sqrt{E_{\mathrm{Data}, 3}} + M_{\leq 2}) && |L T^n \underline{w}(\delta, u)| \lesssim (\sqrt{E_{\mathrm{Data}, 3}} + M_{\leq 2}) \\
    & |T^{n + 1} w(\delta, u)| \lesssim (\sqrt{E_{\mathrm{Data}, 3}} + M_{\leq 2}) \delta && |T^{n + 1} \underline{w}(\delta, u)| \lesssim (\sqrt{E_{\mathrm{Data}, 3}} + M_{\leq 2}) \delta \\
    & |L w(\delta, u)| \lesssim (\sqrt{E_{\mathrm{Data}, 3}} + M_{\leq 2}) && |L \underline{w}(\delta, u)| \lesssim (\sqrt{E_{\mathrm{Data}, 3}} + M_{\leq 2}) \\
    & |T w(\delta, u)| \lesssim (\sqrt{E_{\mathrm{Data}, 3}} + M_{\leq 2}) \delta && |T \underline{w}(\delta, u) + \frac{2}{\gamma - 1}| \lesssim (\sqrt{E_{\mathrm{Data}, 3}} + M_{\leq 2}) \delta \\
\end{aligned}
\end{cases}
\end{equation*}

We introduce a new notation to represent $\sqrt{E_{\mathrm{Data}, 3}} + \sqrt{F_{\mathrm{Data}, 3}} + M_{\leq 2}$, since it will frequently appear later. We include the flux because it is tied together with the energy in the estimates and will contribute to the bounds we obtain.
\begin{equation*}
    M = \sqrt{E_{\mathrm{Data}, 3}} + \sqrt{F_{\mathrm{Data}, 3}} + M_{\leq 2}
\end{equation*}

Taking the propagation equation for $\kappa$ into consideration, we notice that the quantities we need to control are $L T^2 w, T^3 w, L T^2 \underline{w}, T^3 \underline{w}$. Other bounds can be obtained from these. So we write the bootstrap assumptions as follows.
\begin{equation} \label{eq:local-base-bootstrap}
\begin{cases}
    |L T^2 w(t, u)| \leq \mathrm{C} M \\
    |L T^2 \underline{w}(t, u)| \leq \mathrm{C} M \\
    |T^3 w(t, u)| \leq \mathrm{C} M t \\
    |T^3 \underline{w}(t, u)| \leq \mathrm{C} M t, \quad \forall (t, u) \in [\delta, t^*] \times [0, u^*] \\
\end{cases}
\end{equation}

$ \mathrm{C} =  \mathrm{C}(c_0, r_0)$ is a large constant.

Then by repeating the auguments above and using the propagation equation for $\kappa$, we immediately obtain all $L^\infty$-bounds that we need. ($n = 1, 2$)
\begin{equation} \label{eq:local-base-Linfty}
\begin{cases}
\begin{aligned}
    & |L w| \lesssim M && |L \underline{w}| \lesssim M && |L \kappa - 1| \lesssim M t \\
    & |T w| \lesssim M t && |T \underline{w} + \frac{2}{\gamma - 1}| \lesssim M t && |\kappa - t| \lesssim M t^2 \\
    & |L T^n w| \lesssim M && |L T^n \underline{w}| \lesssim M && |L T^n \kappa| \lesssim M t \\
    & |T^{n + 1} w| \lesssim M t && |T^{n + 1} \underline{w}| \lesssim M t && |T^n \kappa| \lesssim M t^2 \\
\end{aligned}
\end{cases}
\end{equation}

We place a smallness requirement on $t^*$.
\begin{equation} \label{eq:local-base-smallness_sim}
    t^* \leq \frac{1}{\mathrm{C}_1' M}
\end{equation}

This is to ensure that the $L^\infty$-bounds above imply the following $\sim$ relations.
\begin{equation} \label{eq:local-sim_relations}
\begin{cases}
\begin{aligned}
    & |c(t, u) - c_0| \leq \frac{1}{2} c_0 && |v(t, u) - v_0| \leq \frac{1}{2} c_0 && |r(t, u) - r_0| \leq \frac{1}{2} r_0 \\
    & |L \kappa(t, u) - 1| \leq \frac{1}{2} && |\kappa(t, u) - t| \leq \frac{1}{2} t && \\
\end{aligned}
\end{cases}
\end{equation}

We now use the energy estimate to improve these assumptions.

Recall the expression of the wave operator in terms of $L, \underline{L}$.
\begin{equation*}
\begin{aligned}
    \mu \Box_g f & = - \frac{1}{2} \mathrm{tr} \underline{\chi} (Lf) - \frac{1}{2} \mathrm{tr} \chi (\underline{L}f) - L(\underline{L}(f)) \\
    & = - \frac{1}{2} \mathrm{tr} \underline{\chi} (Lf) - \frac{1}{2} \mathrm{tr} \chi (\underline{L}f) - \underline{L}(L(f)) - L(\frac{\kappa}{c}) (Lf) \\
    & = - \frac{2 \kappa v}{c r} Lf - 2 \frac{v + c}{r} Tf - \frac{\kappa}{c} LLf - 2 LTf - L(\frac{\kappa}{c}) Lf \\
\end{aligned}
\end{equation*}

Substitute the first order Euler equations for $w, \underline{w}$ into the expression to get the wave equations for $w, \underline{w}$.
\begin{equation} \label{eq:wave_w_wbar}
\begin{aligned}
    \mu \Box_g w & = - \frac{\kappa}{c} \frac{(v - c)}{r} Lw - \frac{v + c}{r} \underline{L}w - L\underline{L}w \\
    & = - \frac{v - c}{r} (\underline{L} - 2T)w + \frac{(v + c)}{r} \frac{\kappa v}{r} + L(\frac{\kappa v}{r}) \\
    & = - \frac{v - c}{r} (- \frac{\kappa v}{r} - 2Tw) + \frac{(v + c)}{r} \frac{\kappa v}{r} + L(\frac{\kappa v}{r}) \\
    & = 2 \frac{\kappa v^2}{r^2} + 2 \frac{v - c}{r} Tw + L(\frac{\kappa v}{r}) \\
    \mu \Box_g \underline{w} & = - \frac{\kappa}{c} \frac{(v - c)}{r} L\underline{w} - \frac{v + c}{r} \underline{L}\underline{w} - \underline{L}L\underline{w} - L(\frac{\kappa}{c}) L\underline{w} \\
    & = \frac{\kappa}{c} \frac{(v - c)}{r} \frac{c v}{r} + \frac{\kappa}{c} \frac{(v + c)}{r} \frac{c v}{r} - 2 \frac{v + c}{r} T\underline{w} + \underline{L}(\frac{c v}{r}) + L(\frac{\kappa}{c}) \frac{c v}{r} \\
    & = 2 \frac{\kappa v^2}{r^2} - 2 \frac{v + c}{r} T\underline{w} + \frac{\kappa}{c} L(\frac{c v}{r}) + 2 T(\frac{c v}{r}) + L(\frac{\kappa}{c}) \frac{c v}{r} \\
    & = 2 \frac{\kappa v^2}{r^2} - 2 \frac{v + c}{r} T\underline{w} + L(\frac{\kappa v}{r}) + 2T(\frac{c v}{r}) \\
\end{aligned}
\end{equation}

We need to commute $T$, so we also write down the expression for the commutator $[\mu \Box_g, T]$.
\begin{equation} \label{eq:commutator}
\begin{aligned}
    [\mu \Box_g, T] f & = \mu \Box_g T f - T (\mu \Box_g f) \\
    & = 2 T (\frac{\kappa v}{c r}) L f + 2 T (\frac{v + c}{r}) T f + T (\frac{\kappa}{c}) L L f + L T (\frac{\kappa}{c}) L f \\
\end{aligned}
\end{equation}

Now we begin to control the integral $\int_{[\delta, t^*] \times [0, u^*]} Q_\psi d\mu_g$ by energy, flux and the $L^\infty$-bounds term by term.

First deal with $\int_{[\delta, t^*] \times [0, u^*]} Q_{1, \psi} + Q_{2, \psi} d\mu_g$ with $\psi = T^3 w, T^3 \underline{w}$.
\begin{equation*}
\begin{aligned}
    \int_{[\delta, t^*] \times [0, u^*]} Q_{1, \psi} + Q_{2, \psi} d\mu_g & = \int_\delta^{t^*} \int_0^{u^*} \int_{S_{t, u}} \mu (Q_{1, \psi} + Q_{2, \psi}) d\mu_{\slashed{g}} du dt \\
    & = \int_{t, u, S_{t, u}} T(\frac{\kappa}{c}) |L \psi|^2 - \frac{2 \kappa v}{c r} L \psi \underline{L} \psi d\mu_{\slashed{g}} du dt \\
    & = \int_{t, u, S_{t, u}} (\frac{T \kappa}{c} - \frac{\kappa T c}{c^2}) |L \psi|^2 - \frac{2 \kappa v}{c r} L \psi \underline{L} \psi d\mu_{\slashed{g}} du dt \\
    (\mathrm{take}\ |\cdot|) & \lesssim \int_{t, u, S_{t, u}} M t^2 |L \psi|^2 + t |L \psi| |\underline{L} \psi| d\mu_{\slashed{g}} du dt\\
    & \lesssim \int_\delta^{t^*} M \mathcal{E}_\psi^{[0, u^*]}(t) dt
\end{aligned}
\end{equation*}

Then we deal with the main part $\int_{[\delta, t^*] \times [0, u^*]} Q_{0, T^3 \psi} d\mu_g$ with $\psi = w, \underline{w}$. We separate it into the commutator part and the source part using \ref{eq:wave_w_wbar}, \ref{eq:commutator}.
\begin{equation} \label{eq:local-base-error_split_1}
\begin{aligned}
    & \int_{[\delta, t^*] \times [0, u^*]} Q_{0, T^3 \psi} d\mu_g \\
    = & \int_\delta^{t^*} \int_0^{u^*} \int_{S_{t, u}} - \mu \tilde{\rho}_{T^3 \psi} (\frac{\kappa}{c} L T^3 \psi + \underline{L} T^3 \psi) d\mu_\slashed{g} du dt \\
    = & \int_{t, u, S_{t, u}} - \mu \Box_g T^3 \psi (\frac{\kappa}{c} L T^3 \psi + \underline{L} T^3 \psi) d\mu_\slashed{g} du dt \\
    = & \int_{t, u, S_{t, u}} - (\underbrace{\sum\limits_{n_1 + n_2 = 2, n_i \geq 0} T^{n_1} [\mu \Box_g, T] T^{n_2} \psi}_{\mathrm{comm}} + \underbrace{T^3 (\mu \tilde{\rho}_\psi)}_{\mathrm{source}}) (\frac{\kappa}{c} L T^3 \psi + \underline{L} T^3 \psi) d\mu_\slashed{g} du dt \\
    = & I_{\mathrm{comm}, \psi} + I_{\mathrm{source}, \psi}
\end{aligned}
\end{equation}

\begin{equation} \label{eq:local-base-error_split_2}
\begin{aligned}
& \begin{aligned}
    I_{\mathrm{comm}, \psi} & = \int_{t, u, S_{t, u}} - (\frac{\kappa}{c} L T^3 \psi + \underline{L} T^3 \psi) \sum\limits_{n_1 + n_2 = 2,\ n_i \geq 0} T^{n_1} [\mu \Box_g, T] T^{n_2} \psi d\mu_\slashed{g} du dt \\
    & = \int_{t, u, S_{t, u}} - (\frac{\kappa}{c} L T^3 \psi + \underline{L} T^3 \psi) \sum\limits_{n_1 + n_2 = 2,\ n_i \geq 0} \mathrm{C}_{n_1, n_2} \\
    & \qquad (\underbrace{2 T^{n_1 + 1} (\frac{\kappa v}{c r}) L T^{n_2} \psi}_{1} + \underbrace{2 T^{n_1 + 1} (\frac{v + c}{r}) T^{n_2 + 1} \psi}_{2} \\
    & \qquad + \underbrace{T^{n_1 + 1} (\frac{\kappa}{c}) L L T^{n_2} \psi}_{3} + \underbrace{L T^{n_1 + 1} (\frac{\kappa}{c}) L T^{n_2} \psi}_{4}) d\mu_\slashed{g} du dt \\
    & = I_{\mathrm{comm}, 1, \psi} + I_{\mathrm{comm}, 2, \psi} + I_{\mathrm{comm}, 3, \psi} + I_{\mathrm{comm}, 4, \psi} \\
\end{aligned} \\
& \begin{aligned}
    I_{\mathrm{source}, \psi} = \int_{t, u, S_{t, u}} - T^3 (\mu \tilde{\rho}_\psi) (\frac{\kappa}{c} L T^3 \psi + \underline{L} T^3 \psi) d\mu_\slashed{g} du dt \\
\end{aligned} \\
\end{aligned}
\end{equation}

We now estimate these term by term. Handle the source term first.
\begin{equation*}
\begin{aligned}
    T^3 (\mu \tilde{\rho}_w) & = T^3 (2 \frac{\kappa v^2}{r^2} + 2 \frac{v - c}{r} Tw + L(\frac{\kappa v}{r})) \\
    (\mathrm{take}\ |\cdot|) & \lesssim (|T^3 \kappa| + M t^2) + (|T^4 w| + M t) + (|L T^3 \kappa| + M |T^3 \kappa| + t |L T^3 v| + M t) \\
    & \lesssim M |T^3 \kappa| + |T^4 w| + |T^4 \underline{w}| + t |LT^3 w| + t |L T^3 \underline{w}| + M t \\
\end{aligned}
\end{equation*}

\begin{equation*}
\begin{aligned}
    I_{\mathrm{source}, w} & = \int_{t, u, S_{t, u}} - T^3 (\mu \tilde{\rho}_w) (\frac{\kappa}{c} L T^3 w + \underline{L} T^3 w) d\mu_\slashed{g} du dt \\
    (\mathrm{take}\ |\cdot|) & \lesssim \int_{t, u, S_{t, u}} (M |T^3 \kappa| + |T^4 w| + |T^4 \underline{w}| + t |LT^3 w| + t |L T^3 \underline{w}| + M t) |\frac{\kappa}{c} L T^3 w + \underline{L} T^3 w| d\mu_\slashed{g} du dt \\
    & \lesssim \int_{t, u, S_{t, u}} M^2 |T^3 \kappa|^2 + M^2 t^2 d\mu_\slashed{g} du dt + \int_\delta^{t^*} \mathcal{E}_{T^3 w}^{[0, u^*]}(t) + \mathcal{E}_{T^3 \underline{w}}^{[0, u^*]}(t) dt \\
    & \lesssim M^2 {t^*}^3 + \int_{t, u, S_{t, u}} M^2 |T^3 \kappa|^2 d\mu_\slashed{g} du dt + \int_\delta^{t^*} \mathcal{E}_{T^3 w}^{[0, u^*]}(t) + \mathcal{E}_{T^3 \underline{w}}^{[0, u^*]}(t) dt \\
\end{aligned}
\end{equation*}

The case $I_{\mathrm{source}, \underline{w}}$ is quite similar.

\begin{equation*}
\begin{aligned}
    T^3 (\mu \tilde{\rho}_{\underline{w}}) & = T^3 (2 \frac{\kappa v^2}{r^2} - 2 \frac{v + c}{r} T\underline{w} + L(\frac{\kappa v}{r}) + 2T(\frac{c v}{r})) \\
    (\mathrm{take}\ |\cdot|) & \lesssim (|T^3 \kappa| + M t^2) + (|T^4 \underline{w}| + M t) \\
    & \qquad + (|L T^3 \kappa| + M |T^3 \kappa| + t |L T^3 v| + M t) \\
    & \qquad + (|T^4 c| + |T^4 v| + |T^3 \kappa| + M t) \\
    & \lesssim M |T^3 \kappa| + |T^4 w| + |T^4 \underline{w}| + t |LT^3 w| + t |L T^3 \underline{w}| + M t \\
\end{aligned}
\end{equation*}

\begin{equation*}
\begin{aligned}
    I_{\mathrm{source}, \underline{w}} & = \int_{t, u, S_{t, u}} - T^3 (\mu \tilde{\rho}_{\underline{w}}) (\frac{\kappa}{c} L T^3 \underline{w} + \underline{L} T^3 \underline{w}) d\mu_\slashed{g} du dt \\
    (\mathrm{take}\ |\cdot|) & \lesssim M^2 {t^*}^3 + \int_{t, u, S_{t, u}} M^2 |T^3 \kappa|^2 d\mu_\slashed{g} du dt + \int_\delta^{t^*} \mathcal{E}_{T^3 w}^{[0, u^*]}(t) + \mathcal{E}_{T^3 \underline{w}}^{[0, u^*]}(t) dt \\
\end{aligned}
\end{equation*}

Now we turn to the commutator terms.

\begin{equation*}
\begin{aligned}
    & |\sum\limits_{n_1 + n_2 = 2,\ n_i \geq 0} \mathrm{C}_{n_1, n_2} T^{n_1 + 1} (\frac{\kappa v}{c r}) L T^{n_2} \psi| \lesssim M |T^3 \kappa| + M t \\
    & |\sum\limits_{n_1 + n_2 = 2,\ n_i \geq 0} \mathrm{C}_{n_1, n_2} T^{n_1 + 1} (\frac{v + c}{r}) T^{n_2 + 1} \psi| \lesssim M t \\
    & \begin{aligned}
        |\sum\limits_{n_1 + n_2 = 2,\ n_i \geq 0} \mathrm{C}_{n_1, n_2} L T^{n_1 + 1} (\frac{\kappa}{c}) L T^{n_2} \psi| & \lesssim M |T^4 w| + M |T^4 \underline{w}| + M^2 |T^3 \kappa| \\
        & \qquad + M t |L T^3 w| + M t |L T^3 \underline{w}| + M \\
    \end{aligned}
\end{aligned}
\end{equation*}

So for $I_{\mathrm{comm}, 1, \psi}, I_{\mathrm{comm}, 2, \psi}, I_{\mathrm{comm}, 4, \psi}$, we only need to estimate $I_{\mathrm{comm}, 4, \psi}$ with this bound.

\begin{equation*}
\begin{aligned}
    & L L w = L (\frac{c}{\kappa} (\underline{L} - 2 T)) w = - 2 L (\frac{c}{\kappa}) T w - 2 \frac{c}{\kappa} L T w - L (\frac{c v}{r}) \\
    & L L \underline{w} = - L (\frac{c v}{r})
\end{aligned}
\end{equation*}

So for $I_{\mathrm{comm}, 3, \psi}$, we only deal with $I_{\mathrm{comm}, 3, w}$ since we will put absolute value on every single term.
\begin{equation*}
\begin{aligned}
    & |\sum\limits_{n_1 + n_2 = 2,\ n_i \geq 0} \mathrm{C}_{n_1, n_2} T^{n_1 + 1} (\frac{\kappa}{c}) L L T^{n_2} w| \\
    = & |\sum\limits_{n_1 + n_2 = 2,\ n_i \geq 0} \mathrm{C}_{n_1, n_2} T^{n_1 + 1} (\frac{\kappa}{c}) T^{n_2} (- 2 L (\frac{c}{\kappa}) T w - 2 \frac{c}{\kappa} L T w - L (\frac{c v}{r}))| \\
    \lesssim & \frac{M}{t} |T^3 \kappa| + |L T^3 w| + M \\
\end{aligned}
\end{equation*}

\begin{equation*}
\begin{aligned}
    & |\sum\limits_{\psi = w, \underline{w},\ 1 \leq i \leq 4} I_{\mathrm{comm}, i, \psi}| \\
    \lesssim & \sum\limits_{\psi = w, \underline{w}} \int_{t, u, S_{t, u}} |\frac{\kappa}{c} L T^3 \psi + \underline{L} T^3 \psi| \\
    & (\frac{M}{t} |T^3 \kappa| + |L T^3 w| + M + M |T^4 w| + M |T^4 \underline{w}| + M^2 |T^3 \kappa| + M t |L T^3 w| + M t |L T^3 \underline{w}| + M) d\mu_\slashed{g} du dt \\
    \lesssim & M \int_\delta^{t^*} \mathcal{E}_{T^3 w}^{[0, u^*]}(t) + \mathcal{E}_{T^3 \underline{w}}^{[0, u^*]}(t) dt \\
    & + \sum\limits_{\psi = w, \underline{w}} \int_{t, u, S_{t, u}} |\frac{\kappa}{c} L T^3 \psi + \underline{L} T^3 \psi| (\frac{M}{t} |T^3 \kappa| + |L T^3 w| + M) d\mu_\slashed{g} du dt \\
\end{aligned}
\end{equation*}

\begin{equation*}
\begin{aligned}
    & \sum\limits_{\psi = w, \underline{w}} \int_{t, u, S_{t, u}} |\frac{\kappa}{c} L T^3 \psi + \underline{L} T^3 \psi| (\frac{M}{t} |T^3 \kappa| + |L T^3 w| + M) d\mu_\slashed{g} du dt \\
    \lesssim & \sum\limits_{\psi = w, \underline{w}} \int_{t, u, S_{t, u}} \frac{M \epsilon}{t} |\frac{\kappa}{c} L T^3 \psi + \underline{L} T^3 \psi|^2 + (\frac{M}{\epsilon t} |T^3 \kappa|^2 + \frac{t}{M \epsilon} |L T^3 w|^2 + \frac{M t}{\epsilon}) d\mu_\slashed{g} du dt \\
    \lesssim & \int_\delta^{t^*} \frac{M \epsilon}{t} (\mathcal{E}_{T^3 w}^{[0, u^*]}(t) dt + \mathcal{E}_{T^3 \underline{w}}^{[0, u^*]}(t)) + \int_{t, u, S_{t, u}} \frac{M}{\epsilon t} |T^3 \kappa|^2 d\mu_\slashed{g} du dt + \frac{1}{M \epsilon} \int_0^{u^*} \mathcal{F}_{T^3 w}^{[\delta, t^*]}(u) du + \frac{M {t^*}^2}{\epsilon} \\
\end{aligned}
\end{equation*}

$\epsilon$ is a small constant to be determined.

Now we sum these up to get the following.
\begin{equation*}
\begin{aligned}
    & |\sum\limits_{\psi = T^3 w, T^3 \underline{w}} \mathcal{E}_\psi^{[0, u^*]}(t^*) - \mathcal{E}_\psi^{[0, u^*]}(\delta) + \mathcal{F}_\psi^{[\delta, t^*]}(u^*) - \mathcal{F}_\psi^{[\delta, t^*]}(0)| \\
    = & |\sum\limits_{\psi = T^3 w, T^3 \underline{w}} \int_{\{\delta \leq t \leq t^*, 0 \leq u \leq u^*\}} Q_\psi d\mu_g| \\
    \lesssim & \sum\limits_{\psi = T^3 w, T^3 \underline{w}} \int_\delta^{t^*} M \mathcal{E}_\psi^{[0, u^*]}(t) dt + M^2 {t^*}^3 + \int_{t, u, S_{t, u}} M^2 |T^3 \kappa|^2 d\mu_\slashed{g} du dt + \sum\limits_{\psi = T^3 w, T^3 \underline{w}} \int_\delta^{t^*} \mathcal{E}_\psi^{[0, u^*]}(t) dt \\
    & + \sum\limits_{\psi = T^3 w, T^3 \underline{w}} \int_\delta^{t^*} \frac{M \epsilon}{t} \mathcal{E}_\psi^{[0, u^*]}(t) dt + \int_{t, u, S_{t, u}} \frac{M}{\epsilon t} |T^3 \kappa|^2 d\mu_\slashed{g} du dt + \frac{1}{M \epsilon} \int_0^{u^*} \mathcal{F}_{T^3 w}^{[\delta, t^*]}(u) du + \frac{M {t^*}^2}{\epsilon} \\
    \lesssim & \sum\limits_{\psi = T^3 w, T^3 \underline{w}} \int_\delta^{t^*} (\frac{M \epsilon}{t} + M) \mathcal{E}_\psi^{[0, u^*]}(t) dt + \int_{t, u, S_{t, u}} \frac{M}{\epsilon t} |T^3 \kappa|^2 d\mu_\slashed{g} du dt + \frac{1}{M \epsilon} \int_0^{u^*} \mathcal{F}_{T^3 w}^{[\delta, t^*]}(u) du + \frac{M {t^*}^2}{\epsilon}
\end{aligned}
\end{equation*}

Now we handle the $\int |T^3 \kappa|^2$ term.
\begin{equation*}
\begin{aligned}
    & \int_{t, u, S_{t, u}} \frac{M}{\epsilon t} |T^3 \kappa|^2 d\mu_\slashed{g} du dt \\
    \lesssim & \int_\delta^{t^*} \int_0^{u^*} \frac{M}{\epsilon t} |T^3 \kappa|^2 du dt \\
    \lesssim & \int_0^{u^*} \int_\delta^{t^*} \frac{M}{\epsilon t} (|T^3 \kappa(\delta, u)| + \int_\delta^t |T^4 w(\tau, u)| + |T^4 \underline{w}(\tau, u)| d\tau)^2 dt du \\
    \lesssim & \int_\delta^{t^*} \frac{M}{\epsilon t} dt \int_0^{u^*} |T^3 \kappa(\delta, u)|^2 du + \int_0^{u^*} \int_\delta^{t^*} \frac{M}{\epsilon t} (\int_\delta^t |T^4 w(\tau, u)| + |T^4 \underline{w}(\tau, u)| d\tau)^2 dt du \\
\end{aligned}
\end{equation*}

\begin{equation*}
\begin{aligned}
    \int_\delta^{t^*} \frac{M}{\epsilon t} dt \int_0^{u^*} |T^3 \kappa(\delta, u)|^2 du & \lesssim \frac{M}{\epsilon} \ln(\frac{t^*}{\delta}) E_{\mathrm{Data}, 3} \delta^4 \\
    & \lesssim \frac{M^2 \delta^3}{\epsilon} \\
\end{aligned}
\end{equation*}

\begin{equation*}
\begin{aligned}
    & \int_0^{u^*} \int_\delta^{t^*} \frac{M}{\epsilon t} (\int_\delta^t |T^4 w(\tau, u)| + |T^4 \underline{w}(\tau, u)| d\tau)^2 dt du \\
    \lesssim & \int_0^{u^*} \int_\delta^{t^*} \frac{M}{\epsilon t} \int_\delta^t \tau d\tau \int_\delta^t \frac{|T^4 w(\tau, u)|^2 + |T^4 \underline{w}(\tau, u)|^2}{\tau} d\tau dt du \\
    \lesssim & \int_0^{u^*} \int_\delta^{t^*} \frac{M t}{\epsilon} \int_\delta^t \frac{|T^4 w(\tau, u)|^2 + |T^4 \underline{w}(\tau, u)|^2}{\tau} d\tau dt du \\
    = & \int_0^{u^*} \int_\delta^{t^*} \frac{|T^4 w(\tau, u)|^2 + |T^4 \underline{w}(\tau, u)|^2}{\tau} \int_\tau^{t^*} \frac{M t}{\epsilon} dt d\tau du \\
    \lesssim & \int_0^{u^*} \int_\delta^{t^*} \frac{M {t^*}^2}{\epsilon \tau} (|T^4 w(\tau, u)|^2 + |T^4 \underline{w}(\tau, u)|^2) d\tau du \\
    \lesssim & \sum\limits_{\psi = T^3 w, T^3 \underline{w}} \int_\delta^{t^*} \frac{M {t^*}^2}{\epsilon t} \mathcal{E}_\psi^{[0, u^*]}(t) dt
\end{aligned}
\end{equation*}

Substitute these into the inequality above to get the desired estimate.
\begin{equation*}
\begin{aligned}
    & |\sum\limits_{\psi = T^3 w, T^3 \underline{w}} \mathcal{E}_\psi^{[0, u^*]}(t^*) - \mathcal{E}_\psi^{[0, u^*]}(\delta) + \mathcal{F}_\psi^{[\delta, t^*]}(u^*) - \mathcal{F}_\psi^{[\delta, t^*]}(0)| \\
    \lesssim & \sum\limits_{\psi = T^3 w, T^3 \underline{w}} \int_\delta^{t^*} (\frac{M \epsilon}{t} + \frac{M {t^*}^2}{\epsilon t} + M) \mathcal{E}_\psi^{[0, u^*]}(t) dt + \frac{M {t^*}^2}{\epsilon} \\
    & + \frac{1}{M \epsilon} \int_0^{u^*} \mathcal{F}_{T^3 w}^{[\delta, t^*]}(u) du \\
\end{aligned}
\end{equation*}

Rewrite it as follows.
\begin{equation} \label{eq:local-base-energy_estimate}
\begin{aligned}
    & \sum\limits_{\psi = T^3 w, T^3 \underline{w}} \mathcal{E}_\psi^{[0, u^*]}(t^*) + \mathcal{F}_\psi^{[\delta, t^*]}(u^*) \\
    \leq & \mathrm{C}_{\mathrm{high}} (\sum\limits_{\psi = T^3 w, T^3 \underline{w}} \int_\delta^{t^*} (\frac{M \epsilon}{t} + \frac{M {t^*}^2}{\epsilon t} + M) \mathcal{E}_\psi^{[0, u^*]}(t) dt \\
    & + \frac{3^2}{M \epsilon} \int_0^{u^*} \mathcal{F}_{T^3 w}^{[\delta, t^*]}(u) du) + \mathrm{C}_{\mathrm{low}, 3} \frac{M {t^*}^2}{\epsilon} + \sum\limits_{\psi = T^3 w, T^3 \underline{w}} \mathcal{E}_\psi^{[0, u^*]}(\delta) + \mathcal{F}_\psi^{[\delta, t^*]}(0)
\end{aligned}
\end{equation}

\begin{remark} \label{rmk:high_low_constants}
We are distinguishing the hidden constants in front of the highest order terms and the lower order terms. This is for the proof of proposition \ref{prop:local-regularity}.

We will perform the Gronwall type arguments below, which will put restrictions on the size of $(t^*, u^*)$. We need to be clear about how these restrictions depend on $N$ and all the $L^\infty$-bounds.

We claim that $\mathrm{C}_{\mathrm{high}}$ will remain the same for all $N$. $3^2$ in front of the flux will become $N^2$. $\mathrm{C}_{\mathrm{low}, N}$ will depend on $N$ and the higher order $L^\infty$-bounds.

If a term is of the highest order, then all $T$ must act on that term, which means that its coefficient will be some fixed low order quantities regardless of how large $N$ is.

This is why we only put the indices $3$ in $\mathrm{C}_{\mathrm{low}, 3}$ and write $\mathrm{C}_{\mathrm{high}}$ without indices. However, there is one exception which comes from the term $\mathrm{C}_{0, 2} T (\frac{\kappa}{c}) L L T^2 w$. Here, $\mathrm{C}_{0, 2} = 3$ ($\mathrm{C}_{0, N} = N$ in general case). So we have a $3^2$ in front of the flux.

The $M$ in the coefficients of the energy and flux will also stay the same (will not change to higher order $L^\infty$-bounds) for arbitrary $N$.
\end{remark}

Now we choose $\epsilon$ and put further restrictions on $(t^*, u^*)$.
\begin{equation} \label{eq:local-base-smallness_Gronwall}
    \epsilon = \frac{1}{10 \mathrm{C}_{\mathrm{high}} M}, \quad t^* \leq \epsilon = \frac{1}{10 \mathrm{C}_{\mathrm{high}} M}, \quad u^* \leq \frac{1}{200 \mathrm{C}_{\mathrm{high}}^2}
\end{equation}

Then the estimate \ref{eq:local-base-energy_estimate} gives the following.
\begin{equation*}
\begin{aligned}
    \sum\limits_{\psi = T^3 w, T^3 \underline{w}} \mathcal{E}_\psi^{[0, u^*]}(t^*) + \mathcal{F}_\psi^{[\delta, t^*]}(u^*) & \leq \sum\limits_{\psi = T^3 w, T^3 \underline{w}} \int_\delta^{t^*} \frac{0.3}{t} \mathcal{E}_\psi^{[0, u^*]}(t) dt + 100 \mathrm{C}_{\mathrm{high}}^2 \int_0^{u^*} \mathcal{F}_{T^3 w}^{[\delta, t^*]}(u) du \\
    & \qquad + \mathrm{C}_{\mathrm{low}, 3} \frac{M {t^*}^2}{\epsilon} + \sum\limits_{\psi = T^3 w, T^3 \underline{w}} \mathcal{E}_\psi^{[0, u^*]}(\delta) + \mathcal{F}_\psi^{[\delta, t^*]}(0)
\end{aligned}
\end{equation*}

Let $\hat{u} \in [0, u^*]$ be the place where the flux attains its maximum.
\begin{equation*}
    \sum\limits_{\psi = T^3 w, T^3 \underline{w}}\mathcal{F}_\psi^{[\delta, t^*]}(\hat{u}) = \max\limits_{u \in [0, u^*]} \sum\limits_{\psi = T^3 w, T^3 \underline{w}}\mathcal{F}_\psi^{[\delta, t^*]}(u)
\end{equation*}

And we can repeat the energy estimate on $(t, u) \in [\delta, t^*] \times [0, \hat{u}]$ to get the following.
\begin{equation*}
\begin{aligned}
    \sum\limits_{\psi = T^3 w, T^3 \underline{w}} \mathcal{E}_\psi^{[0, \hat{u}]}(t^*) + \mathcal{F}_\psi^{[\delta, t^*]}(\hat{u}) & \leq \sum\limits_{\psi = T^3 w, T^3 \underline{w}} \int_\delta^{t^*} \frac{0.3}{t} \mathcal{E}_\psi^{[0, \hat{u}]}(t) dt + 100 \mathrm{C}_{\mathrm{high}}^2 \int_0^{\hat{u}} \mathcal{F}_{T^3 w}^{[\delta, t^*]}(u) du \\
    & \qquad + \mathrm{C}_{\mathrm{low}, 3} \frac{M {t^*}^2}{\epsilon} + \sum\limits_{\psi = T^3 w, T^3 \underline{w}} \mathcal{E}_\psi^{[0, \hat{u}]}(\delta) + \mathcal{F}_\psi^{[\delta, t^*]}(0)
\end{aligned}
\end{equation*}

Drop the energy on the left hand side, use the smallness of $u^*$ and the maximizing property of $\hat{u}$, we get the following.
\begin{equation*}
\begin{aligned}
    \sum\limits_{\psi = T^3 w, T^3 \underline{w}} \mathcal{F}_\psi^{[\delta, t^*]}(\hat{u}) & \leq \sum\limits_{\psi = T^3 w, T^3 \underline{w}} \int_\delta^{t^*} \frac{0.3}{t} \mathcal{E}_\psi^{[0, \hat{u}]}(t) dt + \frac{1}{2} \mathcal{F}_{T^3 w}^{[\delta, t^*]}(\hat{u}) \\
    & \qquad + \mathrm{C}_{\mathrm{low}, 3} \frac{M {t^*}^2}{\epsilon} + \sum\limits_{\psi = T^3 w, T^3 \underline{w}} \mathcal{E}_\psi^{[0, \hat{u}]}(\delta) + \mathcal{F}_\psi^{[\delta, t^*]}(0)
\end{aligned}
\end{equation*}

This gives the uniform bound of flux in terms of energy.
\begin{equation*}
\begin{aligned}
    \max\limits_{u \in [0, u^*]} \sum\limits_{\psi = T^3 w, T^3 \underline{w}} \mathcal{F}_\psi^{[\delta, t^*]}(u) & \leq \sum\limits_{\psi = T^3 w, T^3 \underline{w}} \int_\delta^{t^*} \frac{0.6}{t} \mathcal{E}_\psi^{[0, u^*]}(t) dt + 2 \mathrm{C}_{\mathrm{low}, 3} \frac{M {t^*}^2}{\epsilon} \\
    & \qquad + 2 \sum\limits_{\psi = T^3 w, T^3 \underline{w}} \mathcal{E}_\psi^{[0, u^*]}(\delta) + \mathcal{F}_\psi^{[\delta, t^*]}(0)
\end{aligned}
\end{equation*}

Substitute this into the energy estimate \ref{eq:local-base-energy_estimate} and drop the flux term on the left hand side to get an estimate of the energy.
\begin{equation*}
\begin{aligned}
    \sum\limits_{\psi = T^3 w, T^3 \underline{w}} \mathcal{E}_\psi^{[0, u^*]}(t^*) & \leq \sum\limits_{\psi = T^3 w, T^3 \underline{w}} \int_\delta^{t^*} \frac{0.6}{t} \mathcal{E}_\psi^{[0, u^*]}(t) dt + 2 \mathrm{C}_{\mathrm{low}, 3} \frac{M {t^*}^2}{\epsilon} \\
    & \qquad + 2 \sum\limits_{\psi = T^3 w, T^3 \underline{w}} \mathcal{E}_\psi^{[0, u^*]}(\delta) + \mathcal{F}_\psi^{[\delta, t^*]}(0) \\
    & \leq \sum\limits_{\psi = T^3 w, T^3 \underline{w}} \int_\delta^{t^*} \frac{0.6}{t} \mathcal{E}_\psi^{[0, u^*]}(t) dt + 2 (\frac{\mathrm{C}_{\mathrm{low}, 3} M}{\epsilon} + E_{\mathrm{Data}, 3} + F_{\mathrm{Data}, 3}) {t^*}^2 \\
\end{aligned}
\end{equation*}

Now we apply the standard Gronwall inequality. Define two temporary functions as follows.
\begin{equation*}
    \alpha(t) = \mathrm{C}_\alpha t^2 = 2 (\frac{\mathrm{C}_{\mathrm{low}, 3} M}{\epsilon} + E_{\mathrm{Data}, 3} + F_{\mathrm{Data}, 3}) t^2, \quad \beta(t) = \frac{\mathrm{C}_\beta}{t} = \frac{0.6}{t}
\end{equation*}

Then we have the following estimate for energy.
\begin{equation*}
\begin{aligned}
    \sum\limits_{\psi = T^3 w, T^3 \underline{w}} \mathcal{E}_\psi^{[0, u^*]}(t^*) & \leq \alpha(t^*) + \int_\delta^{t^*} \alpha(t) \beta(t) e^{\int_t^{t^*} \beta(\tau) d\tau} dt \\
    & = \mathrm{C}_\alpha {t^*}^2 + \int_\delta^{t^*} 0.6 \mathrm{C}_\alpha t \cdot e^{0.6 \ln(\frac{t^*}{t})} dt \\
    & \leq \frac{2}{1.4} \mathrm{C}_\alpha {t^*}^2
\end{aligned}
\end{equation*}

From the expression of $\mathrm{C}_\alpha, \epsilon, M$, we conclude that we have proved the following.
\begin{equation*}
    \sum\limits_{\psi = T^3 w, T^3 \underline{w}} \mathcal{E}_\psi^{[0, u^*]}(t) \lesssim M^2 t^2, \quad \forall t \in [\delta, t^*]
\end{equation*}

Perform a similar argument to the beginning of the proof, we have the following.
\begin{equation*}
\begin{cases}
    |L T^2 w(t, u)| \leq \mathrm{C}' M \sqrt{u^*} + M_2 \\
    |L T^2 \underline{w}(t, u)| \leq \mathrm{C}' M \sqrt{u^*} + M_2 \\
    |T^3 w(t, u)| \leq \mathrm{C}' M t \sqrt{u^*} + M_2 t \\
    |T^3 \underline{w}(t, u)| \leq \mathrm{C}' M t \sqrt{u^*}  + M_2 t \\
\end{cases}
\end{equation*}

This time we are not dropping the $\sqrt{u^*}$ term using $u^* \leq 1$. Now we see that the final step to close the whole bootstrap proof is to improve the bootstrap assumption \ref{eq:local-base-bootstrap}. We just need to place another smallness restriction on $u^*$.
\begin{equation} \label{eq:local-base-smallness_close}
    u^* \leq \frac{(\mathrm{C} - 1)^2}{2 \mathrm{C}'^2}
\end{equation}

Now we collect the smallness requirements \ref{eq:local-base-smallness_sim}, \ref{eq:local-base-smallness_Gronwall}, \ref{eq:local-base-smallness_close} that we place on $(t^*, u^*)$ during the proof.
\begin{equation*}
    t^* \leq \min\{\frac{1}{C_1' M}, \frac{1}{10 \mathrm{C}_{\mathrm{high}} M}\}, \quad u^* \leq \min\{\frac{1}{200 \mathrm{C}_{\mathrm{high}}^2}, \frac{\mathrm{C}^2}{2 \mathrm{C}'^2}\}
\end{equation*}

This gives the smallness requirement \ref{eq:local-base-t_u_smallness} in the proposition.

And we finish our proof by listing out the estimates we obtain and showing that we have obtained all the estimates in \ref{eq:local-base-final_estimates}.

By the Gronwall argument we already have the estimate for energy.
\begin{equation*}
    \sum\limits_{\psi = T^3 w, T^3 \underline{w}} \mathcal{E}_\psi^{[0, u^*]}(t) \lesssim M^2 t^2, \quad \forall t \in [\delta, t^*]
\end{equation*}

Notice that the bound for flux is the same as the one for the energy, so we immediately obtain the estimate for flux.
\begin{equation*}
    \max\limits_{u \in [0, u^*]} \sum\limits_{\psi = T^3 w, T^3 \underline{w}} \mathcal{F}_\psi^{[\delta, t]}(u) \lesssim M^2 t^2, \quad \forall t \in [\delta, t^*]
\end{equation*}

Then we obtain the estimate for $T^3 \kappa$ by using the propagation equation for $\kappa$.
\begin{equation*}
\begin{aligned}
    \int_{u, S_{t, u}} |T^3 \kappa|^2 d\mu_{\slashed{g}} du & \lesssim \int_u |T^3 \kappa(t, u)|^2 du \\
    & \lesssim \int_u |T^3 \kappa(\delta, u|^2 du + \int_u (\int_\delta^t |T^4 w(\tau, u)| + |T^4 \underline{w}(\tau, u)| d\tau)^2 du \\
    & \lesssim M^2 \delta^4 + \int_u \frac{1}{2} (t^2 - \delta^2) \int_\delta^t \frac{|T^4 w(\tau, u)|^2 + |T^4 \underline{w}(\tau, u)|^2}{\tau} d\tau du \\
    & \lesssim M^2 \delta^4 + t^2 \int_\delta^t \frac{M^2 \tau^2}{\tau} d\tau \\
    & \lesssim M^2 t^4
\end{aligned}
\end{equation*}

The $L^\infty$-bounds have already been established in \ref{eq:local-base-Linfty}, \ref{eq:local-sim_relations} at the beginning of the proof.

\subsection{Proof of Proposition \ref{prop:local-regularity}} \label{subsec:local-proof-regularity}

The proof is quite similar to the proof of proposition \ref{prop:local-base_energy}.

There are two differences. One is that now we need to be careful about the restrictions on the smallness of $(t^*, u^*)$, which is mentioned in remark \ref{rmk:high_low_constants}. Another is that we do not need to bootstrap, so we can be more casual about the $L^\infty$-bounds we use during the estimates.

We introduce the notation $\lesssim_A$, which means that the hidden constants may also depend on $A$.

We also introduce the notation $M_{\leq n}' := \max\limits_{0 \leq i \leq n} M_i'$ for simplicity.

Now we consider the energy estimates for $\psi = T^N w, T^N \underline{w}$.

The term $\int_{[\delta, t^*] \times [0, u^*]} Q_{1, \psi} + Q_{2, \psi} d\mu_g$ is exactly the same.
\begin{equation*}
    \int_{[\delta, t^*] \times [0, u^*]} Q_{1, \psi} + Q_{2, \psi} d\mu_g \lesssim \int_\delta^{t^*} M \mathcal{E}_\psi^{[0, u^*]}(t)
\end{equation*}

Then we also separate the main part $\int_{[\delta, t^*] \times [0, u^*]} Q_{0, T^N \psi} d\mu_g$ with $\psi = w, \underline{w}$ into the commutator part and the source part as in \ref{eq:local-base-error_split_1}, \ref{eq:local-base-error_split_2}.
\begin{equation*}
\begin{aligned}
    & \int_{[\delta, t^*] \times [0, u^*]} Q_{0, T^N \psi} d\mu_g \\
    = & \int_{t, u, S_{t, u}} - (\underbrace{\sum\limits_{n_1 + n_2 = N - 1, n_i \geq 0} T^{n_1} [\mu \Box_g, T] T^{n_2} \psi}_{\mathrm{comm}} + \underbrace{T^N (\mu \tilde{\rho}_\psi)}_{\mathrm{source}}) (\frac{\kappa}{c} L T^N \psi + \underline{L} T^N \psi) d\mu_\slashed{g} du dt \\
    = & I_{\mathrm{comm}, \psi} + I_{\mathrm{source}, \psi}
\end{aligned}
\end{equation*}

\begin{equation*}
\begin{aligned}
& \begin{aligned}
    I_{\mathrm{comm}, \psi} & = \int_{t, u, S_{t, u}} - (\frac{\kappa}{c} L T^N \psi + \underline{L} T^N \psi) \sum\limits_{n_1 + n_2 = N - 1,\ n_i \geq 0} \mathrm{C}_{n_1, n_2} \\
    & \qquad (\underbrace{2 T^{n_1 + 1} (\frac{\kappa v}{c r}) L T^{n_2} \psi}_{1} + \underbrace{2 T^{n_1 + 1} (\frac{v + c}{r}) T^{n_2 + 1} \psi}_{2} \\
    & \qquad + \underbrace{T^{n_1 + 1} (\frac{\kappa}{c}) L L T^{n_2} \psi}_{3} + \underbrace{L T^{n_1 + 1} (\frac{\kappa}{c}) L T^{n_2} \psi}_{4}) d\mu_\slashed{g} du dt \\
    & = I_{\mathrm{comm}, 1, \psi} + I_{\mathrm{comm}, 2, \psi} + I_{\mathrm{comm}, 3, \psi} + I_{\mathrm{comm}, 4, \psi} \\
\end{aligned} \\
& \begin{aligned}
    I_{\mathrm{source}, \psi} = \int_{t, u, S_{t, u}} - T^N (\mu \tilde{\rho}_\psi) (\frac{\kappa}{c} L T^N \psi + \underline{L} T^N \psi) d\mu_\slashed{g} du dt \\
\end{aligned} \\
\end{aligned}
\end{equation*}

We estimate these term by term just like in the previous subsection. However, now we need to track both the highest order terms and the second highest order terms.

We use $T^N (\mu \tilde{\rho}_w)$ as an example.
\begin{equation*}
\begin{aligned}
    T^N (\mu \tilde{\rho}_w) & = T^N (2 \frac{\kappa v^2}{r^2} + 2 \frac{v - c}{r} Tw + L(\frac{\kappa v}{r})) \\
    (\mathrm{take}\ |\cdot|) & \lesssim_{M_{\leq N - 2}'} (|\frac{2 v^2}{r^2} T^N \kappa| + |\frac{4 \kappa v}{r^2} T^N v| + |(N T (\frac{2 v^2}{r^2}) + \frac{4 \kappa v^2}{r^3}) T^{N - 1} \kappa| + t^2) \\
    & \quad + (|\frac{2 (v - c)}{r} T^{N + 1} w| + |\frac{2 T w}{r} T^N (v - c)| + |N T (\frac{2 (v - c)}{r}) T^N w| + |\frac{2 (v - c) T w}{r^2} T^{N - 1} \kappa| + t^2) \\
    & \quad + (|\frac{v}{r} L T^N \kappa| + |L(\frac{v}{r}) T^N \kappa| + |\frac{\kappa}{r} L T^N v| \\
    & \qquad + |(N T (\frac{v}{r}) + \frac{\kappa v}{r^2}) L T^{N - 1} \kappa| + |(N L T (\frac{v}{r}) + \frac{L (\kappa v)}{r^2}) T^{N - 1} \kappa| \\
    & \qquad + |L (\frac{\kappa}{r}) T^N v| + |N T (\frac{\kappa}{r}) L T^{N - 1} v| + t) \\
    & \lesssim_{M_{\leq N - 2}'} \underbrace{M |T^N \kappa| + |T^{N + 1} (w, \underline{w})| + t |LT^N (w, \underline{w})|}_{\mathrm{hidden\ constants\ are\ the\ same\ as\ } N = 3} \\
    & \qquad + |T^{N - 1} \kappa| + |T^N (w, \underline{w})| + t^2 |L T^{N - 1} (w, \underline{w})| + t \\
\end{aligned}
\end{equation*}

The other terms are estimated in the same way.
\begin{equation*}
\begin{aligned}
    T^N (\mu \tilde{\rho}_{\underline{w}}) & = T^N (2 \frac{\kappa v^2}{r^2} - 2 \frac{v + c}{r} T\underline{w} + L(\frac{\kappa v}{r}) + 2T(\frac{c v}{r})) \\
    & \lesssim_{M_{\leq N - 2}'} M |T^N \kappa| + |T^{N + 1} (w, \underline{w})| + t |LT^N (w, \underline{w})| \\
    & \qquad + |T^{N - 1} \kappa| + |T^N (w, \underline{w})| + t^2 |L T^{N - 1} (w, \underline{w})| + t \\
\end{aligned}
\end{equation*}

\begin{equation*}
\begin{aligned}
    & |\sum\limits_{n_1 + n_2 = N - 1,\ n_i \geq 0} \mathrm{C}_{n_1, n_2} T^{n_1 + 1} (\frac{\kappa v}{c r}) L T^{n_2} \psi| \lesssim_{M_{\leq N - 2}'} M |T^N \kappa| + |T^{N - 1} \kappa| + t |T^N (w, \underline{w})| + t |L T^{N - 1} \psi| + t \\
    & |\sum\limits_{n_1 + n_2 = N - 1,\ n_i \geq 0} \mathrm{C}_{n_1, n_2} T^{n_1 + 1} (\frac{v + c}{r}) T^{n_2 + 1} \psi| \lesssim_{M_{\leq N - 2}'} |T^{N - 1} \kappa| + |T^N (w, \underline{w})| + t^2 \\
    & \begin{aligned}
        |\sum\limits_{n_1 + n_2 = N - 1,\ n_i \geq 0} \mathrm{C}_{n_1, n_2} L T^{n_1 + 1} (\frac{\kappa}{c}) L T^{n_2} \psi| & \lesssim_{M_{\leq N - 2}'} M |T^{N + 1} (w, \underline{w})| + M^2 |T^N \kappa| + M t |L T^N (w, \underline{w})| \\
        & + |T^N (w, \underline{w})| + |T^{N - 1} \kappa| + t |L T^{N - 1} (w, \underline{w})| + |L T^{N - 1} \psi| + 1 \\
    \end{aligned}
\end{aligned}
\end{equation*}

\begin{equation*}
\begin{aligned}
    & |\sum\limits_{n_1 + n_2 = N - 1,\ n_i \geq 0} \mathrm{C}_{n_1, n_2} T^{n_1 + 1} (\frac{\kappa}{c}) L L T^{n_2} w| \\
    = & |\sum\limits_{n_1 + n_2 = N - 1,\ n_i \geq 0} \mathrm{C}_{n_1, n_2} T^{n_1 + 1} (\frac{\kappa}{c}) T^{n_2} (- 2 L (\frac{c}{\kappa}) T w - 2 \frac{c}{\kappa} L T w - L (\frac{c v}{r}))| \\
    \lesssim_{M_{\leq N - 2}'} & \underbrace{\frac{M}{t} |T^N \kappa| + N |L T^N w|}_{\mathrm{hidden\ constants\ don't\ depend\ on\ }N} + \frac{1}{t} |T^{N - 1} \kappa| + |L T^{N - 1} w| + t |L T^{N - 1} \underline{w}| + \frac{1}{t} |T^N w| + |T^N \underline{w}| + 1 \\
\end{aligned}
\end{equation*}

Now, just as in the case of $N = 3$, we sum these up and use the AM-GM inequality to separate the $L^2$-norms and get the following estimate.
\begin{equation*}
\begin{aligned}
    & |\sum\limits_{\psi = T^N w, T^N \underline{w}} \mathcal{E}_\psi^{[0, u^*]}(t^*) - \mathcal{E}_\psi^{[0, u^*]}(\delta) + \mathcal{F}_\psi^{[\delta, t^*]}(u^*) - \mathcal{F}_\psi^{[\delta, t^*]}(0)| \\
    \lesssim_{M_{\leq N - 2}'} & \sum\limits_{\psi = T^N w, T^N \underline{w}} \int_\delta^{t^*} (\frac{M \epsilon}{t} + M) \mathcal{E}_\psi^{[0, u^*]}(t) dt + \int_{t, u, S_{t, u}} \frac{M}{\epsilon t} |T^N \kappa|^2 d\mu_\slashed{g} du dt \\
    & + \frac{1}{M \epsilon} \int_0^{u^*} \mathcal{F}_{T^N w}^{[\delta, t^*]}(u) du + \sum\limits_{\psi = T^{N - 1} w, T^{N - 1} \underline{w}} \int_\delta^{t^*} \frac{1}{t} \mathcal{E}_\psi^{[0, u^*]}(t) dt \\
    & + \int_\delta^{t^*} \frac{1}{t} \int_{u, S_{t, u}} |T^{N - 1} \kappa|^2 d\mu_\slashed{g} du dt + {t^*}^2 \\
\end{aligned}
\end{equation*}

And we also have the same estimate for $\int |T^N \kappa|^2$.
\begin{equation*}
    \int_{t, u, S_{t, u}} \frac{M}{\epsilon t} |T^N \kappa|^2 d\mu_\slashed{g} du dt \lesssim \frac{M E_{\mathrm{Data}, N} \delta^3}{\epsilon} + \sum\limits_{\psi = T^N w, T^N \underline{w}} \int_\delta^{t^*} \frac{M {t^*}^2}{\epsilon t} \mathcal{E}_\psi^{[0, u^*]}(t) dt
\end{equation*}

We choose the same $\epsilon$ as in proposition \ref{prop:local-base_energy} but require $u^* \leq \frac{1}{20 N^2 \mathrm{C}_{\mathrm{high}}^2}$. We also substitute the known lower order energy bounds into the inequality. Then we get the desired energy estimate.
\begin{equation*}
\begin{aligned}
    & \sum\limits_{\psi = T^N w, T^N \underline{w}} \mathcal{E}_\psi^{[0, u^*]}(t^*) + \mathcal{F}_\psi^{[\delta, t^*]}(u^*) \\
    \leq & \sum\limits_{\psi = T^N w, T^N \underline{w}} \int_\delta^{t^*} \frac{0.3}{t} \mathcal{E}_\psi^{[0, u^*]}(t) dt + 10 N^2 \mathrm{C}_{\mathrm{high}}^2 \int_0^{u^*} \mathcal{F}_{T^3 w}^{[\delta, t^*]}(u) du \\
    & + \mathrm{C}_{\mathrm{low}, N} ({t^*}^2 + E_{\mathrm{Data}, N} \delta^3 + E_{N - 1}' {t^*}^2) + \sum\limits_{\psi = T^N w, T^N \underline{w}} \mathcal{E}_\psi^{[0, u^*]}(\delta) + \mathcal{F}_\psi^{[\delta, t^*]}(0)
\end{aligned}
\end{equation*}

Here, $\mathrm{C}_{\mathrm{low}, N}$ depends on $c_0, v_0, r_0, N, M_{\leq N - 2}'$.

We make the same Gronwall type argument to get the final result.
\begin{equation*}
\begin{aligned}
    & \sum\limits_{\psi = T^N w, T^N \underline{w}} \mathcal{E}_\psi^{[0, u^*]}(t) \leq E_N' t^2 \\
    & \sup\limits_{u \in [0, u^*]} \sum\limits_{\psi = T^N w, T^N \underline{w}} \mathcal{F}_\psi^{[\delta, t]}(u) \leq E_N' t^2 \\
    & \int_0^{u^*} \int_{S_{t, u}} |T^N \kappa(t, u)|^2 d\mu_{\slashed{g}} du \leq E_N' t^4
\end{aligned}
\end{equation*}

$E_N'$ depends on $c_0, v_0, r_0, N, M_{\leq N - 2}', E_{N - 1}', E_{\mathrm{Data}, N}, F_{\mathrm{Data}, N}$. Then we integrate from the boundary $C_0$ to obtain the corresponding $L^\infty$-bounds.
\begin{equation*}
\begin{cases}
\begin{aligned}
     & |L T^{N - 1} w| \leq M_{N - 1}' && |L T^{N - 1} \underline{w}| \leq M_{N - 1}' && |L T^{N - 1} \kappa| \leq M_{N - 1}' t \\
    & |T^N w| \leq M_{N - 1}' t && |T^N \underline{w}| \leq M_{N - 1}' t && |T^{N - 1} \kappa| \leq M_{N - 1}' t^2 \\
\end{aligned}
\end{cases}
\end{equation*}

$M_{N - 1}'$ depends on $c_0, r_0, M_{N - 1}, E_N'$.

Since we introduced new smallness requirements $u^* \leq \frac{1}{20 N^2 \mathrm{C}_{\mathrm{high}}^2}$ during the proof, we need to divide the interval $[0, u^*]$ into sub intervals $[0, \frac{1}{N^2} u^*], \cdots, [\frac{(N^2 - 1)}{N^2} u^*, u^*]$ and perform the energy estimates on each sub interval.

This finishes the proof.

\subsection{Convergence to the Rarefaction Wave Solution}

The energy estimates we proved in propositions \ref{prop:local-base_energy} and \ref{prop:local-regularity} provide a conditional uniform lower bound of the lifespan $t^*$ and conditional uniform upper bounds of the norms. We state these as the following corollary.

\begin{corollary}[Conditional Uniform Control of Approximate Solutions] \label{coro:local-control_qpproximate_solutions}
Let $\Psi_{\delta, N}$ defined on $(t, u) \in [\delta, t^*(\delta, u^*, N)] \times [0, u^*]$ be the approximate solution we constructed in subsection \ref{subsec:local-localexist}.

If we have the smallness of $\delta$.
\begin{equation*}
    u^* \leq \frac{1}{\mathrm{C}_2}
\end{equation*}

And we have the following condition on the initial time slice $\Sigma_\delta$, where $3 \leq m \leq N$.
\begin{equation*}
\begin{cases}
    \mathcal{E}_{T^m w}^{[0, u^*]}(\delta) + \mathcal{E}_{T^m \underline{w}}^{[0, u^*]}(\delta) \leq E_{\mathrm{Data}, m} \delta^2 \\
    \int_{[0, u^*]} \int_{S_{\delta, u}} |T^m \kappa|^2 d\mu_{\slashed{g}} du \leq E_{\mathrm{Data}, m} \delta^4 \\
    \mathcal{F}_{T^m w}^{[\delta, t^*]}(0) + \mathcal{F}_{T^m \underline{w}}^{[\delta, t^*]}(0) \leq F_{\mathrm{Data}, m} {t^*}^2 \\
    \begin{aligned}
        & |c(\delta, u) - c_0| \leq \frac{1}{4} c_0 && |v(\delta, u) - v_0| \leq \frac{1}{4} c_0 && |r(\delta, u) - r_0| \leq \frac{1}{4} r_0 \\
        & |L \kappa(\delta, u) - 1| \leq \frac{1}{4} && |\kappa(\delta, u) - \delta| \leq \frac{1}{4} \delta && \\
    \end{aligned}
\end{cases}
\end{equation*}

Then the approximate solution can be extended to the following time.
\begin{equation*}
    t^* = \frac{1}{\mathrm{C}_1 \max\{M_0, M_1, M_2, \sqrt{E_{\mathrm{Data}, 3}}, \sqrt{F_{\mathrm{Data}, 3}}\}}
\end{equation*}

And the solution has the following bounds of norms, where $3 \leq m \leq N, 1 \leq n \leq N - 1$.
\begin{equation*}
\begin{cases}
    \mathcal{E}_{T^m w}^{[0, u^*]}(t) + \mathcal{E}_{T^m \underline{w}}^{[0, u^*]}(t) \leq E_m' t^2 \\
    \max\limits_{u \in [0, u^*]} \mathcal{F}_{T^m w}^{[\delta, t]}(u) + \mathcal{F}_{T^m \underline{w}}^{[\delta, t]}(u) \leq E_m' t^2 \\
    \int_{[0, u^*]} \int_{S_{t, u}} |T^m \kappa|^2 d\mu_{\slashed{g}} du \leq E_m' t^4 \\
    \begin{aligned}
        & |c - c_0| \leq \frac{1}{2} c_0 && |v - v_0| \leq \frac{1}{2} v_0 && |r - r_0| \leq \frac{1}{2} r_0 \\
        & |L w| \leq M_0' && |L \underline{w}| \leq M_0' && |L \kappa - 1| \leq M_0' t \leq \frac{1}{2} \\
        & |T w| \leq M_0' t && |T \underline{w} + \frac{2}{\gamma + 1}| \leq M_0' t && |\kappa - t| \leq M_0' t^2 \leq \frac{1}{2} t \\
        & |L T^n w| \leq M_n' && |L T^n \underline{w}| \leq M_n' && |L T^n \kappa| \leq M_n' t \\
        & |T^{n + 1} w| \leq M_n' t && |T^{n + 1} \underline{w}| \leq M_n' t && |T^n \kappa| \leq M_n' t^2 \\
    \end{aligned}
\end{cases}
\end{equation*}

The constants depend on $E_{\mathrm{Data}, m}$ and the data on $C_0$ (already fixed by the background solution, and determine $F_{\mathrm{Data}, m}$).
\end{corollary}

Now we want to use the family of approximate solutions to pass to the desired rarefaction wave solution by convergence. That is, we want to find a sequence of $\Psi_{\delta_N, N}$ such that
\begin{itemize}
    \item $N \to \infty$
    \item $\delta_N \to 0$
    \item These approximate solutions can be defined on a uniform region of $(t, u) \in [\delta, t^*] \times [0, u^*]$.
    \item We have uniform bounds for all the norms up to $C^\infty$.
\end{itemize}

The corollary \ref{coro:local-control_qpproximate_solutions} reduces these requirements to the problem of obtaining uniform control of the $N$-th order approximation data at time $t = \delta_N$ on a uniform range of $u$.

The main obstacle is that when we add more and more orders of approximations, there will be more and more terms in the Taylor expansion, contributing to the size of the initial data. We overcome this obstacle by observing that the behavior of the solution at the singular point is, in fact, well structured.

We first list the expression \ref{eq:initial_Taylor} for the $N$-th order approximation data at time $t = \delta$.
\begin{equation*}
\begin{cases}
\begin{aligned}
    & w_N(\delta, u) = w(\delta, 0) + \sum\limits_{n = 1}^N \frac{T^n w(\delta, 0)}{n!} u^n && (L w)_N(\delta, u) = - \frac{c_N v_N}{r_N} - 2 \frac{c_N}{\kappa_N} T w_N \\
    & \underline{w}_N(\delta, u) = \underline{w}(\delta, 0) + \sum\limits_{n = 1}^N \frac{T^n \underline{w}(\delta, 0)}{n!} u^n && (L \underline{w})_N(\delta, u) = - \frac{c_N v_N}{r_N} \\
    & r_N(\delta, u) = r(\delta, 0) - \sum\limits_{n = 0}^{N - 1} \frac{T^n \kappa(\delta, 0)}{n!} u^n && L r_N(\delta, u) = c_N + v_N \\
\end{aligned}
\end{cases}
\end{equation*}

Now we state the proposition.

\begin{proposition}[Existence of Uniform Approximation Data] \label{prop:local-uniform_data}
\begin{equation*}
    \exists u_0 = u_0(c_0), \quad \delta_N = \delta_N(N)
\end{equation*}

Such that we have uniform bounds of norms for the $N$-th approximation data. ($3 \leq m \leq N$)
\begin{equation*}
\begin{cases}
    \mathcal{E}_{T^m w_N}^{[0, u_0]}(\delta_N) + \mathcal{E}_{T^m \underline{w}_N}^{[0, u_0]}(\delta_N) \leq E_{\mathrm{Data}, m} \delta_N^2 \\
    \int_{[0, u_0]} \int_{S_{\delta_N, u}} |T^m \kappa_N|^2 d\mu_{\slashed{g}} du \leq E_{\mathrm{Data}, m} \delta_N^4 \\
    \begin{aligned}
        & |c_N(\delta_N, u) - c_0| \leq \frac{1}{4} c_0 && |v_N(\delta_N, u) - v_0| \leq \frac{1}{4} c_0 && |r_N(\delta_N, u) - r_0| \leq \frac{1}{4} r_0 \\
        & |L \kappa_N(\delta_N, u) - 1| \leq \frac{1}{4} && |\kappa_N(\delta_N, u) - \delta_N| \leq \frac{1}{4} \delta_N && \\
    \end{aligned}
\end{cases}
\end{equation*}

Where $E_{\mathrm{Data}, m}$ is uniform in $N$.
\end{proposition}

Now we prove the proposition.

We already know that $T^n w(t, 0), T^n \underline{w}(t, 0), T^n r(t, 0)$ are all smooth functions of $t$, because of the smoothness of the background solution and the fact that these functions are solved from ODEs with smooth coefficients recursively. Then, from the expressions we immediately know that $w_N, \underline{w}_N, r_N$ (then $c_N, v_N$) are smooth functions of $(t, u)$ for each fixed $N$.

We also obtain the same vanishing order as \ref{eq:local-vanishing_C0} for $w_N, \underline{w}_N, r_N$ for each fixed $N$.
\begin{equation*}
\begin{cases}
    \begin{aligned}
    & T w_N = O(t), && T^n w_N = O(t) \\
    & T \underline{w}_N = - \frac{2}{\gamma + 1} + O(t), && T^n \underline{w}_N = O(t) \\
    & \kappa_N = t + O(t^2), && T^{n - 1} \kappa_N = O(t^2)
    \end{aligned}, \quad (n \geq 2)
\end{cases}
\end{equation*}

Here $O$ is uniform in $u$.

The vanishing orders also imply the following.
\begin{equation*}
\begin{cases}
\begin{aligned}
    & T c_N = - \frac{\gamma - 1}{\gamma + 1} + O(t) && T v_N = - \frac{2}{\gamma + 1} + O(t) \\
    & L \kappa_N = 1 + O(t) && T r_N = - t + O(t^2)
\end{aligned}
\end{cases}
\end{equation*}

Then we go down to $t = 0$ and get the value of $c_N(0, u), v_N(0, u), r_N(0, u)$.
\begin{equation*}
\begin{cases}
    c_N(0, u) = c_0 - \frac{\gamma - 1}{\gamma + 1} u \\
    v_N(0, u) = v_0 - \frac{2}{\gamma + 1} u \\
    r_N(0, u) = r_0
\end{cases}
\end{equation*}

Now we give the smallness requirement for $u^*$.
\begin{equation*}
    u^* \leq u_0 = \frac{1}{8} c_0 \min\{|\frac{\gamma + 1}{\gamma - 1}|, |\frac{\gamma + 1}{2}|\}
\end{equation*}

This is to make sure that at $t = 0$ we have the following.
\begin{equation*}
    |c_N(0, u) - c_0| \leq \frac{1}{8} c_0, \quad |v_N(0, u) - v_0| \leq \frac{1}{8} c_0
\end{equation*}

Using the continuity of the functions and the vanishing orders, we also give an upper bound $\delta_N^*$ for $\delta_N$ such that we have the following (one part of the requirements for the initial data on $\Sigma_\delta$).
\begin{equation*}
\begin{cases}
\begin{aligned}
    & |c_N(\delta, u) - c_0| \leq \frac{1}{4} c_0 && |v_N(\delta, u) - v_0| \leq \frac{1}{4} c_0 && |r_N(\delta, u) - r_0| \leq \frac{1}{4} r_0 \\
    & |L \kappa_N(\delta, u) - 1| \leq \frac{1}{4} && |\kappa_N(\delta, u) - \delta| \leq \frac{1}{4} \delta && \forall (\delta, u) \in [0, \delta_N^*] \times [0, u_0] \\
\end{aligned}
\end{cases}
\end{equation*}

With these information, we deduce that $(L w)_N, (L \underline{w})_N, L r_N$ are smooth functions on $(t, u) \in [0, \delta_N^*] \times [0, u_0]$.

Now, for each $m \geq 3$, we need to find uniform (in $N$) bounds for $E_{\mathrm{Data}, m}(\Psi_{\delta_N, N})$, which is the other part of the requirements on the initial data. To achieve this, we use $L^\infty$ to control $L^2$ since everything is smooth now. Check the expressions of the $L^2$ norms of the initial data in the corollary \ref{coro:local-control_qpproximate_solutions}, we find out that we just need to find uniform bounds for the following quantities.
\begin{equation*}
\begin{cases}
\begin{aligned}
    & T^m (L w)_N(\delta_N, u) && T^m (L w)_N(\delta_N, u) && \\
    & \frac{T^{m + 1} w_N(\delta_N, u)}{\delta_N} && \frac{T^{m + 1} w_N(\delta_N, u)}{\delta_N} && \frac{T^m \kappa_N(\delta_N, u)}{\delta_N^2}
\end{aligned}
\end{cases}
\end{equation*}

We do some reductions first. The control of $\frac{T^m \kappa_N(\delta_N, u)}{\delta_N^2}$ can be reduced to control $\frac{L T^m \kappa_N(t, u)}{t}$, which is just a combination of $\frac{T^{m + 1} w_N(t, u)}{t}, \frac{T^{m + 1} w_N(t, u)}{t}$.

Similarly, we reduce the control of $\frac{T^{m + 1} w_N(t, u)}{t}, \frac{T^{m + 1} w_N(t, u)}{t}$ to the control of $T^{m + 1} \partial_t w_N(t, u), T^{m + 1} \partial_t \underline{w}_N(t, u)$.

So we need to find uniform bounds for these quantities.
\begin{equation*}
\begin{cases}
\begin{aligned}
    & T^m (L w)_N(\delta_N, u) && T^m (L w)_N(\delta_N, u) \\
    & T^{m + 1} \partial_t w_N(t, u) && T^{m + 1} \partial_t \underline{w}_N(t, u), \quad \forall t \in [0, \delta_N]
\end{aligned}
\end{cases}
\end{equation*}

Since we can make $\delta_N$ arbitrarily close to zero, and these quantities are all smooth, we only need to find uniform bounds for them at $t = 0$.

We first deal with $T^{m + 1} \partial_t w_N(\delta_N, u), T^{m + 1} \partial_t \underline{w}_N(\delta_N, u)$.
\begin{equation*}
\begin{cases}
\begin{aligned}
    & T^{m + 1} \partial_t w_N(t, u) = \sum\limits_{n = m + 1}^N \frac{L T^n w(t, 0)}{(n - m - 1)!} u^{n - m - 1} \\
    & T^{m + 1} \partial_t \underline{w}_N(t, u) = \sum\limits_{n = m + 1}^N \frac{L T^n \underline{w}(t, 0)}{(n - m - 1)!} u^{n - m - 1}
\end{aligned}
\end{cases}
\end{equation*}

We see that to prove uniform boundedness, we need to show that for some constants $\mathrm{C}_I, \lambda$, we have the following growth speed control.
\begin{equation*}
    |L T^n w(0, 0)|, |L T^n \underline{w}(0, 0)| \leq \mathrm{C}_I \lambda^n n!, \quad \forall n \geq 4
\end{equation*}

To show this, we need to get back to the ODEs \ref{eq:ODE-general}, from which we obtain the following expressions.
\begin{equation*}
\begin{aligned}
    LT^n w & = - L (\frac{v}{2 r} + \frac{Lw}{2 c}) T^{n - 1} \kappa - (\frac{v}{2 r} + \frac{Lw}{2 c}) L T^{n - 1} \kappa \\
    & \qquad - \frac{1}{2} \sum\limits_{\substack{n_1 + n_2 = n - 1 \\ n_i \geq 0,\ n_2 \leq n - 2}} \frac{(n - 1)!}{n_1! n_2!} (L T^{n_1} (\frac{v}{r} + \frac{Lw}{c}) T^{n_2} \kappa + T^{n_1} (\frac{v}{r} + \frac{Lw}{c}) L T^{n_2} \kappa) \\
    L T^n \underline{w} & = \frac{\gamma - 1}{2} (- 2 \frac{\underline{w}}{r} T^n \underline{w} + 2 \frac{w}{r} T^n w - \frac{\underline{w}^2 - w^2}{r^2} T^{n - 1} \kappa) \\
    & \qquad + \frac{\gamma - 1}{2} \sum\limits_{\substack{n_1 + n_2 = n - 1 \\ n_i \geq 0,\ n_2 \leq n - 2}} \frac{(n - 1)!}{n_1! n_2!} (- 2 T^{n_1} \frac{\underline{w}}{r} T^{n_2 + 1} \underline{w} + 2 T^{n_1} \frac{w}{r} T^{n_2 + 1} w - T^{n_1} \frac{\underline{w}^2 - w^2}{r^2} T^{n_2} \kappa) \\
\end{aligned}
\end{equation*}

Using the fact that everything is bounded and smooth in $t$, we take $t = 0$ and drop all the vanishing terms.
\begin{equation*}
\begin{aligned}
    L T^n w(0, 0) & = - \frac{1}{2} T^{n - 1} (\frac{L w}{c}) L \kappa \\
    & = - \frac{1}{2} \sum\limits_{i = 0}^{n - 1} (- 1)^i i! \frac{(T c)^i}{c^{i + 1}} L T^{n - 1 - i} w \\
    & = - \frac{1}{2} \sum\limits_{i = 0}^{n - 1} i! \frac{(\frac{\gamma - 1}{\gamma + 1})^i}{c_0^{i + 1}} L T^{n - 1 - i} w(0, 0) \\
    L T^n \underline{w}(0, 0) & = 0 \\
\end{aligned}
\end{equation*}

So $L T^n \underline{w}(0, 0)$ is not a problem. We only need to deal with $L T^n w(0, 0)$. The calculation above implies the following inequality.
\begin{equation*}
\begin{aligned}
    |\frac{c_0^n}{(\frac{\gamma - 1}{\gamma + 1})^n} \frac{L T^n w(0, 0)}{n!}| & \leq \frac{1}{2 n} \frac{\gamma + 1}{\gamma - 1} \sum\limits_{i = 0}^{n - 1} \frac{i! (n - 1 - i)!}{(n - 1)!} |\frac{c_0^{n - 1 - i}}{(\frac{\gamma - 1}{\gamma + 1})^{n - 1 - i}}\frac{L T^{n - 1 - i} w(0, 0)}{(n - 1 - i)!}| \\
\end{aligned}
\end{equation*}

Temporarily setting $a(n) = |\frac{c_0^n}{(\frac{\gamma - 1}{\gamma + 1})^n} \frac{L T^n w(0, 0)}{n!}| \geq 0$, we can rewrite the inequality as follows.
\begin{equation*}
    a(n) \leq \frac{1}{2 n} \frac{\gamma + 1}{\gamma - 1} \sum\limits_{i = 0}^{n - 1} \frac{i! (n - 1 - i)!}{(n - 1)!} a(n - 1 - i)
\end{equation*}

We know the limit behavior of the coefficients.
\begin{equation*}
    \lim\limits_{n \to \infty} \sum\limits_{i = 0}^{n - 1} \frac{i! (n - 1 - i)!}{(n - 1)!} = 2, \quad \lim\limits_{n \to \infty} \frac{1}{2 n} \frac{\gamma + 1}{\gamma - 1} \sum\limits_{i = 0}^{n - 1} \frac{i! (n - 1 - i)!}{(n - 1)!} = 0 < 1
\end{equation*}

So, we know that $\{a(n)\}$ is a bounded sequence. Then we obtain the desired upper bound for $L T^n w(0, 0)$.
\begin{equation*}
    |L T^n w(0, 0)| \leq \mathrm{C}_I \frac{(\frac{\gamma - 1}{\gamma + 1})^n}{c_0^n} n!
\end{equation*}

Substituting these bounds into the expressions for $T^{m + 1} \partial_t w_N(0, u), T^{m + 1} \partial_t \underline{w}_N(0, u)$, we know that they are uniformly bounded on $u \in [0, \frac{\gamma + 1}{\gamma - 1} c_0 - \epsilon], \forall \epsilon > 0$.

\begin{remark}
$u = \frac{\gamma + 1}{\gamma - 1} c_0$ is the point where the solution reaches vacuum. This means that the convergence radius we obtain is precisely the full range of possible rarefaction waves. Probably this would help when one tries to construct the rarefaction solution connecting the background solution to vacuum.
\end{remark}

We now turn to $T^m (L w)_N(\delta_N, u), T^m (L w)_N(\delta_N, u)$.
\begin{equation*}
\begin{aligned}
    & T^m (L w)_N(t, u) = - T^m (\frac{c_N v_N}{r_N}) - 2 T^m (\frac{c_N}{\kappa_N} T w_N) \\
    & T^m (L \underline{w})_N(t, u) = - T^m (\frac{c_N v_N}{r_N}) \\
\end{aligned}
\end{equation*}

Then we take $t = 0$ and drop all vanishing terms.
\begin{equation*}
\begin{aligned}
    T^m (L w)_N(0, u) & = - 2 c_N(0, u) \lim\limits_{t \to 0} \frac{T^{m + 1} w_N(t, u)}{\kappa_N(t, u)} - 2 m T c_N(0, u) \lim\limits_{t \to 0} \frac{T^m w_N(t, u)}{\kappa_N(t, u)} \\
    & = - 2 c_N(0, u) T^{m + 1} \partial_t w_N(0, u) - 2 m T c_N(0, u) T^m \partial_t w_N(0, u) \\
    & = - 2 (c_0 - \frac{\gamma - 1}{\gamma + 1} u) T^{m + 1} \partial_t w_N(0, u) + 2 m \frac{\gamma - 1}{\gamma + 1} T^m \partial_t w_N(0, u) \\
    T^m (L \underline{w})_N(0, u) & = 0 \\
\end{aligned}
\end{equation*}

They are also uniformly (in $N$) bounded. This finishes the proof.

\begin{remark}
If one only wants a solution of finite regularity, then these uniformly bounded issues can be skipped. Improving regularity after convergence would be hard, because there is no obvious way to propagate regularity from the singular point.
\end{remark}

The convergence is now straightforward.

By corollary \ref{coro:local-control_qpproximate_solutions} and proposition \ref{prop:local-uniform_data}, we find a sequence of approximation solutions $\Psi_{\delta_N, N} = (w_N, \underline{w}_N, r_N)$ defined on a uniform region $(t, u) \in [\delta_N, t^*] \times [0, u^*]$ (except for the initial time $\delta_N$). We require $\delta_N$ to be decreasing to zero.

They belong to the following regularity classes.
\begin{equation*}
\begin{cases}
    w_N \in C_t^0 H_u^{N + 1} \cap C_t^1 H_u^N \\
    \underline{w}_N \in C_t^0 H_u^{N + 1} \cap C_t^1 H_u^N \\
    r_N \in C_t^1 H_u^{N + 1}
\end{cases}
\end{equation*}

And for each $m \geq 3$, they have uniform norm bounds.
\begin{equation*}
\begin{cases}
\begin{aligned}
    & \sup\limits_{N \geq m} \|\frac{w_N - \frac{1}{\gamma - 1} c_0 + \frac{1}{2} v_0}{t}\|_{C_t^0([\delta_N, t^*]) H_u^{m + 1}} < \infty && \sup\limits_{N \geq m} \|L w_N\|_{C_t^0([\delta_N, t^*]) H_u^m} < \infty \\
    & \sup\limits_{N \geq m} \|\frac{\underline{w}_N - \frac{1}{\gamma - 1} c_0 - \frac{1}{2} v_0 + \frac{2}{\gamma + 1} u}{t}\|_{C_t^0([\delta_N, t^*]) H_u^{m + 1}} < \infty && \sup\limits_{N \geq m} \|L \underline{w}_N\|_{C_t^0([\delta_N, t^*]) H_u^m} < \infty \\
    & \sup\limits_{N \geq m} \|\frac{r_N - r_0 - (v_0 + c_0) t + t u}{t^2}\|_{C_t^0([\delta_N, t^*]) H_u^{m + 1}} < \infty && \sup\limits_{N \geq m} \|\frac{L r_N - (v_0 + c_0) + u}{t}\|_{C_t^0([\delta_N, t^*]) H_u^{m + 1}} < \infty
\end{aligned}
\end{cases}
\end{equation*}

Then we adopt the common method of proving strong convergence in the low regularity space and weak convergence in the high regularity space.

The low regularity space is chosen to be $C_{t, u}^1([\delta_m, t^*] \times [0, u^*])$, which corresponds to the classical solution of Euler equations. To obtain strong convergence, we need to show that the $\{\Psi_{\delta_N, N}\}_{N \geq m}$ sequence is a bounded sequence in $C_{t, u}^2([\delta_m, t^*] \times [0, u^*])$. Recall the expression for $LLw, LL\underline{w}$.
\begin{equation*}
    L L w = - 2 L (\frac{c}{\kappa}) T w - 2 \frac{c}{\kappa} L T w - L (\frac{c v}{r}), \quad L L \underline{w} = - L (\frac{c v}{r})
\end{equation*}

Using the fact that the Sobolev space $H_u^1$ is an algebra, and $\kappa \sim t, r \sim r_0$, and that we have uniform bounds for $L^{\leq 1}(w_N, \underline{w}_N, r_N)$, we immediately know that the sequence is bounded in $C_{t, u}^2([\delta_m, t^*] \times [0, u^*])$.

Then the Arzelà-Ascoli theorem gives the compactness of the sequence $\{\Psi_{\delta_N, N}\}_{N \geq m}$ in $C_{t, u}^1([\delta_m, t^*] \times [0, u^*])$.

We can prove this for every $m \geq 3$, so by a diagonal sequence trick, we can find a subsequence strongly converging in $C_{t, u}^1([\delta_m, t^*] \times [0, u^*])$ for every $m$.

We denote the converging subsequence as $\{\Psi_{\delta_{N_i}, N_i}\}_{i \geq 1}$, and denote the limit as $\Psi = (w, \underline{w}, r)$.

Since we are taking limits in $C^1$ solutions to the Euler equations, the limit is also a classical solution to the Euler equations.

Now we discuss the regularity of the limit. We prove its smoothness by weak convergence.

Fix a time slice $\Sigma_t$, on this time slice we have the following strong convergence.
\begin{equation*}
\begin{cases}
    (w_{N_i}, \underline{w}_{N_i}, r_{N_i}) \to (w, \underline{w}, r) \quad \mathrm{in}\ C_u^1 \\
    ((L w)_{N_i}, (L \underline{w})_{N_i}, L r_{N_i}) \to (L w, L \underline{w}, L r) \quad \mathrm{in}\ C_u^0
\end{cases}
\end{equation*}

And we know that the left hand side is bounded in every $H_u^m$ space (need to neglect finitely many terms).

Now we use the property of Hilbert spaces: every bounded sequence has a weakly converging subsequence. This property, together with the uniqueness of the converging limit, results in the following.
\begin{equation*}
\begin{cases}
    (w_{N_i}, \underline{w}_{N_i}, r_{N_i}) \rightharpoonup (w, \underline{w}, r) \quad \mathrm{in}\ H_u^m \\
    ((L w)_{N_i}, (L \underline{w})_{N_i}, L r_{N_i}) \rightharpoonup (L w, L \underline{w}, L r) \quad \mathrm{in}\ H_u^m, \qquad \forall m \geq 0
\end{cases}
\end{equation*}

This tells us that the limit has the same regularity as the converging subsequence, and has the same norm bounds as the subsequence. Since the limit is also a solution to the Euler equations, by the local well-posedness result, the limit (in every regularity class) is continuous in time except at $t = 0$.
\begin{equation*}
\begin{cases}
\begin{aligned}
    & \|\frac{w - \frac{1}{\gamma - 1} c_0 + \frac{1}{2} v_0}{t}\|_{C_t^0((0, t^*]) H_u^{m + 1}} < \infty && \|L w\|_{C_t^0((0, t^*]) H_u^m} < \infty \\
    & \|\frac{\underline{w} - \frac{1}{\gamma - 1} c_0 - \frac{1}{2} v_0 + \frac{2}{\gamma + 1} u}{t}\|_{C_t^0((0, t^*]) H_u^{m + 1}} < \infty && \|L \underline{w}\|_{C_t^0((0, t^*]) H_u^m} < \infty \\
    & \|\frac{r - r_0 - (v_0 + c_0) t + t u}{t^2}\|_{C_t^0((0, t^*]) H_u^{m + 1}} < \infty && \|\frac{L r - (v_0 + c_0) + u}{t}\|_{C_t^0((0, t^*]) H_u^{m + 1}} < \infty \quad \forall m \geq 0
\end{aligned}
\end{cases}
\end{equation*}

Then, by removing the $t$ in the denominators and using the equality $L \underline{w} = - \frac{c v}{r}$, we get the continuity of all the quantities at $t = 0$ except $L w$.
\begin{equation*}
\begin{cases}
    w \in C_t^0([0, t^*]) H_u^{m + 1} \cap C_t^1((0, t^*]) H_u^m \\
    \underline{w} \in C_t^0([0, t^*]) H_u^{m + 1} \cap C_t^1([0, t^*]) H_u^m \\
    r \in C_t^1([0, t^*]) H_u^{m + 1} \qquad \qquad \qquad \qquad \forall m \geq 0 \\
    w, \underline{w}, r \in C_{t, u}^\infty((0,t^*] \times [0, u^*]) \\
\end{cases}
\end{equation*}

This finishes the construction of local rarefaction wave solutions.

We can control the higher order $t$-derivatives when $t \to 0$ by using the transport equations like in the subsection 5.2 of \cite{Luo-Yu-2}.

For $\underline{w}, r$, we use the following equations and commute $L, T$ if needed.
\begin{equation*}
    L \underline{w} = - \frac{c v}{r}, \qquad L r = v + c
\end{equation*}

For $w$, directly using the Euler equations will cause the loss of $t$-orders. We need to commute one more $L$ to obtain the transport equation for $L$-derivatives of $w$.
\begin{equation*}
\begin{aligned}
    & \underline{L} w = - \frac{\kappa v}{r} && \Longrightarrow L w = - \frac{c v}{r} - \frac{2 c}{\kappa} T w \qquad \mathrm{\ Loss\ of\ } t \mathrm{\ order} \\
    & L \underline{L} w = - L(\frac{\kappa v}{r}) && \Longrightarrow \underline{L} L w = L (\frac{\kappa}{c}) L w - L(\frac{\kappa v}{r}) \qquad \mathrm{\ Transport\ equation\ for\ } Lw
\end{aligned}
\end{equation*}

So we use the following three equations.
\begin{equation*}
\begin{cases}
    \underline{L} L w = (L (\frac{\kappa}{c}) + \frac{\kappa}{r}) L w - (\frac{\kappa}{r} L \underline{w} + L(\frac{\kappa}{r}) v) \\
    L \underline{w} = - \frac{c v}{r} \\
    L r = v + c
\end{cases}
\end{equation*}

By commuting $L, T$, we can obtain the boundedness and continuity of the higher $L$-derivatives in the following order.
\begin{equation*}
\begin{aligned}
    & L T^{\geq 0} r \in C_{t, u}^0 ([0, t^*] \times [0, u^*]), \quad L T^{\geq 0} \underline{w} \in C_{t, u}^0([0, t^*] \times [0, u^*]) & \Longrightarrow & L T^{\geq 0} w \in C_{t, u}^0([0, t^*] \times [0, u^*]) \\
    \Longrightarrow & L^2 T^{\geq 0} r \in C_{t, u}^0 ([0, t^*] \times [0, u^*]), \quad L^2 T^{\geq 0} \underline{w} \in C_{t, u}^0([0, t^*] \times [0, u^*]) & \Longrightarrow & L^2 T^{\geq 0} w \in C_{t, u}^0([0, t^*] \times [0, u^*]) \\
    \Longrightarrow & \cdots \\
\end{aligned}
\end{equation*}

Now we have obtained all the regularity results. We list them in the following.
\begin{equation*}
\begin{cases}
    w, \underline{w}, c, v, r \in C_{t, u}^\infty([0,t^*] \times [0, u^*]) \\
    w = \frac{1}{\gamma - 1} c_0 - \frac{1}{2} v_0 + O(t) \\
    \underline{w} = \frac{1}{\gamma - 1} c_0 + \frac{1}{2} v_0 - \frac{2}{\gamma + 1} u + O(t) \\
    c = c_0 - \frac{\gamma - 1}{\gamma + 1} u + O(t) \\
    v = v_0 - \frac{2}{\gamma + 1} u + O(t) \\
    r = r_0 + (v_0 + c_0)t - t u + O(t^2)
\end{cases}
\end{equation*}

We note that the low order parts are exactly the 1-D rarefaction wave solution.

\section{Global in Time Solution for Background Close to Constant State}

Now we turn to look at how far the rarefaction wave solution can extend. There are two directions: extend in $u$ to vacuum, or extend in $t$ to infinity. We now focus on extending in $t$, trying to find when the rarefaction wave solution exists globally in time.

This won't come for free for an arbitrary exterior background solution. We prove that we can extend the rarefaction wave solution globally in time for exterior solutions close to the constant state.

\begin{figure}[H]
\centering
\scalebox{0.6}{\begin{tikzpicture}

% Coordinates
\coordinate (origin) at (0,0);
\coordinate (r-axis) at (6,0);
\coordinate (t-axis) at (0,6);
\coordinate (singular) at (1,1);
\node at (r-axis) [anchor=north]{r};
\node at (t-axis) [anchor=east]{t};
\draw [thick,->] (origin) -- (r-axis);
\draw [thick,->] (origin) -- (t-axis);
\coordinate (rS) at (1,0);
\coordinate (tS) at (0,1);
\node at (rS) [anchor=north]{1};
\node at (tS) [anchor=east]{1};
\draw [dotted,semithick] (rS) -- (singular);
\draw [dotted,semithick] (tS) -- +(6,0);
% Figure
\draw [semithick] (singular) .. controls +(1,1.2) and +(-1,-1) .. +(4.5,4.5) node [anchor=south]{$C_0$};
\draw [dashed,semithick] (singular) .. controls +(-0.2,1.2) and +(-1+0.12,-1) .. +(2.7,4.5) node [anchor=south]{$C_{u^*}$};
\foreach \x in {1,...,5}
{\draw [dotted,semithick] (singular) .. controls +(1-0.2*\x,1.2) and +(-1+0.02^\x,-1) .. +(4.5-0.3*\x,4.5);};
\fill [gray,opacity=0.5] (singular) .. controls +(1,1.2) and +(-1,-1) .. +(4.5,4.5) -- (6,5.5) -- (6,1) -- cycle;
% Text
\node at (3.5,4.5) {$\kappa \sim \ln t$};
\node at (4,2) {Close to Constant};

\end{tikzpicture}}
\caption{Global Rarefaction}
\label{fig:global2}
\end{figure}
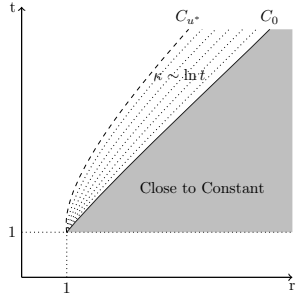

\begin{theorem}[Global Existence of Rarefaction Wave Solutions] \label{thm:global}
Without loss of generality, Assume $c_0 = 1, r_0 = 1$ and assume that the initial singularity lies at $(t, r) = (1, 1)$. If the exterior background solution is close to the constant state, which means that it satisfies \ref{eq:global-perturbation_of_constant}, then the rarefaction wave solution can be extended globally in time on $[0, u^*]$. See figure \ref{fig:global2}.
The rarefaction wave solution satisfies the following.

\begin{itemize}
    \item It satisfies the norm bounds \ref{eq:global-base-final_estimates} and \ref{eq:global-high-final_estimates}.

    \item $u^*$ depends on up to $4$ derivatives of the exterior solution.

    \item The rarefaction wave expands like $\ln t$, which means $\kappa \sim \ln t$.
\end{itemize}
\end{theorem}

\begin{remark}
The $\ln t$ rarefaction wave obtained in this result is similar to the recent result \cite{Wang-Global-Exterior} by Qian Wang. She proved the global existence of exterior solutions for 3-D compressible Euler equations without symmetry (close to spherical symmetry) and showed that if there is enough density drop in the initial data, there will be $\ln t$ rarefaction happening at null infinity.

The initial data is regular and without symmetry in \cite{Wang-Global-Exterior}, and is spherical symmetric but singular in our case. The similar rarefaction wave at long time evolved from different types of initial data may imply some generality of the $\ln t$ rarefaction at null infinity in 3-D.
\end{remark}

We first handle the most advantageous case we can have: the constant state background, $(v, c) \equiv (0, c_0)$, and the rarefaction wave originates from $r = r_0$.

\subsection{Computation of the Data on the Boundary $C_0$}

By the local result, we already have a smooth rarefaction solution in $(0, t^*] \times [0, u^*]$. By the local well-posedness of this problem, the only remaining thing is to follow the bootstrap scheme and close the energy estimates in the region $[t^*, T] \times [0, u^*]$ for arbitrarily large $T$.

The ansatz shall be changed, since the ansatz with positive powers of $t$ is used to deal with the singularity at $t = 0$, and is not suitable for large $t$.

The idea of the new ansatz will come from the computation of the data on the boundary $C_0$. \\

Recall that during the construction of the local solution, the data on $C_0$ is computed by solving ODE systems for $(T^n w, T^n \underline{w}, T^{n - 1} \kappa)$. It begins with \ref{eq:ODE-first} for the $n = 1$ case, then \ref{eq:ODE-general} for the general case. We copy the general ODE system here for reference. 
\begin{equation*}
\begin{cases}
    T^n w = - (\frac{v}{2 r} + \frac{Lw}{2 c}) T^{n - 1} \kappa + \sum\limits_{\substack{n_1 + n_2 = n - 1 \\ n_i \geq 0,\ n_2 \leq n - 2}} \mathrm{C}_{n_1, n_2} T^{n_1} (\frac{v}{r} + \frac{Lw}{c}) T^{n_2} \kappa \\
    L T^n \underline{w} = \frac{\gamma - 1}{2} (- 2 \frac{\underline{w}}{r} T^n \underline{w} + 2 \frac{w}{r} T^n w - \frac{\underline{w}^2 - w^2}{r^2} T^{n - 1} \kappa) \\
    \qquad \qquad + \sum\limits_{\substack{n_1 + n_2 = n - 1 \\ n_i \geq 0,\ n_2 \leq n - 2}} \mathrm{C}_{n_1, n_2} (- 2 T^{n_1} \frac{\underline{w}}{r} T^{n_2 + 1} \underline{w} + 2 T^{n_1} \frac{w}{r} T^{n_2 + 1} w - T^{n_1} \frac{\underline{w}^2 - w^2}{r^2} T^{n_2} \kappa) \\
    L T^{n - 1} \kappa = - \frac{\gamma + 1}{2} T^n \underline{w} + \frac{3 - \gamma}{2} T^n w
\end{cases}
\end{equation*}

The terms in the summations are lower order terms, which would have been known when we solve for $(T^n w, T^n \underline{w}, T^{n - 1} \kappa)$.

The initial data at $t = 0$ is set as in the local case, which comes from the well-known 1-D case. The constant state background gives $(c, v, r)$ on the whole $u = 0$ surface. They are listed as follows: 
\begin{equation*}
\begin{aligned}
    & t = 0 \mathrm{\ Data:} && u = 0 \mathrm{\ Data:} \\
    & \begin{cases}
        Tw = 0 \\
        T\underline{w} = - \frac{2}{\gamma + 1} \\
        \kappa = 0 \\
        T^n w = T^n \underline{w} = T^{n - 1} \kappa = 0 \ (n \geq 2)
    \end{cases} && \begin{cases}
        c = c_0 = \mathrm{Const.} \\
        v = 0 \\
        w = \underline{w} = \frac{1}{\gamma - 1} c_0 = \mathrm{Const.} \\
        r = r_0 + c_0 t
    \end{cases}
\end{aligned}
\end{equation*}

In the constant state background case, there are vanishing which make the solving process much easier and clearer. We have the following on $C_0$.
\begin{equation*}
    v = 0, \quad Lw = 0, \quad w = \underline{w} = \frac{1}{\gamma - 1} c_0
\end{equation*}

So the system simplifies to the following.
\begin{equation*}
\begin{cases}
    T^n w = \mathrm{L.O.T.} \\
    LT^n \underline{w} = - \frac{c_0}{r} T^n \underline{w} + \frac{c_0}{r} T^n w + \mathrm{L.O.T.} \\
    LT^{n - 1} \kappa = - \frac{\gamma + 1}{2} T^n \underline{w} + \frac{3 - \gamma}{2} T^n w
\end{cases}
\end{equation*}

Then we see that after the lower orders are solved, which means that all the L.O.T. in the system are known, solving the system can be split into three steps:

\begin{enumerate}
    \item First, directly compute $T^n w$ from the expression.
    
    \item Second, solve for $T^n \underline{w}$. Since $T^n w$ is also known, the equation for $T^n \underline{w}$ simplifies to an inhomogeneous linear ODE.
    \begin{equation*}
        LT^n \underline{w} = - \frac{c_0}{r} T^n \underline{w} + F \mathrm{\ (known\ terms)}
    \end{equation*}

    Then $T^n \underline{w}$ can be written in an explicit formula:
    \begin{equation*}
    \begin{aligned}
        T^n \underline{w}(t, 0) & = T^n \underline{w}(0, 0) e^{\int_0^t - \frac{c_0}{r(s)} ds} + \int_0^t F(\tau) e^{\int_\tau^t - \frac{c_0}{r(s)} ds} d\tau \\
        & = T^n \underline{w}(0, 0) e^{\int_0^t (- \ln(r)')(s) ds} + \int_0^t F(\tau) e^{\int_\tau^t (- \ln(r)')(s) ds} d\tau \\
        & = T^n \underline{w}(0, 0) \frac{r_0}{r(t)} + \int_0^t F(\tau) \frac{r(\tau)}{r(t)} d\tau
    \end{aligned}
    \end{equation*}

    \item Lastly, solve $T^{n - 1} \kappa$ by direct integration given by the following formula.
    \begin{equation*}
        T^{n - 1}\kappa(t, 0) = \int_0^t - \frac{\gamma + 1}{2} T^n \underline{w}(\tau) + \frac{3 - \gamma}{2} T^n w(\tau) d\tau
    \end{equation*}
\end{enumerate}

Now we try to solve the ODE systems explicitly. The first system gives the following.
\begin{equation*}
\begin{cases}
    Tw = 0 \\
    LT\underline{w} = - \frac{c_0}{r} T\underline{w} + \frac{c_0}{r} Tw\\
    L\kappa = - \frac{\gamma + 1}{2} T\underline{w} + \frac{3 - \gamma}{2} Tw
\end{cases} \Rightarrow \begin{cases}
    Tw = 0 \\
    T\underline{w} = - \frac{2}{\gamma + 1} \frac{r_0}{r} \\
    \kappa = \frac{r_0}{c_0} \ln(\frac{r}{r_0})
\end{cases}
\end{equation*}

Go on to the second system.
\begin{equation*}
\begin{cases}
\begin{aligned}
    T^2 w & = - \frac{\kappa Tv}{2 r} \\
    & = \frac{1}{\gamma + 1} \frac{1}{c_0} (\frac{r_0}{r})^2 \ln(\frac{r}{r_0})
\end{aligned} \\
\begin{aligned}
    LT^2 \underline{w} & = - \frac{c_0}{r} T^2 \underline{w} + \frac{c_0}{r} T^2 w - (\gamma - 1) \frac{(T\underline{w})^2}{r} - 2 (\gamma - 1) \frac{\underline{w} T\underline{w}}{r^2} \kappa \\
    & = - \frac{c_0}{r} T^2 \underline{w} + \frac{1}{\gamma + 1} \frac{1}{r_0} (\frac{r_0}{r})^3 \ln(\frac{r}{r_0}) - \frac{4 (\gamma - 1)}{(\gamma + 1)^2} \frac{1}{r_0} (\frac{r_0}{r})^3 + \frac{4}{\gamma + 1} \frac{1}{r_0} (\frac{r_0}{r})^3 \ln(\frac{r}{r_0}) \\
    & = - \frac{c_0}{r} T^2 \underline{w} + \frac{5}{\gamma + 1} \frac{1}{r_0} (\frac{r_0}{r})^3 \ln(\frac{r}{r_0}) - \frac{4 (\gamma - 1)}{(\gamma + 1)^2} \frac{1}{r_0} (\frac{r_0}{r})^3
\end{aligned} \\
\begin{aligned}
    LT\kappa & = - \frac{\gamma + 1}{2} T^2 \underline{w} + \frac{3 - \gamma}{2} T^2 w \\
    & = - \frac{\gamma + 1}{2} T^2 \underline{w} + \frac{3 - \gamma}{2 (\gamma + 1)} \frac{1}{c_0} (\frac{r_0}{r})^2 \ln(\frac{r}{r_0})
\end{aligned}
\end{cases}
\end{equation*}

Notice that all lower order terms appearing on the right hand side take the form $(\frac{r_0}{r})^a \ln^b (\frac{r}{r_0})$. This means that we need to calculate two types of integration involving such terms.
\begin{equation*}
\begin{cases}
    \int_0^t (\frac{r_0}{r(\tau)})^a \ln^b (\frac{r(\tau)}{r_0}) \frac{r(\tau)}{r(t)} d\tau = \frac{r_0}{r(t)} \int_0^t (\frac{r_0}{r(\tau)})^{a - 1} \ln^b (\frac{r(\tau)}{r_0}) d\tau, \quad \mathrm{from\ solving\ } T^n \underline{w} \\
    \int_0^t (\frac{r_0}{r(\tau)})^a \ln^b (\frac{r(\tau)}{r_0}) d\tau, \quad \mathrm{from\ solving\ } T^{n - 1} \kappa
\end{cases}
\end{equation*}

The integral $\int_0^t (\frac{r_0}{r(\tau)})^a \ln^b (\frac{r(\tau)}{r_0}) d\tau$ can be found using the following formula.
\begin{equation*}
\begin{aligned}
    & \frac{d}{dt}((\frac{r_0}{r})^{a - 1} \ln^b (\frac{r}{r_0})) = - (a - 1) \frac{c_0}{r_0} (\frac{r_0}{r})^a \ln^b (\frac{r}{r_0}) + b \frac{c_0}{r_0} (\frac{r_0}{r})^a \ln^{b - 1} (\frac{r}{r_0}) \\
    & (\frac{r_0}{r})^a \ln^b (\frac{r}{r_0}) = - \frac{1}{a - 1} \frac{r_0}{c_0} \frac{d}{dt}((\frac{r_0}{r})^{a - 1} \ln^b (\frac{r}{r_0})) + \frac{b}{a - 1} (\frac{r_0}{r})^a \ln^{b - 1} (\frac{r}{r_0})
\end{aligned}
\end{equation*}

\begin{equation*}
\begin{aligned}
    & (\frac{r_0}{r})^a \ln^b (\frac{r}{r_0}) = - \frac{r_0}{c_0} \sum\limits_{i = 0}^b \frac{b!}{(b - i)!} \frac{1}{(a - 1)^{i + 1}} \frac{d}{dt} ((\frac{r_0}{r})^{a - 1} \ln^{b - i} (\frac{r}{r_0})), \quad a \geq 2 \\
    & (\frac{r_0}{r}) \ln^b (\frac{r}{r_0}) = \frac{r_0}{c_0} \frac{1}{b + 1} \frac{d}{dt} \ln^{b + 1}(\frac{r}{r_0})
\end{aligned}
\end{equation*}

\begin{equation} \label{eq:integral_formula}
\begin{aligned}
    & \int_0^t (\frac{r_0}{r(\tau)})^a \ln^b (\frac{r(\tau)}{r_0}) d\tau = - \frac{r_0}{c_0} \sum\limits_{i = 0}^b \frac{b!}{(b - i)!} \frac{1}{(a - 1)^{i + 1}} (\frac{r_0}{r(t)})^{a - 1} \ln^{b - i} (\frac{r(t)}{r_0}) + \frac{r_0}{c_0} \frac{1}{(a - 1)^{b + 1}}, \quad a \geq 2 \\
    & \int_0^t (\frac{r_0}{r(\tau)}) \ln^b (\frac{r(\tau)}{r_0}) d\tau = \frac{r_0}{c_0} \frac{1}{b + 1} \ln^{b + 1}(\frac{r(t)}{r_0})
\end{aligned}
\end{equation}

So the solution for the second order system can be written out.
\begin{equation*}
\begin{cases}
\begin{aligned}
    T^2 w = \frac{1}{\gamma + 1} \frac{1}{c_0} (\frac{r_0}{r})^2 \ln(\frac{r}{r_0})
\end{aligned} \\
\begin{aligned}
    T^2 \underline{w} & = \frac{5}{\gamma + 1} \frac{1}{c_0} (- (\frac{r_0}{r})^2 \ln(\frac{r}{r_0}) - (\frac{r_0}{r})^2 + \frac{r_0}{r}) - \frac{4 (\gamma - 1)}{(\gamma + 1)^2} \frac{1}{c_0} (- (\frac{r_0}{r})^2 + \frac{r_0}{r}) \\
    & = \frac{1}{c_0} (\frac{\gamma + 9}{(\gamma + 1)^2} \frac{r_0}{r} - \frac{\gamma + 9}{(\gamma + 1)^2} (\frac{r_0}{r})^2 - \frac{5}{\gamma + 1} (\frac{r_0}{r})^2 \ln(\frac{r}{r_0}))
\end{aligned} \\
\begin{aligned}
    T \kappa = \frac{r_0}{c_0^2} (- \frac{\gamma + 9}{2 (\gamma + 1)} \ln(\frac{r}{r_0}) + \frac{5 \gamma + 17}{2(\gamma + 1)} (1 - \frac{r_0}{r}) - \frac{2 \gamma + 4}{\gamma + 1} \frac{r_0}{r} \ln(\frac{r}{r_0}))
\end{aligned}
\end{cases}
\end{equation*}

Notice that all the terms of $T^2 w, T^2 \underline{w}, T\kappa$ take the form $(\frac{r_0}{r})^a \ln^b (\frac{r}{r_0})$, and their leading terms (the terms which have worst decay when $t \to \infty$) are $(\frac{r_0}{r})^2 \ln(\frac{r}{r_0}), \frac{r_0}{r}, \ln(\frac{r}{r_0})$ respectively. We will prove inductively that these decay estimates hold for higher order terms $T^n w, T^n \underline{w}, T^{n - 1} \kappa$.

We define a terminology that will be used in the proof.

\begin{definition}[Function of Type $r_0^p c_0^q(a, b)$]
We say a function $f(t)$ is of type $r_0^p c_0^q(a, b)$, or $f \in r_0^p c_0^q (a, b)$ if it takes the following form.
\begin{equation*}
    f(t) = r_0^p c_0^q \sum\limits_{\substack{a', b',\ \mathrm{finitely\ many} \\ a' > a \mathrm{\ or\ } a' = a, b' \leq b}} R_{a', b'}(\gamma) (\frac{r_0}{r})^a \ln^b (\frac{r}{r_0})
\end{equation*}

$R_{a', b'}(\gamma)$ is a rational fraction of $\gamma$. This is saying that it can be controlled by $r_0^p c_0^q (\frac{r_0}{r})^a \ln^b (\frac{r}{r_0})$ in long time.
\end{definition}

We now prove the following proposition.

\begin{proposition}[Expression for $T^n w, T^n \underline{w}, T^{n - 1} \kappa$ Corresponding to Constant State Background]
In the case of the constant state background, the quantities take the following form.
\begin{equation*}
\begin{cases}
    T^n w \in \frac{1}{c_0^{n - 1}} (2, 1) \\
    T^n \underline{w} \in \frac{1}{c_0^{n - 1}} (1, 0) \\
    T^{n - 1} \kappa \in \frac{r_0}{c_0^n} (0, 1), \quad n \geq 1 \\
\end{cases}
\end{equation*}
\end{proposition}

Note that the powers of $r_0, c_0$ give the dimensions of these quantities.

Proof by induction.

The cases of $n = 1, 2$ have been verified above. Assume that the cases of $m$ with $m < n$ hold, now consider the case of $n$. We first write down the equations.
\begin{equation*}
\begin{cases}
    T^n w = \sum\limits_{\substack{n_1 + n_2 = n - 1 \\ n_i \geq 0,\ n_2 \leq n - 2}} \mathrm{C}_{n_1, n_2} T^{n_1} (\frac{v}{r} + \frac{Lw}{c}) T^{n_2} \kappa \\
    L T^n \underline{w} = - \frac{c_0}{r} T^n \underline{w} + \frac{c_0}{r} T^n w \\
    \qquad \qquad + \sum\limits_{\substack{n_1 + n_2 = n - 1 \\ n_i \geq 0,\ n_2 \leq n - 2}} \mathrm{C}_{n_1, n_2} (- 2 T^{n_1} \frac{\underline{w}}{r} T^{n_2 + 1} \underline{w} + 2 T^{n_1} \frac{w}{r} T^{n_2 + 1} w - T^{n_1} \frac{\underline{w}^2 - w^2}{r^2} T^{n_2} \kappa) \\
    L T^{n - 1} \kappa = - \frac{\gamma + 1}{2} T^n \underline{w} + \frac{3 - \gamma}{2} T^n w
\end{cases}
\end{equation*}

We have already removed the highest order terms vanishing from the properties of the constant state background.

We first prove some intermediate results for $m \geq 1$, which appear on the right hand side.
\begin{equation*}
    T^m w \in \frac{1}{c_0^{m - 1}} (2, 1) \implies L T^m w \in \frac{1}{c_0^{m - 1}} (3, 1)
\end{equation*}

\begin{equation*}
\begin{cases}
    c = \frac{\gamma - 1}{2} (\underline{w} + w) \\
    v = (\underline{w} - w) \\
    T^m w \in \frac{1}{c_0^{m - 1}} (2, 1) \\
    T^m \underline{w} \in \frac{1}{c_0^{m - 1}} (1, 0) \\
\end{cases} \implies \begin{cases}
    T^m c \in \frac{1}{c_0^{m - 1}} (1, 0) \\
    T^m v \in \frac{1}{c_0^{m - 1}} (1, 0) \\
\end{cases}
\end{equation*}

\begin{equation*}
    T^m(\frac{1}{r}) = \sum\limits_{\substack{m = m_1 + \cdots + m_p + p \\ m_i \geq 0}} \mathrm{C}_{m_1, \cdots, m_p} \frac{T^{m_1}\kappa \cdots T^{m_p}\kappa}{r^{p + 1}} \implies T^m (\frac{1}{r}) \in \frac{1}{r_0 c_0^m} (2, 1)
\end{equation*}

\begin{equation*}
    T^m (\frac{1}{c}) = \sum\limits_{\substack{m = m_1 + \cdots + m_p \\ m_i \geq 1}} \mathrm{C}_{m_1, \cdots, m_p} \frac{T^{m_1} c \cdots T^{m_p} c}{c^{p + 1}} \implies T^m (\frac{1}{c}) \in \frac{1}{c_0^{m + 1}} (1, 0)
\end{equation*}

We substitute these results into the equations to get the following.
\begin{equation*}
\begin{cases}
    \sum\limits_{\substack{n_1 + n_2 = n - 1 \\ n_i \geq 0,\ n_2 \leq n - 2}} \mathrm{C}_{n_1, n_2} T^{n_1} (\frac{v}{r}) T^{n_2} \kappa \in \frac{1}{c_0^{n - 1}} (2, 1) \\
    \sum\limits_{\substack{n_1 + n_2 = n - 1 \\ n_i \geq 0,\ n_2 \leq n - 2}} \mathrm{C}_{n_1, n_2} T^{n_1} (\frac{Lw}{c}) T^{n_2} \kappa \in \frac{1}{c_0^{n - 1}} (4, 2)
\end{cases} \implies T^n w \in \frac{1}{c_0^{n - 1}} (2, 1)
\end{equation*}

\begin{equation*}
\begin{cases}
    \sum\limits_{\substack{n_1 + n_2 = n - 1 \\ n_i \geq 0,\ n_2 \leq n - 2}} - \mathrm{C}_{n_1, n_2} 2 T^{n_1} \frac{\underline{w}}{r} T^{n_2 + 1} \underline{w} \in \frac{1}{r_0 c_0^{n - 2}} (3, 1) \\
    \sum\limits_{\substack{n_1 + n_2 = n - 1 \\ n_i \geq 0,\ n_2 \leq n - 2}} \mathrm{C}_{n_1, n_2} 2 T^{n_1} \frac{w}{r} T^{n_2 + 1} w \in \frac{1}{r_0 c_0^{n - 2}} (5, 2) \\
    \sum\limits_{\substack{n_1 + n_2 = n - 1 \\ n_i \geq 0,\ n_2 \leq n - 2}} - \mathrm{C}_{n_1, n_2} T^{n_1} \frac{\underline{w}^2 - w^2}{r^2} T^{n_2} \kappa \in \frac{1}{r_0 c_0^{n - 2}} (3, 1) \\
    \frac{c_0}{r} T^n w \in \frac{1}{r_0 c_0^{n - 2}} (3, 1) \\
\end{cases} \implies \begin{cases}
    L T^n \underline{w} = - \frac{c_0}{r} T^n \underline{w} + F \\
    F \in \frac{1}{r_0 c_0^{n - 2}} (3, 1)
\end{cases}
\end{equation*}

Now we use the integral formula \ref{eq:integral_formula}.
\begin{equation*}
\begin{aligned}
    T^n \underline{w}(t, 0) & = \int_0^t F(\tau) \frac{r(\tau)}{r(t)} d\tau \\
    & = \frac{r_0}{r(t)} \int_0^t \frac{1}{r_0 c_0^{n - 2}} (2, 1)(\tau) d\tau \\
    & = \frac{r_0}{r(t)} \frac{1}{c_0^{n - 1}} (0, 0)(t) \\
    & \in \frac{1}{c_0^{n - 1}} (1, 0)
\end{aligned}
\end{equation*}

The integration of $(2, 1)$ becomes $(0, 0)$ because there are constant boundary terms coming from $t = 0$.

Finally, use the integration formula for $T^{n - 1} \kappa$ to obtain the expression for $T^{n - 1} \kappa$.
\begin{equation*}
\begin{aligned}
    T^{n - 1}\kappa(t, 0) & = \int_0^t - \frac{\gamma + 1}{2} T^n \underline{w}(\tau) + \frac{3 - \gamma}{2} T^n w(\tau) d\tau \\
    & = \int_0^t \frac{1}{c_0^{n - 1}} (1, 0)(\tau) d\tau \\
    & \in \frac{r_0}{c_0^n} (0, 1)
\end{aligned}
\end{equation*}

This finishes the induction.

So the size of the boundary data is completely clear. We now summarize the result for $u = 0$ data.

\begin{proposition}[Rarefaction Wave Data Corresponding to Constant State Background]
We have the following decay estimates for the quantities on $C_0$. $n \geq 2$, and $\lesssim$ means up to a constant depending on $\gamma, n$.
\begin{equation*}
\begin{aligned}
\begin{cases}
    c = c_0 = \mathrm{Const.} \\
    v = 0 \\
    w = \underline{w} = \frac{1}{\gamma - 1} c_0 = \mathrm{Const.} \\
    Lw = L\underline{w} = 0 \\
    r = r_0 + c_0 t
\end{cases} & \begin{cases}
    Tw = 0 \\
    T\underline{w} = - \frac{2}{\gamma + 1} \frac{r_0}{r} \\
    LTw = 0 \\
    LT\underline{w} = \frac{2}{\gamma + 1} \frac{c_0}{r_0} (\frac{r_0}{r})^2 \\
    \kappa = \frac{r_0}{c_0} \ln(\frac{r}{r_0})
\end{cases} & \begin{cases}
    |T^n w| \lesssim \frac{1}{c_0^{n - 1}} (\frac{r_0}{r})^2 \ln(\frac{r}{r_0}) \\
    |T^n \underline{w}| \lesssim \frac{1}{c_0^{n - 1}} \frac{r_0}{r} \\
    |LT^n w| \lesssim \frac{1}{c_0^{n - 2} r_0} (\frac{r_0}{r})^3 \ln(\frac{r}{r_0}) \\
    |LT^n \underline{w}| \lesssim \frac{1}{c_0^{n - 2} r_0} (\frac{r_0}{r})^2 \\
    |T^{n - 1}\kappa| \lesssim \frac{r_0}{c_0^n} \ln(\frac{r}{r_0})
\end{cases}
\end{aligned}
\end{equation*}

If we apply a linear transformation to the coordinate functions, we can remove all the dimensional constants in the expressions.
\begin{equation} \label{eq:dimension_transform}
    (t, x, u) \mapsto (\frac{c_0}{r_0} t + 1, \frac{1}{r_0} x, \frac{1}{c_0}u)
\end{equation}

In the new coordinates, $c_0 = r_0 = 1$, and the initial singular point lies at $(t, u) = (1, 1)$. The $r_0^p c_0^q (\frac{r_0}{r})^a \ln^b (\frac{r}{r_0})$ terms become $\frac{1}{t^a} \ln^b t$. And the data corresponding to the constant state background have the following bounds.
\begin{equation*}
\begin{aligned}
\begin{cases}
    c = 1 \\
    v = 0 \\
    w = \underline{w} = \frac{1}{\gamma - 1} \\
    Lw = L\underline{w} = 0 \\
    r = t
\end{cases} & \begin{cases}
    Tw = 0 \\
    T\underline{w} = - \frac{2}{\gamma + 1} \frac{1}{t} \\
    LTw = 0 \\
    LT\underline{w} = \frac{2}{\gamma + 1} \frac{1}{t^2} \\
    \kappa = \ln t
\end{cases} & \begin{cases}
    |T^n w| \lesssim \frac{1}{t^2} \ln t \\
    |T^n \underline{w}| \lesssim \frac{1}{t} \\
    |LT^n w| \lesssim \frac{1}{t^3} \ln t \\
    |LT^n \underline{w}| \lesssim \frac{1}{t^2} \\
    |T^{n - 1}\kappa| \lesssim \ln t \\
\end{cases}
\end{aligned}
\end{equation*}
\end{proposition}

\subsection{Energy Estimates Uniform in $T$}

In this section, we prove an energy estimate result that doesn't rely on the specific expression of the $C_0$ boundary data corresponding to the constant state background. We only need the data to satisfy certain boundedness and decay properties.

We adopt the same coordinate transformation as \ref{eq:dimension_transform}. But now we choose $c_0 = \lim\limits_{|x| \to \infty} c$.

Now we can always assume $t \geq 1$, which simplifies the notation since we don't need to introduce $t + 1$ or $\langle t \rangle$.

Since the behavior of the quantities when $t \to \infty$ differs greatly from the local case where we care about $t \to 0$, the multipliers need to be changed. From the rarefaction wave data we explicitly computed above, we can also see that $w, \underline{w}$ have different decay rates, so we will choose different multipliers for them to capture this information.

Now we choose the multipliers to be the following.
\begin{equation*}
    X_w = t^{1 + s} \frac{c}{\kappa} L + \underline{L}, \quad X_{\underline{w}} = \frac{\kappa}{c} L + \underline{L}, \quad (0 < s < 1)
\end{equation*}

We still use the commutator $T$.

The energy estimate now takes the following form.
\begin{equation*}
\begin{cases}
    \mathcal{E}_{T^N}^{[0, u^*]}(T) - \mathcal{E}_{T^N}^{[0, u^*]}(t^*) + \mathcal{F}_{T^N}^{[t^*, T]}(u^*) - \mathcal{F}_{T^N}^{[t^*, T]}(0) = \int_{\{t^* \leq t \leq T, 0 \leq u \leq u^*\}} Q_{T^N} d\mu_g \\
    \mathcal{E}_{T^N}^{[0, u^*]}(t) = \frac{1}{2} \int_{[0, u^*]} \int_{S_{t, u}} t^{1 + s} |LT^N w|^2 + |\underline{L}T^N w|^2 + \frac{\kappa^2}{c^2} |L T^N \underline{w}|^2 + |\underline{L}T^N \underline{w}|^2 d\mu_\slashed{g} du \\
    \mathcal{F}_{T^N}^{[t^*, T]}(u) = \int_{[t^*, T]} \int_{S_{t, u}} t^{1 + s} \frac{c}{\kappa} |LT^N w|^2 + \frac{\kappa}{c} |L T^N \underline{w}|^2 d\mu_\slashed{g} dt \\
    Q_{T^N} = Q_{0, T^N} + Q_{1, T^N} + Q_{2, T^N} \\
    Q_{0, T^N} = - \tilde{\rho}_{T^N w} (t^{1 + s} \frac{c}{\kappa} LT^N w + \underline{L}T^N w) - \tilde{\rho}_{T^N \underline{w}} (\frac{\kappa}{c} L T^N \underline{w} + \underline{L}T^N \underline{w}) \\
    Q_{1, T^N} = \frac{1}{2 \mu} ((1 + s) t^s + t^{1 + s} \underline{L} (\frac{c}{\kappa}) - t^{1 + s} \frac{c}{\kappa} L (\frac{\kappa}{c})) |LT^N w|^2 + \frac{1}{\mu} T(\frac{\kappa}{c}) |L T^N \underline{w}|^2 \\
    Q_{2, T^N} = - \frac{1}{\mu} (t^{1 + s} \frac{c}{\kappa} \frac{v + c}{r} + \frac{\kappa}{c} \frac{v - c}{r}) LT^N w \underline{L}T^N w - \frac{1}{\mu} (\frac{\kappa}{c} \frac{v + c}{r} + \frac{\kappa}{c} \frac{v - c}{r}) LT^N \underline{w} \underline{L}T^N \underline{w}
\end{cases}
\end{equation*}

Just like in the local case, we prove an energy estimate for $N = 3$, then prove the propagation of regularity.

We first list our assumptions on the initial data on $C_0$.
\begin{equation} \label{eq:global-C0_data}
\begin{cases}
\begin{aligned}
    & |c - 1| \leq \frac{\epsilon}{t} && |v| \leq \frac{\epsilon}{t} && |r - t| \leq \epsilon \ln t \\
    & |L w| \leq \frac{\epsilon}{t^{\frac{3 + s}{2}}} && |L \underline{w}| \leq \frac{\epsilon}{t^2} && |L \kappa - \frac{1}{t}| \leq \frac{\epsilon}{t} \\
    & |T w| \leq \frac{\epsilon}{t^{\frac{3 + s}{2}}} && |T \underline{w} + \frac{2}{\gamma + 1} \frac{1}{t}| \leq \frac{\epsilon}{t} && |\kappa - \ln t| \leq \epsilon \\
    & |L T^n w| \leq \frac{M_n}{t^{\frac{3 + s}{2}}} && |L T^n \underline{w}| \leq \frac{M_n}{t^2} && |L T^n \kappa| \leq \frac{M_n}{t} \\
    & |T^{n + 1} w| \leq \frac{M_n}{t^{\frac{3 + s}{2}}} && |T^{n + 1} \underline{w}| \leq \frac{M_n}{t} && |T^n \kappa| \leq M_n \ln t
\end{aligned} \\
    \mathcal{F}_{T^N}^{[t^*, T]}(0) \leq F_{\mathrm{Data}, N}
\end{cases}
\end{equation}

Where $n \geq 1, N \geq 1$. The assumption of finite flux is stronger in some sense, as the infinite integration domain requires roughly $t^{- \frac{1}{2}}$ more decay. It does not come from $L^\infty$-assumptions as in the local case.

These assumptions essentially say that the background solution is close to constant state.

Now we state the propositions.

\begin{proposition}[Base Energy Estimates] \label{prop:global-base_energy}
Let $(w, \underline{w}, r)$ defined on $(t, u) \in [t^*, T] \times [0, u^*], 1 < t^* \leq 2$ be a solution to the Euler equations. They belong to the following regularity class (the same as in the local case \ref{eq:local-base-regularity_class}).
\begin{equation*}
\begin{cases}
    w \in C_t^0 H_u^{4} \cap C_t^1 H_u^3 \\
    \underline{w} \in C_t^0 H_u^4 \cap C_t^1 H_u^3 \\
    r \in C_t^1 H_u^4
\end{cases}
\end{equation*}

The data on $C_0$ satisfy \ref{eq:global-C0_data}. We also assume the following $L^2$ and $L^\infty$-bounds for the data on $\Sigma_{t^*}$.
\begin{equation*}
\begin{cases}
    \mathcal{E}_{T^3}^{[0, u^*]}(t^*) \leq E_{\mathrm{Data}, 3} \\
    \int_0^{u^*} \int_{S_{t^*, u}} |T^3 \kappa|^2 d\mu_{\slashed{g}} du \leq E_{\mathrm{Data}, 3} \\
    \begin{aligned}
        & |c(t^*, u) - 1| \leq \frac{1}{2} && |v(t^*, u)| \leq \frac{1}{2} && |r(t^*, u) - t^*| \leq \frac{1}{2} t^* \\
        & |L \kappa(t^*, u) - \frac{1}{t^*}| \leq \frac{1}{2 t^*} && |\kappa(t^*, u) - \ln t^*| \leq \frac{1}{2} \ln t^*
    \end{aligned}
\end{cases}
\end{equation*}

Then there exists a constant $\mathrm{C}_3$ depending on $t^*, M_1, M_2, E_{\mathrm{Data}, 3}, F_{\mathrm{Data}, 3}$ such that if we have $u^*$ small enough:
\begin{equation*}
    u^* \leq \frac{1}{\mathrm{C}_3}
\end{equation*}

Then we can obtain uniform (in $T$) bounds, where $n = 1, 2$.
\begin{equation} \label{eq:global-base-final_estimates}
\begin{cases}
    \begin{aligned}
        & \sup\limits_{t \in [t^*, T]} \mathcal{E}_{T^3}^{[0, u^*]}(t) \leq E_3' && \sup\limits_{u \in [0, u^*]} \mathcal{F}_{T^3}^{[t^*, T]}(u) \leq E_3' \\
        & \int_0^{u^*} \int_{S_{t, u}} |T^3 \kappa|^2 d\mu_{\slashed{g}} du \leq E_3' t^2 \ln^2 t \\
        & \int_0^{u^*} \int_{S_{t, u}} |T^4 w|^2 d\mu_{\slashed{g}} du \leq E_3' \frac{\ln^2 t}{t^{1 + s}} && |L T^3 \underline{w}| \leq \frac{M_2'}{t^2} \\
    \end{aligned} \\
    \begin{aligned}
        & |c - 1| \leq \frac{1}{2} && |v| \leq \frac{1}{2 t} && |r - t| \leq \frac{1}{2} t \\
        & |L w| \leq \frac{M_0'}{t^{\frac{3 + s}{2}}} && |L \underline{w}| \leq \frac{M_0'}{t^2} && |L \kappa - \frac{1}{t}| \leq \frac{1}{2 t} \\
    & |T w| \leq \frac{M_0'}{t^{\frac{3 + s}{2}}} \ln t && |T \underline{w} + \frac{2}{\gamma + 1} \frac{1}{t}| \leq \frac{1}{\gamma + 1} \frac{1}{t} && |\kappa - \ln t| \leq \frac{1}{2} \ln t \\
    & |L T^n w| \leq \frac{M_n'}{t^{\frac{3 + s}{2}}} && |L T^n \underline{w}| \leq \frac{M_n'}{t^2} && |L T^n \kappa| \leq \frac{M_n'}{t} \\
    & |T^{n + 1} w| \leq \frac{M_n'}{t^{\frac{3 + s}{2}}} \ln t && |T^{n + 1} \underline{w}| \leq \frac{M_n'}{t} && |T^n \kappa| \leq M_n' \ln t
    \end{aligned}
\end{cases}
\end{equation}
\end{proposition}

\begin{proposition}[Propagation of Regularity] \label{prop:global-regularity}
Let $N \geq 4$, $(w, \underline{w}, r)$ defined on $(t, u) \in [t^*, T] \times [0, u^*], t^* > 1$ be the solution in proposition \ref{prop:global-base_energy}. Assume that the data on $\Sigma_{t^*}$ have higher regularity and is compatible with the data on $C_0$, satisfying the following bounds.
\begin{equation*}
\begin{cases}
    \mathcal{E}_{T^N}^{[0, u^*]}(t^*) \leq E_{\mathrm{Data}, N} \\
    \int_0^{u^*} \int_{S_{t^*, u}} |T^N \kappa|^2 d\mu_{\slashed{g}} du \leq E_{\mathrm{Data}, N}
\end{cases}
\end{equation*}

And assume that we have already obtained the lower order $L^2$ and $L^\infty$-bounds, where $1 \leq n \leq N - 2$.
\begin{equation*}
\begin{cases}
    \begin{aligned}
        & \sup\limits_{t \in [t^*, T]} \mathcal{E}_{T^{N - 1}}^{[0, u^*]}(t) \leq E_{N - 1}' && \sup\limits_{u \in [0, u^*]} \mathcal{F}_{T^{N - 1}}^{[t^*, T]}(u) \leq E_{N - 1}' \\
        & \int_0^{u^*} \int_{S_{t, u}} |T^{N - 1} \kappa|^2 d\mu_{\slashed{g}} du \leq E_{N - 1}' t^2 \ln^2 t \\
        & \int_0^{u^*} \int_{S_{t, u}} |T^N w|^2 d\mu_{\slashed{g}} du \leq E_{N - 1}' \frac{\ln^2 t}{t^{1 + s}} && |L T^{N - 1} \underline{w}| \leq \frac{M_{N - 2}'}{t^2} \\
    \end{aligned} \\
    \begin{aligned}
        & |L T^n w| \leq \frac{M_n'}{t^{\frac{3 + s}{2}}} && |L T^n \underline{w}| \leq \frac{M_n'}{t^2} && |L T^n \kappa| \leq \frac{M_n'}{t} \\
        & |T^{n + 1} w| \leq \frac{M_n'}{t^{\frac{3 + s}{2}}} \ln t && |T^{n + 1} \underline{w}| \leq \frac{M_n'}{t} && |T^n \kappa| \leq M_n' \ln t
    \end{aligned}
\end{cases}
\end{equation*}

Then we can propagate the regularity and obtain higher order $L^2$ and $L^\infty$-bounds.
\begin{equation} \label{eq:global-high-final_estimates}
\begin{cases}
    \begin{aligned}
        & \sup\limits_{t \in [t^*, T]} \mathcal{E}_{T^N}^{[0, u^*]}(t) \leq E_N' && \sup\limits_{u \in [0, u^*]} \mathcal{F}_{T^N}^{[t^*, T]}(u) \leq E_N' \\
        & \int_0^{u^*} \int_{S_{t, u}} |T^N \kappa|^2 d\mu_{\slashed{g}} du \leq E_N' t^2 \ln^2 t \\
        & \int_0^{u^*} \int_{S_{t, u}} |T^{N + 1} w|^2 d\mu_{\slashed{g}} du \leq E_N' \frac{\ln^2 t}{t^{1 + s}} && |L T^N \underline{w}| \leq \frac{M_{N - 1}'}{t^2} \\
    \end{aligned} \\
    \begin{aligned}
        & |L T^{N - 1} w| \leq \frac{M_{N - 1}'}{t^{\frac{3 + s}{2}}} && |L T^{N - 1} \underline{w}| \leq \frac{M_{N - 1}'}{t^2} && |L T^{N - 1} \kappa| \leq \frac{M_{N - 1}'}{t} \\
        & |T^N w| \leq \frac{M_{N - 1}'}{t^{\frac{3 + s}{2}}} \ln t && |T^N \underline{w}| \leq \frac{M_{N - 1}'}{t} && |T^{N - 1} \kappa| \leq M_{N - 1}' \ln t
    \end{aligned}
\end{cases}
\end{equation}
\end{proposition}

\subsection{Proof of Proposition \ref{prop:global-base_energy}}

Define two quantities to simplify the expressions.
\begin{equation*}
    M_{\leq 2} = \max\{M_1, M_2\}, \quad M = \max\{M_1, M_2, \sqrt{E_{\mathrm{Data}, 3}}, \sqrt{F_{\mathrm{Data}, 1}}, \sqrt{F_{\mathrm{Data}, 2}}, \sqrt{F_{\mathrm{Data}, 3}}\}
\end{equation*}

Similar to the proof \ref{subsec:local-proof-local_base_energy} in the local case, we use the $L^2$-bound of the initial data on $\Sigma_{t^*}$, integrate from $C_0$ to obtain the following $L^\infty$-bounds on $\Sigma_{t^*}$.
\begin{equation} \label{eq:global-base-Linfty_t*}
\begin{cases}
\begin{aligned}
    & |T^2 \kappa(t^*, u)| \lesssim_{t^*} (M_{\leq 2} + \sqrt{E_{\mathrm{Data}, 3}}) && \\
    & |L T^2 w(t^*, u*)| \lesssim_{t^*} M_{\leq 2} + \sqrt{E_{\mathrm{Data}, 3}} && |L T^2 \underline{w}(t^*, u)| \lesssim_{t^*} M_{\leq 2} + \sqrt{E_{\mathrm{Data}, 3}} \\
    & |T^3 w(t^*, u)| \lesssim_{t^*} M_{\leq 2} + \sqrt{E_{\mathrm{Data}, 3}} && |T^3 \underline{w}(t^*, u)| \lesssim_{t^*} M_{\leq 2} + \sqrt{E_{\mathrm{Data}, 3}}
\end{aligned}
\end{cases}
\end{equation}

Then we write down the bootstrap assumptions.
\begin{equation} \label{eq:global-base-bootstrap_assumptions}
\begin{cases}
    |L T^2 w(t, u)| \leq \frac{\mathrm{C}}{t^{\frac{3 + s}{2}}} \\
    |T^3 w(t, u)| \leq \frac{\mathrm{C}}{t} \\
    |T^3 \underline{w}(t, u)| \leq \frac{\mathrm{C}}{t}
\end{cases}
\end{equation}

Where $\mathrm{C} > M$ is a large constant depending on $t^*, M$, chosen such that the bounds \ref{eq:global-base-Linfty_t*} are stronger than \ref{eq:global-base-bootstrap_assumptions} restricted on $t = t^*$. This is to ensure that we have a valid place to start.

\begin{remark}
This bootstrap assumption seems to be weaker than the decay rate we claim to have, and also incomplete. This is because we are only using the energy estimate to obtain part of the $L^\infty$-estimates (the ones listed in the bootstrap assumption), the remaining part will be obtained using the first order Euler equations and the propagation equation for $\kappa$.
\end{remark}

Then we integrate from $C_0$ to obtain the following estimates.
\begin{equation*}
\begin{cases}
\begin{aligned}
    & |L T w| \lesssim \frac{\mathrm{C}}{t^{\frac{3 + s}{2}}} && |T^2 w| \lesssim \frac{\mathrm{C}}{t} && |T^2 \underline{w}| \lesssim \frac{\mathrm{C}}{t} \\
    & |L w| \lesssim \frac{\epsilon + \mathrm{C} u}{t^{\frac{3 + s}{2}}} && |T w| \lesssim \frac{\epsilon + \mathrm{C} u}{t} && |T \underline{w} + \frac{2}{\gamma + 1} \frac{1}{t}| \lesssim \frac{\epsilon + \mathrm{C} u}{t} \\
    & && |T c + \frac{\gamma - 1}{\gamma + 1} \frac{1}{t}| \lesssim \frac{\epsilon + \mathrm{C} u}{t} && |T v + \frac{2}{\gamma + 1} \frac{1}{t}| \lesssim \frac{\epsilon + \mathrm{C} u}{t} \\
    & && |c - 1| \lesssim \frac{\epsilon + (1 + \mathrm{C} u) u}{t} && |v| \lesssim \frac{\epsilon + (1 + \mathrm{C} u) u}{t}
\end{aligned}
\end{cases}
\end{equation*}

Then we use $L T^n \kappa = - \frac{\gamma + 1}{2} T^{n + 1} \underline{w} + \frac{3 - \gamma}{2} T^{n + 1} w$ to obtain bounds for $L T^n \kappa$, integrate from $\Sigma_{t^*}$ to obtain the bound for $T^2 \kappa$, and integrate from $C_0$ to obtain the bounds for $T \kappa, \kappa$.
\begin{equation*}
\begin{cases}
\begin{aligned}
    & |L T^2 \kappa| \lesssim \frac{\mathrm{C}}{t} && |T^2 \kappa| \lesssim_{t^*} \mathrm{C} \ln t \\
    & |L T \kappa| \lesssim \frac{\mathrm{C}}{t} && |T \kappa| \lesssim_{t^*} \mathrm{C} \ln t \\
    & |L \kappa - \frac{1}{t}| \lesssim \frac{\epsilon + \mathrm{C} u}{t} && |\kappa - \ln t| \lesssim_{t^*} (\epsilon + \mathrm{C} u) \ln t \\
    & && |r - t| \lesssim_{t^*} (\epsilon + (1 + \mathrm{C} u) u) \ln t
\end{aligned}
\end{cases}
\end{equation*}

We place a smallness requirement on $u^*$.
\begin{equation} \label{eq:global-base-smallness_sim}
    u^* \leq \frac{1}{\mathrm{C}_3^{(1)}(t^*, \mathrm{C})} \leq \frac{1}{\mathrm{C}}
\end{equation}

This is to ensure that the $L^\infty$-bounds above imply the following $\sim$ relations.
\begin{equation*}
\begin{cases}
\begin{aligned}
    & |c - 1| \leq \frac{1}{2} && |v| \leq \frac{1}{2 t} && |r - t| \leq \frac{1}{2} t \\
    & |L \kappa - \frac{1}{t} | \leq \frac{1}{2 t} && |\kappa - \ln t| \leq \frac{1}{2} \ln t && |T \underline{w} + \frac{2}{\gamma + 1} \frac{1}{t}| \leq \frac{1}{\gamma + 1} \frac{1}{t} \\
\end{aligned}
\end{cases}
\end{equation*}

Now we use the first order Euler equations to obtain the remaining $L^\infty$-bounds. We begin with $L T^2 \underline{w}$.
\begin{equation*}
\begin{aligned}
    L T^2 \underline{w} & = \frac{\gamma - 1}{2} (- 2 \frac{\underline{w}}{r} T^2 \underline{w} + 2 \frac{w}{r} T^2 w - \frac{\underline{w}^2 - w^2}{r^2} T \kappa) \\
    & \qquad + \frac{\gamma - 1}{2} (- 2 T \frac{\underline{w}}{r} T \underline{w} + 2 T \frac{w}{r} T w - T \frac{\underline{w}^2 - w^2}{r^2} \kappa) \\
    (\mathrm{take}\ |\cdot|) & \lesssim_{t^*} \frac{\mathrm{C}}{t^2} \\
\end{aligned}
\end{equation*}

Then we integrate from $C_0$ to get lower order bounds.
\begin{equation*}
|L T \underline{w}| \lesssim_{t^*} \frac{\mathrm{C}}{t^2}, \quad |L \underline{w}| \lesssim_{t^*} \frac{1}{t^2}
\end{equation*}

Now we deal with $T^3 w$.
\begin{equation*}
\begin{aligned}
    T^3 w & = - (\frac{v}{2 r} + \frac{Lw}{2 c}) T^2 \kappa - T^2 (\frac{v}{2 r} + \frac{Lw}{2 c}) \kappa - T (\frac{v}{r} + \frac{Lw}{c}) T \kappa \\
    (\mathrm{take}\ |\cdot|) & \lesssim_{t^*} \frac{\mathrm{C}^2}{t^{\frac{3 + s}{2}}} \ln t \\
\end{aligned}
\end{equation*}

Then we integrate from $C_0$ to get lower order bounds.
\begin{equation*}
    |T^2 w| \lesssim_{t^*} \frac{\mathrm{C}}{t^{\frac{3 + s}{2}}} \ln t, \quad |T w| \lesssim_{t^*} \frac{1}{t^{\frac{3 + s}{2}}} \ln t
\end{equation*}

Now we have obtained all the $L^\infty$-bounds in the result of the proposition \ref{prop:global-base_energy}, the remaining work is to prove a uniform (in $T$) energy estimate and improve the bootstrap assumption. We will no longer track the powers of $\mathrm{C}$ in the following proof.

Before we begin, we first list two admissible terms, which will appear many times in the estimates.
\begin{equation*}
\begin{aligned}
    L T^3 \underline{w} & = - T^3 (\frac{c v}{r}) \\
    (\mathrm{take}\ |\cdot|) & \lesssim_{t^*, \mathrm{C}} \frac{1}{t^2}
\end{aligned}
\end{equation*}

\begin{equation*}
\begin{aligned}
    & \quad\  \int_{t, u, S_{t, u}} \frac{1}{t^{3 + \epsilon}} |T^3 \kappa|^2 d\mu_{\slashed{g}} du dt \\
    & = \int_{t, u} \frac{4 \pi r^2}{t^{3 + \epsilon}} |T^3 \kappa|^2 du dt \\
    & \lesssim \int_u |T^3 \kappa(t^*, u)|^2 \int_t \frac{1}{t^{1 + \epsilon}} dt du + \int_{t, u} \frac{4 \pi r^2}{t^{3 + \epsilon}} (\int_{t^*}^t |T^4 w(\tau, u)| + |T^4 \underline{w}(\tau, u)| d\tau)^2 du dt \\
    & \lesssim_{t^*} E_{\mathrm{Data}, 3} + \int_{t, u} \frac{4 \pi r^2}{t^{3 + \epsilon}} (\int_{t^*}^t |T^4 w(\tau, u)| + |T^4 \underline{w}(\tau, u)| d\tau)^2 du dt \\
    & \lesssim E_{\mathrm{Data}, 3} + \int_{t, u} \frac{4 \pi r^2}{t^{3 + \epsilon}} (\int_{t^*}^t \tau (|T^4 w(\tau, u)|^2 + |T^4 \underline{w}(\tau, u)|^2) d\tau) \int_{t^*}^t \frac{1}{\tau} d\tau du dt \\
    & \lesssim_{t^*} E_{\mathrm{Data}, 3} + \int_{\tau, u} \tau (|T^4 w(\tau, u)|^2 + |T^4 \underline{w}(\tau, u)|^2) \int_\tau^T \frac{\ln t}{t^{1 + \epsilon}} dt du d\tau \\
    & \lesssim E_{\mathrm{Data}, 3} + \int_{\tau, u} \tau^{1 - \epsilon} \ln \tau (|T^4 w(\tau, u)|^2 + |T^4 \underline{w}(\tau, u)|^2) du d\tau \\
    & \lesssim E_{\mathrm{Data}, 3} + \int_t \frac{\ln t}{t^{1 + \epsilon}} \mathcal{E}_{T^3}^{[0, u^*]}(t) dt \\
    & \lesssim_{\mathrm{C}} 1 + \int_t \frac{\ln t}{t^{1 + \epsilon}} \mathcal{E}_{T^3}^{[0, u^*]}(t) dt \\
\end{aligned}
\end{equation*}

\begin{equation*}
\begin{aligned}
    T^4 w & = - (\frac{v}{2 r} + \frac{Lw}{2 c}) T^3 \kappa - T^3 (\frac{v}{2 r} + \frac{Lw}{2 c}) \kappa - T^2 (\frac{3 v}{2 r} + \frac{3 Lw}{2 c}) T \kappa - T (\frac{3 v}{2 r} + \frac{3 Lw}{2 c}) T^2 \kappa \\
    (\mathrm{take}\ |\cdot|) & \lesssim_{t^*, \mathrm{C}} \frac{1}{t^{\frac{3 + s}{2}}} \ln t + \frac{1}{t^{\frac{3 + s}{2}}} |T^3 \kappa| + \ln t |L T^3 w| \\
\end{aligned}
\end{equation*}

\begin{equation*}
\begin{aligned}
    & \quad\  \int_{t, u, S_{t, u}} t^{s - \epsilon} |T^4 w|^2 d\mu_{\slashed{g}} du dt \\
    (\mathrm{take}\ |\cdot|)& \lesssim_{t^*, \mathrm{C}} \int_{t, u, S_{t, u}} \frac{1}{t^{3 + \epsilon}} \ln^2 t + \frac{1}{t^{3 + \epsilon}} |T^3 \kappa|^2 + t^{s - \epsilon} \ln^2 t |L T^3 w|^2 d\mu_{\slashed{g}} du dt \\
    & \lesssim_{t^*} E_{\mathrm{Data}, 3} + 1 + \int_t \frac{\ln t}{t^{1 + \epsilon}} \mathcal{E}_{T^3}^{[0, u^*]}(t) dt \\
    & \lesssim_{\mathrm{C}} 1 + \int_t \frac{\ln t}{t^{1 + \epsilon}} \mathcal{E}_{T^3}^{[0, u^*]}(t) dt \\
\end{aligned}
\end{equation*}

The same estimate holds for $\underline{L} T^3 w$.

We deal with $\int Q_{T^3}$ term by term.
\begin{equation*}
\begin{aligned}
    \int_{t, u, S_{t, u}} Q_{1, T^3} d\mu_g & = \int_{t, u, S_{t, u}} \frac{1}{2} ((1 + s) t^s + t^{1 + s} \underline{L} (\frac{c}{\kappa}) - t^{1 + s} \frac{c}{\kappa} L (\frac{\kappa}{c})) |LT^3 w|^2 + T(\frac{\kappa}{c}) |L T^3 \underline{w}|^2 d\mu_{\slashed{g}} du dt \\
    (\mathrm{take}\ |\cdot|) & \lesssim_{t^*, \mathrm{C}} \int_{t, u, S_{t, u}} \frac{t^{1 + s}}{\ln t} |L T^3 w|^2 + \ln t |L T^3 \underline{w}|^2 d\mu_{\slashed{g}} du dt \\
    & \lesssim \int_u \mathcal{F}_{T^3}^{[t^*, T]}(u) du
\end{aligned}
\end{equation*}

\begin{equation*}
\begin{aligned}
    \int_{t, u, S_{t, u}} Q_{2, T^3} d\mu_g & = \int_{t, u, S_{t, u}} - (t^{1 + s} \frac{c}{\kappa} \frac{v + c}{r} + \frac{\kappa}{c} \frac{v - c}{r}) LT^3 w \underline{L}T^3 w - \frac{\kappa}{c} \frac{2 v}{r} LT^3 \underline{w} \underline{L}T^3 \underline{w} d\mu_{\slashed{g}} du dt \\
     (\mathrm{take}\ |\cdot|) & \lesssim \int_{t, u, S_{t, u}} \frac{t^s}{\ln t} |L T^3 w| |\underline{L} T^3 w| + \frac{\ln t}{t^2} |L T^3 \underline{w}| |\underline{L} T^3 \underline{w|} d\mu_{\slashed{g}} du dt \\
     & \lesssim \int_{t, u, S_{t, u}} \frac{t^{- 1 + s}}{\ln t} |\underline{L} T^3 w|^2 d\mu_{\slashed{g}} du dt + \int_u \mathcal{F}_{T^3}^{[t^*, T]}(u) du + \int_t \frac{1}{t^2} \mathcal{E}_{T^3}^{[0, u^*]}(t) dt \\
     & \lesssim_{t^*, \mathrm{C}} 1 + \int_t \frac{1}{t^2} \mathcal{E}_{T^3}^{[0, u^*]}(t) dt + \int_u \mathcal{F}_{T^3}^{[t^*, T]}(u) du \\
\end{aligned}
\end{equation*}

\begin{equation*}
\begin{aligned}
    & \int_{[\delta, t^*] \times [0, u^*]} Q_{0, T^3 \psi} d\mu_g \\
    = & \int_{t, u, S_{t, u}} - (\underbrace{\sum\limits_{n_1 + n_2 = 2, n_i \geq 0} T^{n_1} [\mu \Box_g, T] T^{n_2} \psi}_{\mathrm{comm}} + \underbrace{T^3 (\mu \tilde{\rho}_\psi)}_{\mathrm{source}}) (a_\psi L T^3 \psi + b_\psi \underline{L} T^3 \psi) d\mu_\slashed{g} du dt \\
    = & I_{\mathrm{comm}, \psi} + I_{\mathrm{source}, \psi}
\end{aligned}
\end{equation*}

For $Q_{0, T^3} = Q_{0, T^3 w} + Q_{0, T^3 \underline{w}}$, we split them as in \ref{eq:local-base-error_split_1}, \ref{eq:local-base-error_split_2}.
\begin{equation*}
\begin{aligned}
& \begin{aligned}
    I_{\mathrm{comm}, \psi} & = \int_{t, u, S_{t, u}} - (a_\psi L T^3 \psi + b_\psi \underline{L} T^3 \psi) \sum\limits_{n_1 + n_2 = 2,\ n_i \geq 0} \mathrm{C}_{n_1, n_2} \\
    & \qquad (\underbrace{2 T^{n_1 + 1} (\frac{\kappa v}{c r}) L T^{n_2} \psi}_{1} + \underbrace{2 T^{n_1 + 1} (\frac{v + c}{r}) T^{n_2 + 1} \psi}_{2} \\
    & \qquad + \underbrace{T^{n_1 + 1} (\frac{\kappa}{c}) L L T^{n_2} \psi}_{3} + \underbrace{L T^{n_1 + 1} (\frac{\kappa}{c}) L T^{n_2} \psi}_{4}) d\mu_\slashed{g} du dt \\
    & = I_{\mathrm{comm}, 1, \psi} + I_{\mathrm{comm}, 2, \psi} + I_{\mathrm{comm}, 3, \psi} + I_{\mathrm{comm}, 4, \psi} \\
\end{aligned} \\
& \begin{aligned}
    I_{\mathrm{source}, \psi} = \int_{t, u, S_{t, u}} - T^3 (\mu \tilde{\rho}_\psi) (a_\psi L T^3 \psi + b_\psi \underline{L} T^3 \psi) d\mu_\slashed{g} du dt \\
\end{aligned} \\
\end{aligned}
\end{equation*}

First, we handle the source terms.

\begin{equation*}
\begin{aligned}
    T^3 (\mu \tilde{\rho}_{\underline{w}}) & = T^3 (2 \frac{\kappa v^2}{r^2} - 2 \frac{v + c}{r} T\underline{w} + L(\frac{\kappa v}{r}) + 2T(\frac{c v}{r})) \\
    & = T^3 (2 \frac{\kappa v^2}{r^2} + L(\frac{\kappa v}{r}) - 2 \frac{v + c}{r} T\underline{w} + 2\frac{c}{r} T v + 2 v T(\frac{c}{r})) \\
    & = T^3 (2 \frac{\kappa v^2}{r^2} + L(\frac{\kappa v}{r}) - 2 \frac{v}{r} T\underline{w} - 2\frac{c}{r} T w + 2 v T(\frac{c}{r})) \\
    (\mathrm{take}\ |\cdot|) & \lesssim_{t^*, \mathrm{C}} (\frac{1}{t^4} |T^3 \kappa| + \frac{\ln t}{t^4}) + (\frac{1}{t^2} |L T^3 \kappa| + \frac{1}{t^{\frac{5 + s}{2}}} |T^3 \kappa| + \frac{\ln t}{t} |L T^3 v| + \frac{\ln t}{t^{\frac{5 + s}{2}}}) \\
    & \qquad + (\frac{1}{t^2} |T^4 \underline{w}| + \frac{1}{t^3}) + (\frac{1}{t} |T^4 w| + \frac{\ln^2 t}{t^{\frac{7 + s}{2}}}) + (\frac{1}{t^2} |T^4 c| + \frac{1}{t^3} |T^3 \kappa| + \frac{\ln t}{t^3}) \\
    & \lesssim \frac{1}{t^{\frac{5 + s}{2}}} |T^3 \kappa| + \frac{1}{t} |T^4 w| + \frac{1}{t^2} |T^4 \underline{w}| + \frac{\ln t}{t} |LT^3 w| + \frac{\ln t}{t} |L T^3 \underline{w}| + \frac{\ln t}{t^{\frac{5 + s}{2}}} \\
    & \lesssim_{t^*, \mathrm{C}} \frac{1}{t^{\frac{5 + s}{2}}} |T^3 \kappa| + \frac{1}{t} |T^4 w| + \frac{1}{t^2} |T^4 \underline{w}| + \frac{\ln t}{t} |LT^3 w| + \frac{\ln t}{t^{\frac{5 + s}{2}}} \\
\end{aligned}
\end{equation*}

\begin{equation*}
\begin{aligned}
    I_{\mathrm{source}, \underline{w}} & = \int_{t, u, S_{t, u}} - T^3 (\mu \tilde{\rho}_{\underline{w}}) (\frac{\kappa}{c} L T^3 \underline{w} + \underline{L} T^3 \underline{w}) d\mu_\slashed{g} du dt \\
    (\mathrm{take}\ |\cdot|) & \lesssim_{t^*, \mathrm{C}} \int_{t, u, S_{t, u}} |\frac{\kappa}{c} L T^3 \underline{w} + \underline{L} T^3 \underline{w}| (\frac{1}{t^{\frac{5 + s}{2}}} |T^3 \kappa| + \frac{1}{t} |T^4 w| + \frac{1}{t^2} |T^4 \underline{w}| + \frac{\ln t}{t} |LT^3 w| + \frac{\ln t}{t^{\frac{5 + s}{2}}}) d\mu_\slashed{g} du dt \\
    & \lesssim \int_{t, u ,S_{t, u}} \frac{1}{t^{1 + s}} |\frac{\kappa}{c} L T^3 \underline{w} + \underline{L} T^3 \underline{w}|^2 d\mu_\slashed{g} du dt \\
    & \qquad + \int_{t, u ,S_{t, u}} (\frac{1}{t^4} |T^3 \kappa|^2 + \frac{1}{t^{1 - s}} |T^4 w|^2 + \frac{1}{t^{3 - s}} |T^4 \underline{w}|^2 + \frac{\ln^2 t}{t^{1 - s}} |LT^3 w|^2 + \frac{\ln^2 t}{t^4}) d\mu_\slashed{g} du dt \\
    & \lesssim_{t^*} 1 + \int_t \frac{1}{t^{1 + s}} \mathcal{E}_{T^3}^{[0, u^*]}(t) dt + \int_u \mathcal{F}_{T^3}^{[t^*, T]}(u) du
\end{aligned}
\end{equation*}

\begin{equation*}
\begin{aligned}
    T^3 (\mu \tilde{\rho}_w) & = T^3 (2 \frac{\kappa v^2}{r^2} + 2 \frac{v - c}{r} Tw + L(\frac{\kappa v}{r})) \\
    (\mathrm{take}\ |\cdot|) & \lesssim_{t^*, \mathrm{C}} \frac{1}{t^{\frac{5 + s}{2}}} |T^3 \kappa| + \frac{1}{t} |T^4 w| + \frac{1}{t^2} |T^4 \underline{w}| + \frac{\ln t}{t} |LT^3 w| + \frac{\ln t}{t} |L T^3 \underline{w}| + \frac{\ln t}{t^{\frac{5 + s}{2}}} \\
\end{aligned}
\end{equation*}

\begin{equation*}
\begin{aligned}
    I_{\mathrm{source}, w} & = \int_{t, u, S_{t, u}} - T^3 (\mu \tilde{\rho}_w) (t^{1 + s} \frac{c}{\kappa} L T^3 w + \underline{L} T^3 w) d\mu_\slashed{g} du dt \\
    (\mathrm{take}\ |\cdot|) & \lesssim_{t^*, \mathrm{C}} \int_{t, u, S_{t, u}} |t^{1 + s} \frac{c}{\kappa} L T^3 w + \underline{L} T^3 w| \\
    & \qquad (\frac{1}{t^{\frac{5 + s}{2}}} |T^3 \kappa| + \frac{1}{t} |T^4 w| + \frac{1}{t^2} |T^4 \underline{w}| + \frac{\ln t}{t} |LT^3 w| + \frac{\ln t}{t^{\frac{5 + s}{2}}}) d\mu_\slashed{g} du dt \\
    & \lesssim \int_{t, u, S_{t, u}} \frac{t^{1 + s}}{\ln t} |L T^3 w|^2 + \frac{1}{t^{1 + s}} |\underline{L} T^2 w|^2 d\mu_\slashed{g} du dt \\
    & \qquad + \int_{t, u ,S_{t, u}} (\frac{1}{t^4} |T^3 \kappa|^2 + \frac{1}{t^{1 - s}} |T^4 w|^2 + \frac{1}{t^{3 - s}} |T^4 \underline{w}|^2 + \frac{\ln^2 t}{t^{1 - s}} |LT^3 w|^2 + \frac{\ln^2 t}{t^4}) d\mu_\slashed{g} du dt \\
    & \lesssim_{t^*} 1 + \int_t \frac{1}{t^{1 + s}} \mathcal{E}_{T^3}^{[0, u^*]}(t) dt + \int_u \mathcal{F}_{T^3}^{[t^*, T]}(u) du
\end{aligned}
\end{equation*}

Then we handle the commutator terms. From how we deal with the source terms we can see that a general way is to use the $L^\infty$-estimates to bound the lower order terms, and use AM-GM inequality to separate the remaining highest order terms, obtaining $\mathcal{E}, \mathcal{F}$.

The integrand is $- (a_\psi L T^3 \psi + b_\psi \underline{L} T^3 \psi) \sum\limits_{n_1 + n_2 = 2,\ n_i \geq 0} T^{n_1} [\mu \Box_g, T] T^{n_2} \psi$, and we have the following for both $\psi = w, \underline{w}$.
\begin{equation*}
    \int_{t, u, S_{t, u}} \frac{1}{t^{1 + s}} |a_\psi L T^3 \psi + b_\psi \underline{L} T^3 \psi|^2 d\mu_{\slashed{g}} du dt \lesssim \int_t \frac{1}{t^{1 + s}} \mathcal{E}_{T^3}^{[0, u^*]}(t) dt + \int_u \mathcal{F}_{T^3}^{[t^*, T]}(u) du
\end{equation*}

The remaining part is $\int t^{1 + s} |\mathrm{terms\ in\ } \sum\limits_{n_1 + n_2 = 2,\ n_i \geq 0} T^{n_1} [\mu \Box_g, T] T^{n_2} \psi|^2$.

This would be the method we take to handle normal terms. However, there will be one exception that needs special treatment. We will see it later.
\begin{equation*}
\begin{aligned}
    & |\sum\limits_{n_1 + n_2 = 2,\ n_i \geq 0} \mathrm{C}_{n_1, n_2} T^{n_1 + 1} (\frac{\kappa v}{c r}) L T^{n_2} \psi| \lesssim_{t^*, \mathrm{C}} \frac{1}{t^{\frac{7 + s}{2}}} |T^3 \kappa| + \frac{\ln t}{t^{\frac{7 + s}{2}}} \\
    & |\sum\limits_{n_1 + n_2 = 2,\ n_i \geq 0} \mathrm{C}_{n_1, n_2} T^{n_1 + 1} (\frac{v + c}{r}) T^{n_2 + 1} \psi| \lesssim_{t^*, \mathrm{C}} \frac{\ln t}{t^3} \\
    & \begin{aligned}
        |\sum\limits_{n_1 + n_2 = 2,\ n_i \geq 0} \mathrm{C}_{n_1, n_2} L T^{n_1 + 1} (\frac{\kappa}{c}) L T^{n_2} \psi| & \lesssim \frac{1}{t^{\frac{3 + s}{2}}} |T^4 w| + \frac{1}{t^{\frac{3 + s}{2}}} |T^4 \underline{w}| + \frac{1}{t^{3 + s}} |T^3 \kappa| \\
        & \qquad + \frac{\ln t}{t^{\frac{3 + s}{2}}} |L T^3 w| + \frac{\ln t}{t^{\frac{3 + s}{2}}} |L T^3 \underline{w}| + \frac{1}{t^{\frac{5 + s}{2}}} \\
    \end{aligned}
\end{aligned}
\end{equation*}

These are all normal terms. So, we have finished the estimate of $I_{\mathrm{comm}, 1, \psi}, I_{\mathrm{comm}, 2, \psi}, I_{\mathrm{comm}, 4, \psi}$.
\begin{equation*}
    |\sum\limits_{\psi = w, \underline{w},\ i = 1, 2, 4} I_{\mathrm{comm}, i, \psi}| \lesssim_{t^*, \mathrm{C}} 1 + \int_t \frac{1}{t^{1 + s}} \mathcal{E}_{T^3}^{[0, u^*]}(t) dt + \int_u \mathcal{F}_{T^3}^{[t^*, T]}(u) du
\end{equation*}

Finally, we deal with $I_{\mathrm{comm}, 3, \psi}$.
\begin{equation*}
\begin{aligned}
    |\sum\limits_{n_1 + n_2 = 2,\ n_i \geq 0} \mathrm{C}_{n_1, n_2} T^{n_1 + 1} (\frac{\kappa}{c}) L L T^{n_2} \underline{w}| & = |\sum\limits_{n_1 + n_2 = 2,\ n_i \geq 0} - \mathrm{C}_{n_1, n_2} T^{n_1 + 1} (\frac{\kappa}{c}) L T^{n_2} (\frac{c v}{r})| \\
    & \lesssim_{t^*, \mathrm{C}} \frac{1}{t^{\frac{5 + s}{2}}} |T^3 \kappa| + \frac{\ln t}{t^{\frac{5 + s}{2}}} \\
\end{aligned}
\end{equation*}

\begin{equation*}
\begin{aligned}
    & |\sum\limits_{n_1 + n_2 = 2,\ n_i \geq 0} \mathrm{C}_{n_1, n_2} T^{n_1 + 1} (\frac{\kappa}{c}) L L T^{n_2} w| \\
    & = |\sum\limits_{n_1 + n_2 = 2,\ n_i \geq 0} \mathrm{C}_{n_1, n_2} T^{n_1 + 1} (\frac{\kappa}{c}) T^{n_2} (- 2 L (\frac{c}{\kappa}) T w - 2 \frac{c}{\kappa} L T w - L (\frac{c v}{r}))| \\
\end{aligned}
\end{equation*}

\begin{equation*}
    |\sum\limits_{n_1 + n_2 = 2,\ n_i \geq 0} \mathrm{C}_{n_1, n_2} T^{n_1 + 1} (\frac{\kappa}{c}) T^{n_2} (- 2 L (\frac{c}{\kappa}) T w - L (\frac{c v}{r}))| \lesssim_{t^*, \mathrm{C}} \frac{1}{t^{\frac{5 + s}{2}}} |T^3 \kappa| + \frac{\ln t}{t^{\frac{5 + s}{2}}}
\end{equation*}

\begin{equation*}
\begin{aligned}
    & |\sum\limits_{\psi = w, \underline{w}} I_{\mathrm{comm}, 3, \psi}| \\
    & \lesssim_{t^*, \mathrm{C}} |\int_{t, u, S_{t, u}} (t^{1 + s} \frac{c}{\kappa} L T^3 w + \underline{L} T^3 w) \sum\limits_{n_1 + n_2 = 2,\ n_i \geq 0} \mathrm{C}_{n_1, n_2} T^{n_1 + 1} (\frac{\kappa}{c}) T^{n_2} (\frac{c}{\kappa} L T w) d\mu_{\slashed{g}} du dt| \\
    & \qquad + 1 + \int_t \frac{1}{t^{1 + s}} \mathcal{E}_{T^3}^{[0, u^*]}(t) dt + \int_u \mathcal{F}_{T^3}^{[t^*, T]}(u) du \\
\end{aligned}
\end{equation*}

The first term is the special one that needs to be treated carefully. We need the flux to control it, just like how we control $a_w L T^3 w$. However, when $\kappa$ has the highest order derivative, extra treatment is required.

\begin{remark}
This is where we need the decay rate of $L T^n w$ to be $\frac{1}{t^{\frac{3}{2}+}}$. All the previous terms only require the decay rate to be $\frac{1}{t^{1 +}}$. One can check that there is room in the powers.
\end{remark}

\begin{equation*}
\begin{aligned}
    & \sum\limits_{n_1 + n_2 = 2,\ n_i \geq 0} \mathrm{C}_{n_1, n_2} T^{n_1 + 1} (\frac{\kappa}{c}) T^{n_2} (\frac{c}{\kappa} L T w) = \frac{L T w}{\kappa} T^3 \kappa + \mathrm{R} \\
    & |\mathrm{R}| \lesssim_{t^*, \mathrm{C}} |L T w| + |L T^2 w| + |L T^3 w|
\end{aligned}
\end{equation*}

\begin{equation*}
\begin{aligned}
    & \quad\  |\int_{t, u, S_{t, u}} (t^{1 + s} \frac{c}{\kappa} L T^3 w + \underline{L} T^3 w) \sum\limits_{n_1 + n_2 = 2,\ n_i \geq 0} \mathrm{C}_{n_1, n_2} T^{n_1 + 1} (\frac{\kappa}{c}) T^{n_2} (\frac{c}{\kappa} L T w) d\mu_{\slashed{g}} du dt| \\
    & \lesssim_{t^*, \mathrm{C}} |\int_{t, u, S_{t, u}} (t^{1 + s} \frac{c}{\kappa} L T^3 w + \underline{L} T^3 w) \frac{L T w}{\kappa} T^3 \kappa d\mu_{\slashed{g}} du dt| + \int_t \frac{\ln t}{t^{1 + s}} \mathcal{E}_{T^3}^{[0, u^*]}(t) dt + \int_u \mathcal{F}_{T^3}^{[t^*, T]}(u) du \\
    & \qquad + \int_{t, u, S_{t, u}} \frac{t^{1 + s}}{\ln t} (|L T w|^2 + |L T^2 w|^2 + |L T^3 w|^2) d\mu_{\slashed{g}} du dt
\end{aligned}
\end{equation*}

To control $L T w, L T^2 w$, we integrate from $C_0$.
\begin{equation*}
\begin{aligned}
    & \quad\  \int_{t, u, S_{t, u}} \frac{t^{1 + s}}{\ln t} |L T^2 w|^2 d\mu_{\slashed{g}} du dt \\
    & \lesssim \int_{t ,u} \frac{t^{3 + s}}{\ln t} |L T^2 w(t, 0)|^2 du dt + \int_{t ,u} \frac{t^{3 + s}}{\ln t} (\int_0^u |L T^3 w(t, u')| du')^2 du dt \\
    & \lesssim F_{\mathrm{Data}, 2} + \int_{t, u} \frac{t^{3 + s}}{\ln t} \int_0^u |L T^3 w(t, u')|^2 du' du dt \\
    & \lesssim F_{\mathrm{Data}, 2} + \int_{t, u'} \frac{t^{3 + s}}{\ln t} |L T^3 w(t, u')|^2 du' dt \\
    & \lesssim F_{\mathrm{Data}, 2} + \int_u \mathcal{F}_{T^3}^{[t^*, T]}(u) du
\end{aligned}
\end{equation*}

Similar results hold for $L T w$.
\begin{equation*}
    \int_{t, u, S_{t, u}} \frac{t^{1 + s}}{\ln t} |L T w|^2 d\mu_{\slashed{g}} du dt \lesssim F_{\mathrm{Data}, 1} + F_{\mathrm{Data}, 2} + \int_u \mathcal{F}_{T^3}^{[t^*, T]}(u) du
\end{equation*}

To control $\int \underline{L} T^3 w \frac{L T w}{\kappa} T^3 \kappa$, we use the improved estimate obtained from the first order Euler equation.
\begin{equation*}
\begin{aligned}
    & \quad\  |\int_{t, u, S_{t, u}} \underline{L} T^3 w \frac{L T w}{\kappa} T^3 \kappa d\mu_{\slashed{g}} du dt| \\
    & \lesssim_{\mathrm{C}} \int_{t, u, S_{t, u}} \frac{1}{\ln t} |\underline{L} T^3 w|^2 d\mu_{\slashed{g}} du dt + \int_{t, u, S_{t, u}} \frac{1}{t^{3 + s} \ln t} |T^3 \kappa|^2 d\mu_{\slashed{g}} du dt \\
    & \lesssim_{t^*, \mathrm{C}} 1 + \int_t \frac{1}{t^{1 + s}} \mathcal{E}_{T^3}^{[0, u^*]}(t) dt
\end{aligned}
\end{equation*}

Now we are left with the last term $\int t^{1 + s} \frac{c}{\kappa^2} L T^3 w L T w T^3 \kappa$. We will perform integration by parts to handle this.
\begin{equation*}
\begin{aligned}
& \begin{aligned}
    L L T^2 w & = T^2 (- 2 L (\frac{c}{\kappa}) T w - 2 \frac{c}{\kappa} L T w - L (\frac{c v}{r})) \\
    & = - 2 T^2 (\frac{c}{\kappa} L T w) + \mathrm{R} \\
    & = - 2 \frac{c}{\kappa} L T^3 w - 4 T (\frac{c}{\kappa}) L T^2 w - 2 T^2 (\frac{c}{\kappa}) L T w + \mathrm{R} \\
\end{aligned} \\
& |\mathrm{R}| \lesssim_{t^*, \mathrm{C}} \frac{1}{t^{\frac{5 + s}{2}}}
\end{aligned}
\end{equation*}

For the lower order terms, we integrate by parts in $u$.
\begin{equation*}
\begin{aligned}
    & \quad\  |\int_{t, u, S_{t, u}} t^{1 + s} \frac{c}{\kappa^2} L T^3 w L T w T^3 \kappa d\mu_{\slashed{g}} du dt| \\
    & = |\int_{t, u, S_{t, u}} t^{1 + s} \frac{1}{\kappa} L T w T^3 \kappa (- \frac{1}{2} L L T^2 w + 2 T(\frac{c}{\kappa}) L T^2 w + T^2(\frac{c}{\kappa}) L T w - \frac{1}{2} \mathrm{R}) d\mu_{\slashed{g}} du dt| \\
    & \quad\  |\int_{t, u, S_{t, u}} t^{1 + s} \frac{1}{\kappa} L T w T^3 \kappa \mathrm{R} d\mu_{\slashed{g}} du dt| \\
    & \lesssim_{t^*, \mathrm{C}} \int_{t, u, S_{t, u}} \frac{t^{1 + s}}{\ln t} |L T w|^2 + \frac{1}{t^4 \ln t} |T^3 \kappa|^2 \\
    & \lesssim_{t^*, \mathrm{C}} 1 + \int_t \frac{1}{t^2} \mathcal{E}_{T^3}^{[0, u^*]}(t) dt + \int_u \mathcal{F}_{T^3}^{[t^*, T]}(u) du \\
\end{aligned}
\end{equation*}

\begin{equation*}
\begin{aligned}
    & \quad\  \int_{t, u, S_{t, u}} t^{1 + s} \frac{1}{\kappa} L T w T^3 \kappa T(\frac{c}{\kappa}) L T^2 w d\mu_{\slashed{g}} du dt \\
    & = \int_{t, u} 4 \pi r^2 t^{1 + s} \frac{1}{\kappa} T(\frac{c}{\kappa}) T^3 \kappa L T w L T^2 w du dt \\
    & = \int_t 4 \pi r^2 t^{1 + s} \frac{1}{\kappa} T(\frac{c}{\kappa}) T^2 \kappa L T w L T^2 w du|_{u = 0}^{u = u^*} - \int_{t, u} 4 \pi T(r^2 t^{1 + s} \frac{1}{\kappa} T(\frac{c}{\kappa}) L T w L T^2 w) T^2 \kappa du dt \\
    (\mathrm{take}\ |\cdot|) & \lesssim_{t^*, \mathrm{C}} \int_{t, S_{t, u^*}} \frac{t^{1 + s}}{\ln t} |L T w| |L T^2 w| d\mu_{\slashed{g}} dt + \int_{t, S_{t, 0}} \frac{t^{1 + s}}{\ln t} |L T w| |L T^2 w| d\mu_{\slashed{g}} dt \\
    & \qquad \qquad + \int_{t, u, S_{t, u}} \frac{t^{1 + s}}{\ln t} (|L T w| + |L T^2 w|) (|L T^2 w| + |L T^3 w|) d\mu_{\slashed{g}} du dt \\
    & \lesssim F_{\mathrm{Data}, 1} + F_{\mathrm{Data}, 2} + \int_u \mathcal{F}_{T^3}^{[t^*, T]}(u) du \\
\end{aligned}
\end{equation*}

Using the same technique on $\int t^{1 + s} \frac{1}{\kappa} L T w T^3 \kappa T^2 (\frac{c}{\kappa}) L T w$ yields the same result. The only thing to note is that one needs to do a pairing for the $T^2 \kappa$ term in $T^2 (\frac{c}{\kappa})$.
\begin{equation*}
    \int_{t, u, S_{t, u}} t^{1 + s} \frac{1}{\kappa} L T w T^3 \kappa \frac{c T^2 \kappa}{\kappa^2} L T w d\mu_{\slashed{g}} du dt = \int_{t, u, S_{t, u}} t^{1 + s} \frac{1}{\kappa} L T w T(\frac{1}{2} T^2 \kappa)^2 \frac{c}{\kappa^2} L T w d\mu_{\slashed{g}} du dt
\end{equation*}

For the highest order term, we integrate by parts in $t$ first.
\begin{equation*}
\begin{aligned}
    & \quad\  \int_{t, u, S_{t, u}} t^{1 + s} \frac{1}{\kappa} L L T^2 w L T w T^3 \kappa d\mu_{\slashed{g}} du dt \\
    & = \int_{t, u} 4 \pi r^2 t^{1 + s} \frac{1}{\kappa} L L T^2 w L T w T^3 \kappa du dt \\
    & = \int_u 4 \pi r^2 t^{1 + s} \frac{1}{\kappa} L T^2 w L T w T^3 \kappa du|_{t = t^*}^{t = T} - \int_{t, u} 4 \pi L T^2 w L(r^2 t^{1 + s} \frac{1}{\kappa} L T w T^3 \kappa) du dt \\
\end{aligned}
\end{equation*}

\begin{equation*}
\begin{aligned}
    & \quad\  |\int_u 4 \pi r^2 t^{1 + s} \frac{1}{\kappa} L T^2 w L T w T^3 \kappa du|_{t = t^*}| \lesssim_{t^*, \mathrm{C}} 1 \\
    & \quad\  \int_u 4 \pi r^2 t^{1 + s} \frac{1}{\kappa} L T^2 w L T w T^3 \kappa du|_{t = T} \\
    & = 4 \pi r^2 t^{1 + s} \frac{1}{\kappa} L T^2 w L T w T^2 \kappa|_{t = T, u = 0}^{t = T, u = u^*} - \int_u 4 \pi T(r^2 t^{1 + s} \frac{1}{\kappa} L T^2 w L T w) T^2 \kappa du|_{t = T} \\
    (\mathrm{take}\ |\cdot|) & \lesssim_{t^*, \mathrm{C}} 1 + \int_{u, S_{T, u}} T^{1 + s} \frac{1}{T^{\frac{3 + s}{2}}} (|L T^3 w| + \frac{1}{T^{\frac{3 + s}{2}}}) d\mu_\slashed{g} du \\
    & \lesssim 1 + \epsilon \mathcal{E}_{T^3}^{[0, u^*]}(T) + \frac{1}{\epsilon}
\end{aligned}
\end{equation*}

We will take $\epsilon$ small enough, so this energy term will be absorbed by the left hand side.
\begin{equation*}
\begin{aligned}
    & \quad\  \int_{t, u} 4 \pi L T^2 w L(r^2 t^{1 + s} \frac{1}{\kappa} L T w T^3 \kappa) du dt \\
    & = \int_{t, u} 4 \pi L T^2 w r^2 t^{1 + s} \frac{1}{\kappa} L L T w T^3 \kappa du dt + \int_{t, u} 4 \pi L T^2 w r^2 t^{1 + s} \frac{1}{\kappa} L T w L T^3 \kappa du dt \\
    & \qquad + \int_{t, u} 4 \pi L T^2 w L(r^2 t^{1 + s} \frac{1}{\kappa}) L T w T^3 \kappa du dt \\
\end{aligned}
\end{equation*}

For $\int_{t, u} 4 \pi L T^2 w r^2 t^{1 + s} \frac{1}{\kappa} L L T w T^3 \kappa du dt$, we use the same method as above, writing $L L T w$ as $L T^n w + \mathrm{R}$.
\begin{equation*}
\begin{aligned}
& \begin{aligned}
    L L T w & = T (- 2 L (\frac{c}{\kappa}) T w - 2 \frac{c}{\kappa} L T w - L (\frac{c v}{r})) \\
    & = - 2 T (\frac{c}{\kappa} L T w) + \mathrm{R} \\
    & = - 2 \frac{c}{\kappa} L T^2 w - 2 T (\frac{c}{\kappa}) L T w + \mathrm{R} \\
\end{aligned} \\
& |\mathrm{R}| \lesssim_{t^*, \mathrm{C}} \frac{1}{t^{\frac{5 + s}{2}}}
\end{aligned}
\end{equation*}

Then we proceed as above to get the following.
\begin{equation*}
\begin{aligned}
    & \quad\  |\int_{t, u} 4 \pi L T^2 w r^2 t^{1 + s} \frac{1}{\kappa} L L T w T^3 \kappa du dt| \\
    & \lesssim_{t^*, \mathrm{C}} 1 + F_{\mathrm{Data}, 1} + F_{\mathrm{Data}, 2} + \int_t \frac{1}{t^2} \mathcal{E}_{T^3}^{[0, u^*]}(t) dt + \int_u \mathcal{F}_{T^3}^{[t^*, T]}(u) du \\
\end{aligned}
\end{equation*}

The other two parts are normal.
\begin{equation*}
\begin{aligned}
    |\int_{t, u} 4 \pi L T^2 w r^2 t^{1 + s} \frac{1}{\kappa} L T w L T^3 \kappa du dt| & \lesssim_{t^*, \mathrm{C}} \int_{t, u, S_{t, u}} \frac{1}{t^{\frac{1 - s}{2}} \ln t} |L T^2 w| |L T^3 \kappa| d\mu_{\slashed{g}} du dt \\
    & \lesssim_{t^*, \mathrm{C}} 1 + \int_t \frac{1}{t^2} \mathcal{E}_{T^3}^{[0, u^*]}(t) dt + \int_u \mathcal{F}_{T^3}^{[t^*, T]}(u) du \\
\end{aligned}
\end{equation*}

\begin{equation*}
\begin{aligned}
    |\int_{t, u} 4 \pi L T^2 w L(r^2 t^{1 + s} \frac{1}{\kappa}) L T w T^3 \kappa du dt| & \lesssim_{t^*, \mathrm{C}} \int_{t, u, S_{t, u}} \frac{1}{t^{\frac{3 - s}{2}} \ln t} |L T^2 w| |T^3 \kappa| d\mu_{\slashed{g}} du dt \\
    & \lesssim_{t^*, \mathrm{C}} 1 + \int_t \frac{1}{t^2} \mathcal{E}_{T^3}^{[0, u^*]}(t) dt + \int_u \mathcal{F}_{T^3}^{[t^*, T]}(u) du
\end{aligned}
\end{equation*}

We now collect all the estimates together to obtain the desired energy estimate. We absorb $F_{\mathrm{Data}, 1}, F_{\mathrm{Data}, 2}$ into $\mathrm{C}$ as well.
\begin{equation} \label{eq:global-base-energy_estimate}
    \mathcal{E}_{T^3}^{[0, u^*]}(T) + \mathcal{F}_{T^3}^{[t^*, T]}(u^*) \leq \mathrm{C}_{\mathrm{high}} (\int_{t^*}^T \frac{\ln t}{t^{1 + s}} \mathcal{E}^{[0, u^*]}(t) dt + \int_0^{u^*} \mathcal{F}_{T^3}^{[t^*, T]}(u) du) + \mathrm{C}_{\mathrm{low}}
\end{equation}

Then we place a smallness requirement on $u^*$ as in the local case.
\begin{equation} \label{eq:global-base-smallness_Gronwall}
    u^* \leq \frac{1}{2 \mathrm{C}_{\mathrm{high}}} = \frac{1}{\mathrm{C}_3^{(2)}}
\end{equation}

Then we look for $\hat{u}$ which maximizes the flux and get the following.
\begin{equation*}
    \max\limits_{u \in [0, u^*]} \mathcal{F}_{T^3}^{[t^*, T]}(u) = \mathcal{F}_{T^3}^{[t^*, T]}(\hat{u}) \leq 2 \mathrm{C}_{\mathrm{high}} \int_{t^*}^T \frac{\ln t}{t^{1 + s}} \mathcal{E}^{[0, u^*]}(t) dt + 2 \mathrm{C}_{\mathrm{low}}
\end{equation*}

Substituting this back into the energy estimate to get the following.
\begin{equation*}
    \mathcal{E}_{T^3}^{[0, u^*]}(T) \leq 2 \mathrm{C}_{\mathrm{high}} \int_{t^*}^T \frac{\ln t}{t^{1 + s}} \mathcal{E}^{[0, u^*]}(t) dt + 2 \mathrm{C}_{\mathrm{low}}
\end{equation*}

Then use the Gronwall inequality to obtain the uniform boundedness of energy and flux in $T$.
\begin{equation*}
\begin{cases}
    \max\limits_{t \in [t^*, T]} \mathcal{E}_{T^3}^{[0, u^*]}(t) \leq 2 \mathrm{C}_{\mathrm{low}} e^{\int_{t^*}^\infty \frac{2 \mathrm{C}_{\mathrm{high}} \ln \tau}{\tau^{1 + s}} d\tau} \\
    \max\limits_{u \in [0, u^*]} \mathcal{F}_{T^3}^{[t^*, T]}(u) \leq 2 \mathrm{C}_{\mathrm{low}} e^{\int_{t^*}^\infty \frac{2 \mathrm{C}_{\mathrm{high}} \ln \tau}{\tau^{1 + s}} d\tau} \\
\end{cases}
\end{equation*}

Finally, we integrate from $C_0$ to obtain the new $L^\infty$-bounds for the quantities in the bootstrap assumption \ref{eq:global-base-bootstrap_assumptions}.
\begin{equation*}
\begin{cases}
    |L T^2 w(t, u)| \leq \frac{\mathrm{C}' \sqrt{u^*} + M_2}{t^{\frac{3 + s}{2}}} \\
    |T^3 w(t, u)| \leq \frac{\mathrm{C}' \sqrt{u^*} + M_2}{t} \\
    |T^3 \underline{w}(t, u)| \leq \frac{\mathrm{C}' \sqrt{u^*} + M_2}{t} \\
\end{cases}
\end{equation*}

We place the last smallness requirement on $u^*$ to close the bootstrap argument.
\begin{equation} \label{eq:global-base-smallness_close}
    u^* \leq \frac{(\mathrm{C} - M_2)^2}{2 \mathrm{C}'^2}
\end{equation}

In the end, we list the extra $L^2$ and $L^\infty$-bounds we can obtain by using the first order Euler equations.
\begin{equation*}
\begin{aligned}
    \int_{u, S_{t, u}} |T^3 \kappa|^2 d\mu_{\slashed{g}} du & \lesssim t^2 \int_u |T^3 \kappa|^2 du \\
    & \lesssim t^2 \int_u |T^3 \kappa(t^*, u)|^2 + (\int_{t^*}^t |T^4 w| + |T^4 \underline{w}| d\tau)^2 du \\
    & \lesssim_{\mathrm{C}} t^2 + t^2 \int_{\tau, u} \frac{\tau^2}{\tau} (|T^4 w|^2 + |T^4 \underline{w}|^2) d\tau \int_{\tau} \frac{1}{\tau} d\tau \\
    & \lesssim_{t^*} t^2 + t^2 \ln^2 t \max\limits_{\tau \in [t^*, t]} \mathcal{E}_{T^3}^{[0, u^*]}(\tau) \\
    & \lesssim_{t^*, \mathrm{C}} t^2 \ln^2 t
\end{aligned}
\end{equation*}

\begin{equation*}
\begin{aligned}
    \int_{u, S_{t, u}} |T^4 w|^2 d\mu_{\slashed{g}} du & \lesssim_{t^*, \mathrm{C}} \int_{u, S_{t, u}} \frac{\ln^2 t}{t^{3 + s}} + \frac{1}{t^{3 + s}} |T^3 \kappa|^2 + \ln^2 t |L T^3 w|^2 d\mu_{\slashed{g}} du \\
    & \lesssim_{t^*, \mathrm{C}} \frac{\ln^2 t}{t^{1 + s}}
\end{aligned}
\end{equation*}

\begin{equation*}
    |L T^3 \underline{w}| \lesssim_{t^*, \mathrm{C}} \frac{1}{t^2}
\end{equation*}

This completes the proof.

\subsection{Proof of Proposition \ref{prop:global-regularity}}

The majority part of this proof is just a tweak of the proof of proposition \ref{prop:global-base_energy}, like the proof of proposition \ref{prop:local-regularity} compared to the proof of proposition \ref{prop:local-base_energy}, so we won't bother copying them again. We recall that we chose $N = 3$ as the starting point to ensure that all the terms in $N \geq 4$ contain at most one highest or second highest order term.

However, the problematic term $\sum\limits_{n_1 + n_2 = N - 1,\ n_i \geq 0} \mathrm{C}_{n_1, n_2} T^{n_1 + 1} (\frac{\kappa}{c}) T^{n_2} (\frac{c}{\kappa} L T w)$ which cost us many pages should still be treated seriously. We use a different technique now, which no longer involves integration by parts.

We observe that we can obtain a slightly stronger bound for the lower order flux.
\begin{equation*}
    \max\limits_{u \in [0, u^*]} \mathcal{F}_{T^{N - 1}}^{[t^*, T]}(u) \leq E_{N - 1}'
\end{equation*}

\begin{equation*}
\begin{aligned}
    & \quad\  \int_{t, S_{t, u}} \frac{t^{1 + s}}{\ln t} \max\limits_{u \in [0, u^*]} |L T^{N - 2} w(t, u)|^2 d\mu_{\slashed{g}} dt \\
    & \lesssim \int_{t, S_{t, u}} \frac{t^{1 + s}}{\ln t} (|L T^{N - 2} w(t, 0)|^2 + \int_0^{u^*} |L T^{N - 1} w(t, u')|^2 du') d\mu_{\slashed{g}} dt \\
    & \lesssim F_{\mathrm{Data}, N - 2} + \int_{t, u', S_{t, u'}} \frac{t^{1 + s}}{\ln t} |L T^{N - 1} w(t, u')|^2 d\mu_{\slashed{g}} du' dt \\
    & \lesssim F_{\mathrm{Data}, N - 2} + E_{N - 1}'
\end{aligned}
\end{equation*}

Similarly, for $L T^m w$, we have the following.
\begin{equation*}
    \int_{t, S_{t, u}} \frac{t^{1 + s}}{\ln t} \max\limits_{u \in [0, u^*]} |L T^m w(t, u)|^2 d\mu_{\slashed{g}} dt \lesssim F_{\mathrm{Data}, m} + E_{m + 1}' \lesssim_{\mathrm{lower}} 1
\end{equation*}

This means that we can first take the $L^\infty$ in $u$ and then integrate in $t$, and still have uniform boundedness. We apply this observation to the problematic term.
\begin{equation*}
\begin{aligned}
    & \quad\  \sum\limits_{n_1 + n_2 = N - 1,\ n_i \geq 0} \mathrm{C}_{n_1, n_2} T^{n_1 + 1} (\frac{\kappa}{c}) T^{n_2} (\frac{c}{\kappa} L T w) \\
    & = \mathrm{C}_{N - 1, 0} (\frac{T^N \kappa}{\kappa} L T w - N T^{N - 1} \kappa \frac{T c}{c \kappa} L T w - \frac{T^N c}{c} L T w) \\
    & \quad + \mathrm{C}_{N - 2, 1} \frac{T^{N - 1} \kappa}{c} T (\frac{c}{\kappa} L T w) + \mathrm{C}_{1, N - 2} T^2 (\frac{\kappa}{c}) \frac{c}{\kappa} L T^{N - 1} w \\
    & \quad + \mathrm{C}_{0, N - 1} (T (\frac{\kappa}{c}) \frac{c}{\kappa} L T^N w + N T (\frac{\kappa}{c}) T (\frac{c}{\kappa}) L T^{N - 1} w - T (\frac{\kappa}{c}) \frac{c T^{N - 1} \kappa}{\kappa^2} L T w) \\
    & \quad + \mathrm{R} \\
    & = \mathrm{H} + \mathrm{R} \\
    & \lesssim_{t^*, \mathrm{C}, N} \frac{|T^{N - 1} \kappa|}{\ln t} |L T^2 w| + (\frac{|T^N \kappa| + |T^{N - 1} \kappa|}{\ln t} + |T^N c|) |L T w| \\
    & \quad + |L T^N w| + |L T^{N - 1} w| + |R| \\
    & |\mathrm{R}| \lesssim_{\mathrm{lower}} |L T w| + \cdots + |L T^{N - 2} w|
\end{aligned}
\end{equation*}

\begin{equation*}
\begin{aligned}
    & \quad\  |\int_{t, u, S_{t, u}} (t^{1 + s} \frac{c}{\kappa} L T^N w + \underline{L} T^N w) \mathrm{R} d\mu_{\slashed{g}} du dt| \\
    & \lesssim \int_{t, u, S_{t, u}} \frac{t^{1 + s}}{\ln t} |L T^N w|^2 + \frac{\ln t}{t^{1 + s}} |\underline{L} T^N w|^2 + \frac{t^{1 + s}}{\ln t} |R|^2 d\mu_{\slashed{g}} du dt \\
    & \lesssim_{\mathrm{lower}} \int_t \frac{\ln t}{t^{1 + s}} \mathcal{E}_{T^N}^{[0, u^*]}(t) dt + \int_u \mathcal{F}_{T^N}^{[t^*, T]}(u) du + 1
\end{aligned}
\end{equation*}

\begin{equation*}
\begin{aligned}
    & \quad\  |\int_{t, u, S_{t, u}} (t^{1 + s} \frac{c}{\kappa} L T^N w + \underline{L} T^N w) \mathrm{H} d\mu_{\slashed{g}} du dt| \\
    & \lesssim \int_{t, u, S_{t, u}} \frac{t^{1 + s}}{\ln t} |L T^N w|^2 + \frac{\ln t}{t^{1 + s}} |\underline{L} T^N w|^2 + \frac{t^{1 + s}}{\ln t} |H|^2 d\mu_{\slashed{g}} du dt \\
    & \lesssim_{t^*, \mathrm{C}, N} \int_t \frac{\ln t}{t^{1 + s}} \mathcal{E}_{T^N}^{[0, u^*]}(t) dt + \int_u \mathcal{F}_{T^N}^{[t^*, T]}(u) du + \int_{t, u, S_{t, u}} \frac{t^{1 + s}}{\ln t} (|L T^N w|^2 + |L T^{N - 1} w|^2) d\mu_{\slashed{g}} du dt \\
    & \quad + \int_{t, u, S_{t, u}} \frac{t^{1 + s}}{\ln t} (\frac{|T^{N - 1} \kappa|^2}{\ln^2 t} |L T^2 w|^2 + (\frac{|T^N \kappa|^2 + |T^{N - 1} \kappa|^2}{\ln^2 t} + |T^N c|^2) |L T w|^2) d\mu_{\slashed{g}} du dt \\
    & \lesssim \int_t \frac{\ln t}{t^{1 + s}} \mathcal{E}_{T^N}^{[0, u^*]}(t) dt + \int_u \mathcal{F}_{T^N}^{[t^*, T]}(u) du + E_{N - 1}' \\
    & \quad + \int_{t, u, S_{t, u}} \frac{t^{1 + s}}{\ln t} (\frac{|T^{N - 1} \kappa|^2}{\ln^2 t} |L T^2 w|^2 + (\frac{|T^N \kappa|^2 + |T^{N - 1} \kappa|^2}{\ln^2 t} + |T^N c|^2) |L T w|^2) d\mu_{\slashed{g}} du dt \\
\end{aligned}
\end{equation*}

\begin{equation*}
\begin{aligned}
    & \quad\  \int_{t, u, S_{t, u}} \frac{t^{1 + s}}{\ln t} (\frac{|T^{N - 1} \kappa|^2}{\ln^2 t} |L T^2 w|^2 + (\frac{|T^N \kappa|^2 + |T^{N - 1} \kappa|^2}{\ln^2 t} + |T^N c|^2) |L T w|^2) d\mu_{\slashed{g}} du dt \\
    & \lesssim \int_t \frac{t^{1 + s}}{\ln t} \max\limits_{u \in [0, u^*]} |L T^2 w|^2 \int_{u, S_{t, u}} \frac{|T^{N - 1} \kappa|^2}{\ln^2 t} d\mu_{\slashed{g}} du dt \\
    & \quad + \int_t \frac{t^{1 + s}}{\ln t} \max\limits_{u \in [0, u^*]} |L T w|^2 \int_{u, S_{t, u}} \frac{|T^N \kappa|^2 + |T^{N - 1} \kappa|^2}{\ln^2 t} + |T^N c|^2 d\mu_{\slashed{g}} du dt \\
    & \lesssim E_{N - 1}' \int_t \frac{t^{3 + s}}{\ln t} \max\limits_{u \in [0, u^*]} (|L T^2 w|^2 + |L T w|^2) dt + \int_{t, u} \frac{t^{3 + s}}{\ln^3 t} \max\limits_{u \in [0, u^*]} |L T w|^2 |T^N \kappa|^2 du dt \\
\end{aligned}
\end{equation*}

\begin{equation*}
    E_{N - 1}' \int_t \frac{t^{3 + s}}{\ln t} \max\limits_{u \in [0, u^*]} (|L T^2 w|^2 + |L T w|^2) dt \lesssim_{\mathrm{lower}} 1
\end{equation*}

For the last term, we temporarily define a function $f$.
\begin{equation*}
    f(t) := \frac{t^{3 + s}}{\ln t} \max\limits_{u \in [0, u^*]} |L T w|^2
\end{equation*}

We list the properties of $f$, which come from the decay estimates and the boundedness of flux.
\begin{equation*}
\begin{cases}
    f(t) \lesssim_{t^*, \mathrm{C}} \frac{1}{\ln t} \\
    \int_{t^*}^\infty f(t) dt \lesssim_{t^*, \mathrm{C}} 1
\end{cases}
\end{equation*}

Then we prove that this term can be controlled by the integration of energy with an integrable coefficient.
\begin{equation*}
\begin{aligned}
    & \quad\  \int_{t, u} \frac{t^{3 + s}}{\ln^3 t} \max\limits_{u \in [0, u^*]} |L T w|^2 |T^N \kappa|^2 du dt \\
    & = \int_{t, u} \frac{f(t)}{\ln^2 t} |T^N \kappa|^2 du dt \\
    & \lesssim \int_{t, u} \frac{f(t)}{\ln^2 t} |T^N \kappa(t^*, u)|^2 du dt + \int_{t, u} \frac{f(t)}{\ln^2 t} (\int_{t^*}^t |T^{N + 1} w(\tau, u)| + |T^{N + 1} \underline{w}(\tau, u)| d\tau)^2 du dt \\
    & \lesssim_{t^*, \mathrm{C}} 1 + \int_{t, u} \frac{f(t)}{\ln^2 t} (\int_{t^*}^t |T^{N + 1} w(\tau, u)| + |T^{N + 1} \underline{w}(\tau, u)|)^2 du dt \\
\end{aligned}
\end{equation*}

\begin{equation*}
\begin{aligned}
    & \quad\  \int_{t, u} \frac{f(t)}{\ln^2 t} (\int_{t^*}^t |T^{N + 1} w(\tau, u)| + |T^{N + 1} \underline{w}(\tau, u)|)^2 du dt \\
    & \lesssim \int_{t, u} \frac{f(t)}{\ln^2 t} \int_{t^*}^t \tau^2 \frac{|T^{N + 1} w(\tau, u)|^2 + |T^{N + 1} \underline{w}(\tau, u)|^2}{\tau} d\tau \int_{t^*}^t \frac{1}{\tau} d\tau du dt \\
    & \lesssim_{t^*} \int_{t^*}^T \frac{f(t)}{\ln t} \int_0^{u^*} \int_{t^*}^t \tau^2 \frac{|T^{N + 1} w(\tau, u)|^2 + |T^{N + 1} \underline{w}(\tau, u)|^2}{\tau} d\tau du dt \\
    & = \int_{\tau, u} (\int_\tau^T \frac{f(t)}{\ln t} dt) \tau^2 \frac{|T^{N + 1} w(\tau, u)|^2 + |T^{N + 1} \underline{w}(\tau, u)|^2}{\tau} du d\tau \\
    & \lesssim \int_\tau (\frac{1}{\tau} \int_\tau^\infty \frac{f(t)}{\ln t} dt) \mathcal{E}_{T^N}^{[0, u^*]}(\tau) d\tau \\
    & =: \int_\tau \tilde{f}(\tau) \mathcal{E}_{T^N}^{[0, u^*]}(\tau) d\tau \\
\end{aligned}
\end{equation*}

Now we need to prove the coefficient is integrable.
\begin{equation*}
\begin{aligned}
    \int_{t^*}^\infty \tilde{f}(\tau) d\tau & = \int_{t^*}^\infty \frac{1}{\tau} \int_\tau^\infty \frac{f(t)}{\ln t} dt d\tau \\
    & = \ln \tau \int_\tau^\infty \frac{f(t)}{\ln t} dt|_{\tau = t^*}^{\tau = \infty} + \int_{t^*}^\infty \ln \tau \frac{f(\tau)}{\ln \tau} d\tau \\
    & = \lim\limits_{\tau \to \infty} \ln \tau \int_\tau^\infty \frac{f(t)}{\ln t} dt + \ln t^* \int_{t^*}^\infty \frac{f(t)}{\ln t} dt + \int_{t^*}^\infty f(\tau) d\tau \\
    & \leq \lim\limits_{\tau \to \infty} \frac{\ln \tau}{\ln \tau} \int_\tau^\infty f(t) dt + \frac{\ln t^*}{\ln t^*} \int_{t^*}^\infty f(t) dt + \int_{t^*}^\infty f(\tau) d\tau \\
    & = 2 \int_{t^*}^\infty f(\tau) d\tau \\
\end{aligned}
\end{equation*}

Using this result, we can obtain an energy estimate slightly different from \ref{eq:global-base-energy_estimate}.
\begin{equation*}
    \mathcal{E}_{T^N}^{[0, u^*]}(T) + \mathcal{F}_{T^N}^{[t^*, T]}(u^*) \leq \mathrm{C}_{\mathrm{high}, N} (\int_{t^*}^T (\tilde{f}(t) + \frac{\ln t}{t^{1 + s}}) \mathcal{E}^{[0, u^*]}(t) dt + \int_0^{u^*} \mathcal{F}_{T^3}^{[t^*, T]}(u) du)  + \mathrm{C}_{\mathrm{low}, N}
\end{equation*}

Where $\mathrm{C}_{\mathrm{high}, N}$ depends on $t^*, \mathrm{C}, N$, and $\mathrm{C}_{\mathrm{low}, N}$ depends on all the lower order bounds $E_n', M_n'$.

Like in the local case, we divide $[0, u^*]$ into smaller intervals and apply the energy estimate result to them inductively to finish the proof.

\subsection{Application to Background Close to Constant State}

In this subsection, we show that we can construct global in time rarefaction solutions for the background which is a perturbation of the constant state.

We still assume that the initial singularity is located at $t = 1, r_0 = 1$, and the constant state is $c_0 = 1, v_0 = 0$. We prove the following proposition.

\begin{proposition}
Suppose that the background solution is a perturbation of the constant state, which means that we require the data on $C_0$ to satisfy the following, with $\delta < \frac{1}{2}$ being a small constant and $n \geq 2$.
\begin{equation} \label{eq:global-perturbation_of_constant}
\begin{cases}
\begin{aligned}
    & |c - 1| \leq \frac{\epsilon}{t} && |v| \leq \frac{\epsilon}{t} && |r - t| \leq \epsilon \ln t \\
    & |L c| \leq \frac{\epsilon}{t^{2 - \delta}} && |L v| \leq \frac{\epsilon}{t^{2 - \delta}} && \\
    & |L^n c| \leq \frac{M_n}{t^{3 - 2 \delta}} && |L^n v| \leq \frac{M_n}{t^{3 - 2 \delta}} && \\
\end{aligned}
\end{cases}
\end{equation}

Then the data on $C_0$ solved from the ODE systems satisfy \ref{eq:global-C0_data}, for some $s > 0$ and different $\epsilon, M_n$. This means that we can construct a global in time rarefaction wave solution for this background by applying the propositions \ref{prop:global-base_energy} and \ref{prop:global-regularity}.
\end{proposition}

\begin{remark}
The decay rates we put on the background solution should be reasonable. $\epsilon$ stands for the initial status being close to the constant state, and $\frac{1}{t}$ decay is the decay rate we expect for 3-D wave equations, with the $L$-direction (tangential to the outgoing null cone) providing $\frac{1}{t^{1 - \delta}}$ extra decay.
\end{remark}

We prove this proposition by analyzing the ODE systems.

First, observe that the assumptions \ref{eq:global-perturbation_of_constant} already imply the desired bounds for $(c, v, r, L w, L \underline{w})$, where one need to use the Euler equation $L \underline{w} = - \frac{c v}{r}$ for $L \underline{w}$.

Now we begin to analyze the first ODE system \ref{eq:ODE-first}.
\begin{equation*}
\begin{cases}
    T w = - (\frac{v}{2 r} + \frac{L w}{2 c}) \kappa \\
    L T \underline{w} = \frac{\gamma - 1}{2} (- 2 \frac{\underline{w}}{r} T \underline{w} + 2 \frac{w}{r} T w - \frac{\underline{w}^2 - w^2}{r^2} \kappa) \\
    L \kappa = - \frac{\gamma + 1}{2} T \underline{w} + \frac{3 - \gamma}{2} T w
\end{cases}
\end{equation*}

We substitute $T w$ into the other two equations and extract the size of the coefficients, rewriting them as follows.
\begin{equation*}
\begin{cases}
\begin{aligned}
    L T \underline{w} & = - \frac{(\gamma - 1) \underline{w}}{r} T \underline{w} - \frac{\gamma - 1}{2} (\frac{w}{r} (\frac{v}{r} + \frac{L w}{c}) + \frac{\underline{w}^2 - w^2}{r^2}) \kappa \\
    & = - \frac{1}{t} (1 + O(\frac{\epsilon \ln t}{t})) T \underline{w} + O(\frac{\epsilon}{t^{3 - \delta}}) \kappa \\
\end{aligned} \\
\begin{aligned}
    L \kappa & = - \frac{\gamma + 1}{2} T \underline{w} - \frac{3 - \gamma}{2} (\frac{v}{2 r} + \frac{L w}{2 c}) \kappa \\
    & = - \frac{\gamma + 1}{2} T \underline{w} + O(\frac{\epsilon}{t^{2 - \delta}}) \kappa \\
\end{aligned}
\end{cases}
\end{equation*}

Here we are abusing the $O$ notation. We use it to represent the $\lesssim$ bounds. Then by integrating the equation we obtain the following.
\begin{equation*}
\begin{aligned}
    T \underline{w}(t) & = - \frac{2}{\gamma - 1} e^{\int_1^t - \frac{(\gamma - 1) \underline{w}(\tau)}{r(\tau)} d\tau} + \int_1^t O(\frac{\epsilon}{\tau^{3 - \delta}}) \kappa(\tau) e^{\int_\tau^t - \frac{(\gamma - 1) \underline{w}(\tau')}{r(\tau')} d\tau'} d\tau \\
    & = - \frac{2}{\gamma - 1} e^{- \ln t + O(\epsilon)} + \int_1^t O(\frac{\epsilon}{\tau^{3 - \delta}}) \kappa(\tau) e^{- \ln \frac{t}{\tau} + O(\epsilon)} d\tau \\
    & = - \frac{2}{\gamma - 1} \frac{1}{t} (1 + O(\epsilon)) + \frac{1}{t} \int_1^t O(\frac{\epsilon}{\tau^{2 - \delta}}) \kappa(\tau) d\tau \\
\end{aligned}
\end{equation*}

\begin{equation*}
\begin{aligned}
    \kappa & = \int_1^t - \frac{\gamma - 1}{2} T\underline{w}(\tau) e^{\int_\tau^t O(\frac{\epsilon}{\tau'^{2 - \delta}}) d\tau'} d\tau \\
    & = \int_1^t - \frac{\gamma - 1}{2} T\underline{w}(\tau) (1 + O(\epsilon)) d\tau \\
\end{aligned}
\end{equation*}

We combine these two formulas to obtain an equation for $T \underline{w}$.
\begin{equation*}
    T \underline{w}(t) = - \frac{2}{\gamma - 1} \frac{1}{t} (1 + O(\epsilon)) + \frac{1}{t} \int_1^t O(\frac{\epsilon}{\tau^{2 - \delta}}) \int_1^\tau - \frac{\gamma - 1}{2} T\underline{w}(\tau') (1 + O(\epsilon)) d\tau' d\tau
\end{equation*}

Now we can bootstrap. We assume $|T \underline{w}| \leq \frac{4}{\gamma - 1}$ and substitute into the equation twice to obtain the correct estimate for $T \underline{w}$.
\begin{equation*}
\begin{aligned}
    T \underline{w}(t) & = - \frac{2}{\gamma - 1} \frac{1}{t} (1 + O(\epsilon)) + \frac{1}{t} \int_1^t O(\frac{\epsilon}{\tau^{2 - \delta}}) \int_1^\tau - \frac{\gamma - 1}{2} T\underline{w}(\tau') (1 + O(\epsilon)) d\tau' d\tau \\
    & = - \frac{2}{\gamma - 1} \frac{1}{t} (1 + O(\epsilon)) + \frac{1}{t} \int_1^t O(\frac{\epsilon}{\tau^{1 - \delta}}) d\tau \\
    & = O(\frac{\epsilon}{t^{1 - \delta}}) \\
    (\mathrm{take}\ |\cdot|) & \leq \frac{4}{\gamma - 1}
\end{aligned}
\end{equation*}

\begin{equation*}
\begin{aligned}
    T \underline{w}(t) & = - \frac{2}{\gamma - 1} \frac{1}{t} (1 + O(\epsilon)) + \frac{1}{t} \int_1^t O(\frac{\epsilon}{\tau^{2 - \delta}}) \int_1^\tau - \frac{\gamma - 1}{2} T\underline{w}(\tau') (1 + O(\epsilon)) d\tau' d\tau \\
    & = - \frac{2}{\gamma - 1} \frac{1}{t} (1 + O(\epsilon)) + \frac{1}{t} \int_1^t O(\frac{\epsilon}{\tau^{2  - 2 \delta}}) d\tau \\
    & = - \frac{2}{\gamma - 1} \frac{1}{t} (1 + O(\epsilon)) + O(\frac{\epsilon}{t^{2 - 2 \delta}}) \\
    & = - \frac{2}{\gamma - 1} \frac{1}{t} (1 + O(\epsilon))
\end{aligned}
\end{equation*}

Then we use the formula for $\kappa$ to obtain the desired bound for $\kappa$.
\begin{equation*}
\begin{aligned}
    \kappa & = \int_1^t - \frac{\gamma - 1}{2} T\underline{w}(\tau) (1 + O(\epsilon)) d\tau \\
    & = \int_1^t - \frac{\gamma - 1}{2} (- \frac{2}{\gamma - 1} \frac{1}{\tau} (1 + O(\epsilon))) (1 + O(\epsilon)) d\tau \\
    & = \ln t (1 + O(\epsilon))
\end{aligned}
\end{equation*}

Then use the expression for $T w$.
\begin{equation*}
    T w = - (\frac{v}{2 r} + \frac{L w}{2 c}) \kappa = O(\frac{\epsilon \ln t}{t^{2 - \delta}})
\end{equation*}

We complete the estimate of the order one quantities $(T w, T \underline{w}, \kappa)$ by reading the size of their $L^n$-derivatives from the equations in the following order.
\begin{equation*}
    L \kappa \to L T \underline{w} \to L T w \to L^2 \kappa \to \cdots
\end{equation*}

\begin{equation*}
\begin{cases}
\begin{aligned}
    & L \kappa = \frac{1}{t} (1 + O(\epsilon)) && L T \underline{w} = O(\frac{1}{t^2}) && L T w = O(\frac{\ln t}{t^{3 - 2 \delta}}) \\
    & L^n \kappa = O(\frac{1}{t^2}) && L^n T \underline{w} = O(\frac{1}{t^2}) && L^n T w = O(\frac{\ln t}{t^{3 - 2 \delta}})
\end{aligned}
\end{cases}
\end{equation*}

We don't need the optimal decay for $L^n$-derivatives with $n \geq 2$. The decay rate listed above is enough for induction. Also note that $L, T$-derivatives always bring better (at least not worse) decay in $t$.

Then we move on to estimate quantities with more $T$-derivatives. We need the general ODE system \ref{eq:ODE-general}.
\begin{equation*}
\begin{cases}
    T^n w = - (\frac{v}{2 r} + \frac{Lw}{2 c}) T^{n - 1} \kappa + \sum\limits_{\substack{n_1 + n_2 = n - 1 \\ n_i \geq 0,\ n_2 \leq n - 2}} \mathrm{C}_{n_1, n_2} T^{n_1} (\frac{v}{r} + \frac{Lw}{c}) T^{n_2} \kappa \\
    L T^n \underline{w} = \frac{\gamma - 1}{2} (- 2 \frac{\underline{w}}{r} T^n \underline{w} + 2 \frac{w}{r} T^n w - \frac{\underline{w}^2 - w^2}{r^2} T^{n - 1} \kappa) \\
    \qquad \qquad + \sum\limits_{\substack{n_1 + n_2 = n - 1 \\ n_i \geq 0,\ n_2 \leq n - 2}} \mathrm{C}_{n_1, n_2} (- 2 T^{n_1} \frac{\underline{w}}{r} T^{n_2 + 1} \underline{w} + 2 T^{n_1} \frac{w}{r} T^{n_2 + 1} w - T^{n_1} \frac{\underline{w}^2 - w^2}{r^2} T^{n_2} \kappa) \\
    L T^{n - 1} \kappa = - \frac{\gamma + 1}{2} T^n \underline{w} + \frac{3 - \gamma}{2} T^n w
\end{cases}
\end{equation*}

Our induction assumption is that we already obtained the lower order estimates for $1 \leq m \leq n - 2$.
\begin{equation*}
\begin{cases}
\begin{aligned}
    & T^m \kappa = O(\ln t) && T^{m + 1} \underline{w} = O(\frac{1}{t}) && T^{m + 1} w = O(\frac{\ln t}{t^{2 - \delta}}) \\
    & L T^m \kappa = O(\frac{1}{t}) && L T^{m + 1} \underline{w} = O(\frac{1}{t^2}) && L T^{m + 1} w = O(\frac{\ln^{m + 1} t}{t^{3 - 2 \delta}}) \\
    & L^n T^m \kappa = O(\frac{1}{t^2}) && L^n T^{m + 1} \underline{w} = O(\frac{1}{t^2}) && L^n T^{m + 1} w = O(\frac{\ln^{m + 1} t}{t^{3 - 2 \delta}}) \\
\end{aligned}
\end{cases}
\end{equation*}

We use the same method as above. The difference is that we don't need to care about $\epsilon$, and we have inhomogeneous terms on the right hand side. We use the induction assumptions to obtain estimates for the inhomogeneous terms.
\begin{equation*}
\begin{aligned}
    f_n & := \sum\limits_{\substack{n_1 + n_2 = n - 1 \\ n_i \geq 0,\ n_2 \leq n - 2}} \mathrm{C}_{n_1, n_2} T^{n_1} (\frac{v}{r} + \frac{Lw}{c}) T^{n_2} \kappa = O(\frac{\ln t}{t^2}) \\
    g_n & := \sum\limits_{\substack{n_1 + n_2 = n - 1 \\ n_i \geq 0,\ n_2 \leq n - 2}}  \mathrm{C}_{n_1, n_2} (- 2 T^{n_1} \frac{\underline{w}}{r} T^{n_2 + 1} \underline{w} + 2 T^{n_1} \frac{w}{r} T^{n_2 + 1} w - T^{n_1} \frac{\underline{w}^2 - w^2}{r^2} T^{n_2} \kappa) = O(\frac{\ln t}{t^3})
\end{aligned}
\end{equation*}

Then we have the following.
\begin{equation*}
\begin{aligned}
    & L T^n \underline{w} = - \frac{1}{t} (1 + O(\frac{\epsilon \ln t}{t})) T^n \underline{w} + O(\frac{\epsilon}{t^{3 - \delta}}) T^{n - 1} \kappa + O(\frac{\ln t}{t^3}) \\
    & L T^{n - 1} \kappa = - \frac{\gamma + 1}{2} T^n \underline{w} + O(\frac{\epsilon}{t^{2 - \delta}}) T^{n - 1} \kappa + O(\frac{\ln t}{t^2}) \\
\end{aligned}
\end{equation*}

\begin{equation*}
\begin{aligned}
    & T^n \underline{w}(t) = \frac{1}{t} \int_1^t O(\frac{\epsilon}{\tau^{2 - \delta}}) T^{n - 1} \kappa(\tau) + O(\frac{\ln \tau}{\tau^2}) d\tau \\
    & T^{n - 1} \kappa = \int_1^t - \frac{\gamma - 1}{2} T^n \underline{w}(\tau) (1 + O(\epsilon)) + O(\frac{\ln \tau}{\tau^2}) d\tau \\
\end{aligned}
\end{equation*}

Then we obtain the equation for $T^n \underline{w}$, and the inhomogeneous terms only contribute an extra $\frac{1}{t}$ term.
\begin{equation*}
    T^n \underline{w}(t) = \frac{1}{t} (O(1) + \int_1^t O(\frac{\epsilon}{\tau^{2 - \delta}}) \int_1^t - \frac{\gamma - 1}{2} T^n \underline{w}(\tau') (1 + O(\epsilon)) d\tau' d\tau)
\end{equation*}

We can still begin from a loose bound and use the equation to improve it, same as above.
\begin{equation*}
    |T^n \underline{w}| \lesssim 1 \implies |T^n \underline{w}| \lesssim \frac{1}{t^{1 - \delta}} \implies |T^n \underline{w}| \lesssim \frac{1}{t}
\end{equation*}

The remaining parts are the same. We note that the growth of the power of logarithms in $L^n T^m w$ is caused by our choice to truncate the extra decay rate of high order $L$-derivatives. The following term causes the accumulation of logarithm powers, but logarithm powers won't break the $\frac{1}{t^{3 - 2 \delta}}$ decay, which is the main part.
\begin{equation*}
    L T^n w = - \frac{L L T^{n - 1} w}{2 c} \kappa + \cdots
\end{equation*}

This finishes the proof.

\printbibliography

\end{document}